\documentclass{article}

\usepackage{PRIMEarxiv}

\usepackage[utf8]{inputenc} % allow utf-8 input
\usepackage[T1]{fontenc}    % use 8-bit T1 fonts    
\usepackage{url}            % simple URL typesetting
\usepackage{booktabs}       % professional-quality tables
\RequirePackage{amsthm,amsmath,amsfonts,amssymb}     % blackboard math symbols
\usepackage{nicefrac}       % compact symbols for 1/2, etc.
\usepackage{microtype}      % microtypography
\usepackage{lipsum}
\usepackage{fancyhdr}       % header

\RequirePackage[numbers]{natbib}
\RequirePackage[colorlinks,citecolor=blue,urlcolor=blue]{hyperref}
\RequirePackage{graphicx}
\usepackage{subfigure}
\usepackage{caption}
\usepackage{subcaption}  % 用于创建子图
\usepackage{makecell}

\theoremstyle{plain}

\newtheorem{theorem}{Theorem}[section]
\newtheorem{lemma}[theorem]{Lemma}
\newtheorem{corollary}[theorem]{Corollary}
\newtheorem{proposition}[theorem]{Proposition}

\theoremstyle{remark}
\newtheorem{definition}[theorem]{Definition}

\newtheorem{assumption}{Assumption}
\newtheorem{remark}[theorem]{Remark}
\newcommand{\bbE}{\mathbb{E}}

\newcommand{\bbR}{\mathbb{R}}
\newcommand{\bbS}{\mathbb{S}}

\newcommand{\calB}{\mathcal{B}}
\newcommand{\calC}{\mathcal{C}}
\newcommand{\calD}{\mathcal{D}}

\newcommand{\calF}{\mathcal{F}}

\newcommand{\calH}{\mathcal{H}}

\newcommand{\calN}{\mathcal{N}}

\newcommand{\calP}{\mathcal{P}}

\newcommand{\calX}{\mathcal{X}}

\newcommand*{\gf}{\mathtt{GF}}
\newcommand*{\krr}{\mathtt{KRR}}

\newcommand{\colorblue}{}

\title{On the Pinsker bound of inner product kernel  regression in large dimensions}
\author{
  Weihao Lu \\
  Department of Statistics and Data Science \\
  Tsinghua University \\
  Beijing, China \\
  \texttt{luwh19@mails.tsinghua.edu.cn} \\
  \And
  Jialin Ding \\
  School of Mathematical Sciences \\
  Peking University \\
  Beijing, China \\
  \texttt{djl123456789@stu.pku.edu.cn} \\
  \And
  Haobo Zhang \\
  Department of Statistics and Data Science \\
  Tsinghua University \\
  Beijing, China \\
  \texttt{zhang-hb21@mails.tsinghua.edu.cn} \\
  \And
  Qian Lin\thanks{Corresponding author} \\
  Department of Statistics and Data Science \\
  Tsinghua University \\
  Beijing, China \\
  \texttt{qianlin@tsinghua.edu.cn} \\
}

\begin{document}
\maketitle

\begin{abstract}
Building on recent studies of large-dimensional kernel regression, particularly those involving inner product kernels on the sphere $\mathbb{S}^{d}$, we investigate the Pinsker bound for inner product kernel regression in such settings.
Specifically, we address the scenario where the sample size   $n$ is given by $\alpha d^{\gamma}(1+o_{d}(1))$ for some $\alpha, \gamma>0$.
We have determined the exact minimax risk for kernel regression in this setting, not only identifying the minimax rate but also the exact constant, known as the Pinsker constant, associated with the excess risk.
\end{abstract}

% keywords can be removed

\section{Introduction}\label{sec:intro}

For a fixed integer $m$ and  a non-decreasing sequence $\{a_j = (\pi j)^{2m}(1+o(1)), j=1,2,... \}$, Pinsker %\cite{pinsker1980optimal} 
considered the following Gaussian sequence model:
%In 1980, Pinsker's celebrated work (\cite{pinsker1980optimal}) established the first exact minimax risk for the problem of estimating the parameter of the Gaussian sequence model. 
$$z_j=\theta_j+ \varepsilon\xi_j, j=1,2, \cdots$$
where $\xi_j$ are i.i.d. $\calN\left(0, 1\right)$ and 
the sequence $\theta=\left(\theta_j\right)$ belongs to an ellipsoid 
$$\Theta_{R}=\left\{\theta: \sum\nolimits_{j} a_j \theta_j^2 \leq R\right\}.$$
In his celebrated work \citep{pinsker1980optimal}, he not only illustrated that the minimax rate of the risk $R(\widehat{\theta},\theta) := \bbE_{\theta}\|\widehat{\theta}-\theta\|^{2}_{\ell^{2}}$ is $\varepsilon^{\frac{4m}{2m+1}}$, but also demonstrated that
\begin{equation}\label{eqn:pinsker_original_result}
\begin{aligned}
    \inf_{\hat{\theta}} \sup _{\theta \in \Theta_{R}} \mathbb{E}_\theta \|\hat{\theta} - \theta\|_{\ell^2}^2
    &=
    \beta(m, R) \cdot \varepsilon^{\frac{4m}{2m+1}}(1+o(1)), 
    %\quad \text{ as } \varepsilon \to 0,
\end{aligned}
\end{equation}
where {\colorblue $\hat{\theta}$ is any estimator of $\theta$, measurable with respect to the observed data set $\{z_j\}_{j=1}^{\infty}$, } $\beta(m, R) = \left( \frac{m}{\pi(m+1)} \right)^{2m/(2m+1)} \left( R(2m+1) \right)^{1/(2m+1)}$. 
Later, Nussbaum \cite{nussbaum1985spline} 
considered the following nonparametric regression model:
$$
x_i = i/n, ~ y_i=f_{\star}(x_i)+\sigma\xi_i, \quad i \leq n,
$$
where $\xi_i$ are i.i.d. $\calN\left(0, 1\right)$ and the regression function $f_{\star}$ is in a subset of the Sobolev space
$W_2^m(R):=\left\{f \in L^2([0,1]) ;~\left\|D^m f\right\|^2 \leq R\right\}$. 
Interestingly,  Nussbaum   \cite{nussbaum1985spline} observed that the following exact asymptotic of the minimax risk for spline regression
\begin{equation}
\label{eqn:asymptotic_equiv_minimax_risk_nussbaum}
    \inf_{\hat{f}}\sup_{f_{\star} \in W_2^m(R)}\mathbb{E}_{f_{\star}}\|\hat{f} - f_{\star}\|_{L^2}^2 = \beta(m, R)\sigma^{\frac{4m}{2m+1}}n^{-\frac{2m}{2m+1}}(1+o(1)),
    %\quad \text{ as } n \to \infty.
\end{equation}
{\colorblue where  $\hat{f}$ is any estimator of $f_{\star}$, measurable with respect to the observed data set $\{(x_i, y_i)\}_{i=1}^{n}$. }
One can easily verify that the exact risk presented in Equation \eqref{eqn:pinsker_original_result} is equivalent to that in Equation \eqref{eqn:asymptotic_equiv_minimax_risk_nussbaum} when the noise level $\varepsilon$ is set to $\varepsilon = n^{-1/2}\sigma$, where $\sigma$ denotes the standard deviation of the noise.
This intriguing phenomenon, where the two asymptotics are equal, was rigorously justified by the seminal work on Le Cam  equivalence. These work established the asymptotic equivalence between Gaussian sequence models, the white noise model, and certain nonparametric regression models (see, e.g., \cite{brown1996asymptotic, le2012asymptotic, le2000asymptotics}).
Since then, subsequent studies have established similar exact risks for a variety of nonparametric estimation problems. These include density estimation, regression models with non-Gaussian noise or random designs, analysis of Besov bodies, and wavelet estimation (e.g., \cite{efroimovich1981estimation, golubev1994nonparametric, golubev1990risk, nussbaum1985spline, efromovich1996nonparametric, donoho1994minimax, donoho1990minimax, korostelev1993exact, tsybakov1997asymptotically}). For a detailed review of these developments, one can refer to \cite{nussbaum1999minimax} and the references therein.
Constants akin to $\beta(m, R)$, now often referred to as the Pinsker constant, play an indispensable role in studying the super-efficiency phenomenon observed in nonparametric problems. This phenomenon has been the subject of extensive investigation (e.g., \cite{le1953some, brown1997superefficiency, 
van1997superefficiency,
cai2005nonparametric}).

Recently, the strong theoretical links between the training dynamics within wide neural networks and the corresponding neural tangent kernel in regression have motivated substantial research into understanding the performance of spectral algorithms, such as kernel ridge regression and kernel gradient descent, in the context of kernel regression problems (see, e.g., \cite{Jacot_NTK_2018, Arora_on_2019, Du_gradient_2019_b, Du_gradient_2019_a, jianfa2022generalization, li2023statistical}). Modern approaches to kernel regression posit that the regression function \( f_{\star} \) is assumed to lie within the interpolation space \( [\mathcal{H}]^s \) of the Reproducing Kernel Hilbert Space (RKHS) \( \mathcal{H} \), where \( s \geq 0 \), rather than simply being an element of \( \mathcal{H} \). While kernel regression with a fixed data dimension \( d \) has been extensively studied, leading to insights on the minimax rate of the excess risk (\cite{Caponnetto2006OptimalRF, caponnetto2007optimal, raskutti2014early, Lin_Optimal_2020, zhang2023optimality}), the consistency of kernel interpolation (\cite{pmlr-v99-rakhlin19a, beaglehole2022kernel,
buchholz2022_KernelInterpolation, li2023kernel}), and the learning curves of spectral algorithms (\cite{Bordelon_Spectrum_2020, Cui2021GeneralizationER, jin2021learning, li2023asymptotic, li2024generalization}), there is an emerging interest in the performance of these algorithms when dealing with large-dimensional data (\cite{Liang_Just_2019, ghorbani2021linearized, mei2021learning, 
Ghosh_three_2021,
mei2022generalization, misiakiewicz_learning_2021, aerni2023strong, barzilai2023generalization, liang2020multiple, zhang2024phase, lu2023optimal, zhang2024optimal, lu2024saturation}). This shift in focus has been largely driven by the desire to better comprehend the intriguing phenomena observed in empirical studies of neural networks, such as double descent behavior and benign overfitting. 
{\colorblue
\cite{Karoui_spectrum_2010} studied the spectral properties of both inner-product and Euclidean distance kernels for general data distribution; based on this, \cite{bartlett2021deep} proved the polynomial barrier and asymptotic risk of kernel ridge regression (KRR) when $n \asymp d$; \cite{pandit2024universality} then proved the polynomial barrier and asymptotic risk of KRR when $n \asymp d^2$ for general data distribution; \cite{misiakiewicz2024non} proved the non-asymptotic deterministic equivalence of prediction risks for KRR; \cite{wang2023overparameterized, wang2024deformed, Montanari_interpolation_2021} proved the learning curves and polynomial approximation barrier of NTK regression for various data distributions.
}
Despite the growing interest in kernel regression, there remains a notable absence of Pinsker bounds for these problems, especially when the data dimensions are large.

Inspired by Pinsker's seminal work and the recent resurgence in kernel regression, we explore the Pinsker bound problem for kernel regression models that incorporate large-dimensional inner product kernels defined on the sphere $\mathbb S^{d}$. More precisely, 
%{\color{red} state your problem and results}
%Kernel regression problem considers $n$ i.i.d. samples $(x_{i}, y_{i})$ from a joint distribution supported on $\mathbb{R}^{d+1}\times \mathbb{R}$ and aims to find an estimator $\hat{f}_{n}$ based on a certain algorithm and given samples, such that the {\it excess risk}
%\begin{align*}
%\|\hat{f}_{n} - f_{\star}\|_{L^2}^2=\mathbb{E}_{x}(\hat{f}_{n}(x) - f_{\star}(x))^{2},
%\end{align*}
%$\|\hat{f}_{n} - f_{\star}\|_{L^2}^2$ is small, where $f_{\star}(x)=\mathbb{E}[y\vert x]$ is the {\it regression function} that falls into $[\calH]^{s}$, the interpolation space of the reproducing kernel Hilbert space (RKHS) $\calH$ with $s \geq 0$.
we address the scenario where the sample size $n$ is given by $\alpha d^{\gamma}(1+o_d(1))$ for some $\alpha, \gamma>0$.
We consider any RKHS $\calH$ associated with an inner product kernel, and we assume that the regression function falls into $\sqrt{R}[\calB]^{s}$, the ball in the interpolation space $[\calH]^{s}$ with radius $\sqrt{R}$.
Then, as stated in Theorem \ref{thm:main_pinsker_constant}, we establish the following exact minimax risk bound, known as the Pinsker bound:
\begin{equation*}
        \inf_{\hat{f}} \sup_{
        \rho_{f_{\star}} \in \mathcal{P}
        %f_{\star} \in \sqrt{R}[\calB]^{s}
        } \mathbb{E}_{(X, Y) \overset{\mathcal{D}}{\sim} \rho_{f_{\star}}^{\otimes n}}\left[
    \|\hat{f}-f_{\star}\|_{L^2}^2
    \right]
    =
    \calC^{\star} d^{-\zeta}(1+o_d(1)), 
    %\quad \text{ as } d \to \infty,
    \end{equation*}
where {\colorblue $\hat{f}$ is any estimator of $f_{\star}$, measurable with respect to the observed data set $(X, Y)$}, and $\mathcal{P}$ consist of all the distributions $\rho_{f_{\star}}$ on $\mathcal{X} \times \mathcal{Y}$ given by (\ref{equation:true_model}) such that Assumption 
\ref{assump_asymptotic}, 
\ref{assu:coef_of_inner_prod_kernel},  and \ref{assump_function_calss} hold for some $\alpha, \gamma>0$.

\subsection{Related works}

%In recent years, a special kind of nonparametric regression problem related to neural networks, the kernel regression problem, has received much attention.

%Inspired by the uniform convergence of excess risk between neural networks and kernel gradient flow, there has been a renewed interest in kernel regression problems and spectral algorithms with sample size $n = \Theta_d(d^{\gamma})$ for $\gamma>0$. 
%Recently, many new phenomena exhibited in large-dimensional (e.g., $n = \Theta_d(d^{\gamma})$ for $\gamma>0$.  ) kernel regression problems  have been found, and we review some of them as follows.

Recently, many new phenomena have been observed in large-dimensional kernel regression problems,  where the sample size $n$ is proportional to $d^{\gamma}$ for some $\gamma>0$. We review some of these phenomena as follows.

\paragraph*{Polynomial approximation barrier}
Early work on the polynomial approximation barrier phenomenon (e.g., \cite{ghorbani2021linearized, Donhauser_how_2021, mei2022generalization, xiao2022precise, misiakiewicz_spectrum_2022, hu2022sharp}) found that for any fixed square-integrable regression function, KRR and kernel gradient flow are consistent if and only if the regression function is a polynomial with degree $\leq \gamma$.  
Note that if $K$, the kernel function associated with $\calH$, is continuous, and if the eigenfunctions of $K$ form an orthonormal basis of $L^2$, then we have $[\calH]^{0} =  L^{2}$ (see, e.g., \cite{steinwart2012mercer, fischer2020_SobolevNorm}).
Hence, their results can also be interpreted in the following way: when $s=0$, and $\gamma$ is not an integer, the excess risks of spectral algorithms ( e.g.,  KRR, kernel gradient descent, etc.) lower bounded by some constants with high probability.
%Subsequent work (\cite{zhang2024optimal, lu2024saturation}) tried to explain the polynomial approximation barrier phenomenon by determining the convergence rate of the excess risk of spectral algorithms when $s \to 0$.
%They found that the convergence rate of the conditional excess risk of spectral algorithms drops when $\gamma$ exceeds any integer and remains constant for most of the other $\gamma$. 
We will provide a detailed discussion and comparison of these results in Section \ref{sec:comparison}.

\paragraph*{Optimal convergence rate for kernel regression}
Another line of work focused on the convergence rate of the minimax risk of kernel regression problems with any $s>0$ (\cite{lu2023optimal, zhang2024optimal, lu2024saturation}). 
Their results can be summarized as follows:
{\colorblue
\begin{itemize}
    \item Let $p=\lfloor \gamma/(s+1)\rfloor$. The minimax risk of kernel regression problems is bounded below by
    \begin{equation}\label{eqn:minimax_lower_early_our}
        \inf_{\hat{f}} \sup_{
        \rho_{f_{\star}} \in \mathcal{P}
        %f_{\star} \in \sqrt{R}[\calB]^{s}
        } \mathbb{E}_{(X, Y) \overset{\mathcal{D}}{\sim} \rho_{f_{\star}}^{\otimes n}}\left[
    \|\hat{f}-f_{\star}\|_{L^2}^2
    \right]
            =
            \left.
            \Omega_d\left(
        d^{-\zeta
        }
        \right) 
        \right/ 
        \text{poly}\left(\ln(d)\right),
    \end{equation}
    where $\hat{f}$ is any estimator of $f_{\star}$, measurable with respect to the observed data set $(X, Y)$, and $\zeta = \min\left\{
        \gamma-p, s(p+1)
        \right\}$ equals the one in Theorem \ref{thm:main_pinsker_constant}.

    \item If we fix a regression function $f_{\star}$ exactly falling into $[\calH]^{s}$, that is, we have $f_{\star} \in [\mathcal{H}]^{s}$ and $f_{\star} \notin [\mathcal{H}]^{s^{\prime}}$ for any $s^{\prime} > s$, then, there exists $t^{\star}>0$, such that for the estimator $\hat{f}_{t^{\star}}^{\gf}$ of kernel gradient flow and the estimator $\hat{f}_{t^{\star}}^{\krr}$ of kernel ridge regression, we have
  \begin{equation}\label{eqn:minimax_upper_early_our}
  \begin{aligned}
      \mathbb{E} \left( \left\|\hat{f}_{t^{\star}}^{\gf}  - f_{\star} \right\|^2_{L^2} \;\Big|\; X \right)  
             =&~ \Theta_{d, \mathbb{P}}\left(d^{-\zeta}\right)
             \cdot \text{poly}\left(\ln(d)\right)\\
        \mathbb{E} \left( \left\|\hat{f}_{t^{\star}}^{\krr}  - f_{\star} \right\|^2_{L^2} \;\Big|\; X \right)  
             =&~ \left\{\begin{matrix}
       & \Theta_{d, \mathbb{P}}(d^{-\zeta})
             \cdot \text{poly}\left(\ln(d)\right), \quad  s \leq 1;\\
        & \Theta_{d, \mathbb{P}}(d^{-\zeta^{\prime}})
             \cdot \text{poly}\left(\ln(d)\right), \quad  s > 1;\\
\end{matrix}\right.     
  \end{aligned}
    \end{equation}
    where $\tilde{s}=\min\{s, 2\}$, $\zeta^{\prime} = \min\left\{
        \gamma-p, \frac{\tau(\gamma-p+1)+p\tilde{s}}{\tau+1}, \tilde{s}(p+1)
        \right\}$, and $\Theta_{d, \mathbb{P}}$ is probability versions of the asymptotic notation $\Theta_{d}$.
\end{itemize}
}
The above results strongly suggest that the exact convergence rate of the minimax risk is \( d^{-\zeta} \), and this is one of the main foci of the current work.

%However, the above results are not enough to conclude that the convergence rate of the minimax risk is $d^{-\zeta}$, since (i) the upper bound holds for fixed regression function with high probability, and (ii) the additional logarithm terms.

\paragraph*{Periodic plateau behavior}

It has been observed that for any fixed function \( f_{\star} \in L^2 \), the excess risk experiences periodic reductions. This interesting phenomenon has been confirmed by the above results \cite{ghorbani2021linearized,
lu2023optimal, zhang2024optimal, lu2024saturation}. For instance, as shown in Fig.~\ref{figure_1_exact_rate}, when \( s = 3 \), the convergence rate of the excess risk remains constant for \( \gamma \) within intervals such as [3,4] and [7,8]. This phenomenon is referred to as the {\it periodic plateau} behavior of large-dimensional spectral algorithms. Based on this observation, it has been concluded that to improve the rate of excess risk for these spectral algorithms, it is necessary to increase the sample size beyond a certain threshold.

\subsection{Notations} 
We first introduce some absolute positive constants, and all other constants defined in the remainder of this paper only depend on these absolute positive constants.

\begin{definition}\label{def:abs_constants}
    We list all the absolute positive constants used in this paper:
    \begin{itemize}
    
    \item \(\alpha\), \(\gamma\), \(c_1\), \(c_2\): Constants in the asymptotic framework (\ref{Asym}). 

    \vspace{3pt}

    \item \(\sigma\): Upper bound on variance of the noise in (\ref{equation:true_model}).

    \vspace{3pt}
    
    \item {\colorblue \(K_{\max}\): maximum value of the kernel function in (\ref{eqn:maximum_of_kernel}).}

    \vspace{3pt}

    \item \(s\), \(R\): Constants representing the source condition and the upper bound on the norm of regression functions in the function class (\ref{eqn:function_calss}).

    \vspace{3pt}

    \item \(a_0, a_1, \ldots, a_{\lfloor \gamma \rfloor + 3}\) : The first $(\lfloor \gamma \rfloor+4)$ coefficients of the Taylor expansion of \(\Phi(\cdot)\) as specified in Assumption \ref{assu:coef_of_inner_prod_kernel}.

    \end{itemize}
\end{definition}

Let's denote the norm in $L^2:=L^2(\mathcal{X}, \rho_{\mathcal{X}})$ as $\|\cdot\|_{L^2}$.
For any integer $\ell \geq 0$, denote $\mathrm{P}_{>\ell}$ as the projection onto polynomials with degree $>\ell$.
We use asymptotic notations $O_d(\cdot),~o_d(\cdot),~\Omega_d(\cdot)$ and $\Theta_d(\cdot)$.
For instance, we say two (deterministic) quantities $U(d), V(d)$ satisfy $U(d) =  o_{d}(V(d)) $ if and only if for any $\varepsilon > 0$, there exists a constant $D_{\varepsilon}$ that only depends on $\varepsilon$ and the absolute positive constants listed in Definition \ref{def:abs_constants}, such that for any $d > D_{\varepsilon}$, we have $U(d)< \varepsilon V(d)$. Furthermore, we use the {\it asymptotically equivalence} notation $U(d) \sim V(d)$ if and only if we have $U(d) = V(d) (1+o_d(1))$.
We use 
%\(z \overset{\mathcal{D}}{=} \calN(0, \eta^2)\) to denote that \(z\) is a random variable following a normal distribution, and 
$z \overset{\mathcal{D}}{\sim} \rho$ to denote that $z$ follows the distribution $\rho$.

\section{Problem setting}\label{sec:prelimi}

We are interested in Pinsker's problem of kernel regression in the large-dimensional setting. To clarify any potential ambiguities and for future research purposes, we provide a detailed discussion of the problem settings in this section.

Suppose that we have observed $n$ i.i.d. samples $(x_i, y_i), i =1,2,\ldots, n$ from the model:
\begin{equation}\label{equation:true_model}
    y=f_{\star}(x)+\epsilon,
\end{equation}
where $x_i$'s are sampled from $\rho_{\calX}$,  which is the uniform distribution on $\mathcal{X} =\mathbb{S}^{d}\subset \bbR^{d+1}$,
$y \in \mathcal{Y} \subset \mathbb{R}$,
$f_{\star}$ is the regression function defined on $\mathcal{X}$, and  $\epsilon_1, \cdots, \epsilon_n \mid (x_1, \cdots, x_n)$ are mutually independent zero-mean variables with variances no greater than $\sigma^{2}$.
%Moreover, $\epsilon_i$ is independent of $x_i$ for $i \leq n$.
\iffalse
\begin{displaymath}
    \epsilon \overset{\mathcal{D}}{=} \mathcal{N}(0,\sigma^{2}).
  \end{displaymath} 
\fi
Denote the $n\times1$ data vector of ${y_i}$'s  and the $n\times d$ data matrix of $x_i$'s by $Y$ and $X$, respectively.
Moreover, let the sample size satisfy the following assumption:
\begin{assumption}\label{assump_asymptotic}
    We assume that there exist positive absolute constants $\alpha \in [c_1, c_2]$ and $\gamma>0$, such that the sample size satisfies
\begin{align}\label{Asym}
    n = \alpha d^{\gamma}(1+o_{d}(1)).
\end{align}
\end{assumption}

\subsection{Inner product kernels}
An inner product kernel $K$ defined on $\bbS^{d}$ is given by 
\begin{equation*}
%\label{eqn:def_inner_product_kernel}
    K(x, x^\prime) = \Phi(\left\langle x, x^\prime \right\rangle), \forall~ x, x^\prime \in \mathbb S^{d},
\end{equation*}
where  $\Phi:[-1,1] \to \mathbb{R}$ is a continuous function independent of $d$.
To avoid unnecessary notation, let us make the following assumption on the function $\Phi$.

\begin{assumption}\label{assu:coef_of_inner_prod_kernel} 
$\Phi(t) \in \mathcal{C}^{\infty} \left([-1,1]\right)$ is a fixed function independent of $d$ and there exists a non-negative sequence of absolute constants $\{a_j \geq 0\}_{j \geq 0}$ such that
    \begin{displaymath}
        \Phi(t) = \sum\nolimits_{j=0}^\infty a_j t^j, 
        %\quad a_{j} \geq 0, ~\text{for any}~ j \geq p+4,
    \end{displaymath}
    where 
    %$p$ is given in Definition \ref{def:abs_constants} and 
    $a_{j} > 0$ for any $j \leq \lfloor \gamma \rfloor+3$.
\end{assumption}

Assumption \ref{assu:coef_of_inner_prod_kernel} implies that the kernel function $K$ is bounded:
\begin{equation}\label{eqn:maximum_of_kernel}
    K_{\max}:= \sup_{x \in \mathcal X} K(x, x) \leq \sum\nolimits_{j=0}^{\infty} a_j < \infty.
\end{equation}
The purpose of assuming $\{a_0, \cdots, a_{\lfloor \gamma \rfloor+3}\}$ are positive is to maintain the clarity and simplicity of the main results and proofs.
%keep the main results and proofs clean. 
Note that, according to Theorem 1.b in \cite{gneiting2013strictly}, the inner product kernel $K$ on the sphere is positive-definite for all dimensions if and only if all coefficients $\{a_{j},j=0,1,2,...\}$ are non-negative.
Moreover, one can check that our main results, Theorem \ref{thm:main_pinsker_constant}, only depend on the former $\lfloor \gamma \rfloor+4$ coefficients $\{a_0, \cdots, a_{\lfloor \gamma \rfloor+3}\}$.
Therefore, the values of $\{a_{j} \geq 0\}_{j \geq \lfloor \gamma \rfloor+4}$ do not affect our results.
Furthermore, our main results can be extended when certain coefficients in $\{a_{j}\}_{j \geq 0}$ are zero. 
For example, one can consider the two-layer NTK defined as in \cite{Bietti_on_2019}, where $a_i=0$ for any $i=3,5,7, \cdots$).

Notice that the inner product kernel $K$ satisfying Assumption \ref{assu:coef_of_inner_prod_kernel} is positive-definite, hence
%For a positive definite kernel $K$ satisfying Assumption \ref{assu:coef_of_inner_prod_kernel}, 
the integral operator
$$T_{K}(f)(x)=\int K(x, x^{\prime}) f(x^{\prime}) ~\mathsf{d} \rho_{\calX}(x^{\prime})$$
 is a positive, self-adjoint, trace-class, and a compact operator (\cite{steinwart2012mercer}).  
The celebrated Mercer's theorem further assures that 
\begin{align}\label{eqn:mercer_decomp}
    K(x, x^{\prime})=\sum\nolimits_{j}\lambda_{j}\phi_{j}(x)\phi_{j}( x^{\prime}),
\end{align}
where the eigenvalues $\{\lambda_{j},j=1,2, \ldots\}$ form a non-increasing sequence, and the corresponding eigenfunctions of $\lambda_j$ is $\phi_{j}(\cdot)$, $j=1, 2, \ldots$.
%form an orthonormal basis in the $L^2(\calX, \rho_{\calX})$ function space.
%Calculating these eigenvalues is generally not straightforward. For example, \cite{jianfa2022generalization, li2023statistical} determined the eigenvalue decay rates $\lambda_j \asymp j^{-(d+1)/d}$ of neural-tangent kernels (NTK) on $\mathbb{R}^{d}$ with fixed dimension $d$ under certain conditions. Fortunately, 
Furthermore, since $K$ is an inner product kernel defined on the sphere, the Funk-Hecke formula provides a more concrete decomposition:
\begin{equation}\label{spherical_decomposition_of_inner_main}
\begin{aligned}
{K}(x,x^\prime) = \sum_{k=0}^{\infty} \mu_{k} \sum_{j=1}^{N(d, k)} Y_{k, j}(x) Y_{k, j}\left(x^\prime\right),
\end{aligned}
\end{equation}
where $Y_{k, j}$ for $j=1, \cdots, N(d, k)$ are spherical harmonic polynomials of degree $k$ and $\mu_{k}$'s are the eigenvalues of $K$  with multiplicity $N(d, k),k=0,1,\cdots$.  Here
$N(d,0)=1$; $N(d, k) = \frac{2k+d-1}{k} \cdot \frac{(k+d-2)!}{(d-1)!(k-1)!}, k =1,2,\cdots$.  We have to emphasize that $\mu_{k}$'s are not necessarily non-increasing.
For more details of the inner product kernels, readers can refer to \cite{gallier2009notes, ghorbani2021linearized}.

\iffalse
We will provide the correspondence between the (larger) eigenvalues \(\lambda_j\)'s and \(\mu_k\)'s in Lemma \ref{lemma:calcul_N}, along with the relationship between the eigenfunctions \(\phi_j\)'s corresponding to the eigenvalues \(\lambda_j\)'s and the spherical harmonics \(Y_{k,j}\)'s.
\fi

\begin{remark}\label{remark:sphere_data}
Most works analyzing spectral algorithms in large-dimensional settings focus on inner product kernels on spheres 
\citep[etc.]{ ghorbani2021linearized,  misiakiewicz_spectrum_2022, xiao2022precise, lu2023optimal, zhang2024optimal, lu2024saturation}.
On one hand, harmonic analysis on the sphere is clearer and more concise. For example, the properties of spherical harmonic polynomials are simpler than those of orthogonal series on general domains. This clarity makes Mercer’s decomposition of the inner product more explicit, avoiding several abstract assumptions (e.g., \cite{https://doi.org/10.1002/cpa.22008}).
On the other hand, very few results are available for Mercer’s decomposition of kernels on general domains, especially when considering the domain’s dimension. Although some studies have attempted to relax the spherical assumption (e.g., \cite{liang2020multiple, aerni2023strong, barzilai2023generalization}), most of them either (i) adopt a near-spherical assumption, (ii) impose strong assumptions on the regression function (e.g., $f_{\star}(x) = x[1] x[2] \cdots x[L]$ for an integer $L>0$, where $x[i]$ denote the $i$-th component of $x$), or (iii) cannot determine the convergence rate of the spectral algorithm's excess risk.
\end{remark}

\iffalse
We consider the inner product kernels on the sphere mainly because harmonic analysis on the sphere is clearer and more concise ( e.g., the properties of spherical harmonic polynomials are simpler than those of orthogonal series on general domains). This clarity makes Mercer’s decomposition of the inner product more explicit, avoiding several abstract assumptions ( e.g., \cite{https://doi.org/10.1002/cpa.22008}).
    Additionally, very few results are available for Mercer’s decomposition of kernels on general domains, especially when considering the domain's dimension.
%For example, even the eigen-decay rate of neural tangent kernels is only determined for the spheres. 
Due to these technical reasons, most works analyzing spectral algorithms in large-dimensional settings focus on inner product kernels on spheres \citep[etc.]{liang2020multiple, ghorbani2021linearized,  misiakiewicz_spectrum_2022, xiao2022precise, lu2023optimal}. 
\fi

\subsection{The interpolation space}\label{subsec:interpolation_space}

The interpolation space $[\mathcal{H}]^s$ (associated with the inner product kernel $K$) with source condition $s \geq 0$ is defined as
\begin{equation*}
  [\mathcal{H}]^s := 
  \Big\{
  \sum\nolimits_{j=1}^{\infty} b_j \lambda_j^{s / 2}\phi_{j}(\cdot): \left(b_j\right)_{j} \in \ell_2 
  \Big\} 
  \subseteq L^{2}(\mathcal{X}, \rho_{\mathcal{X}}),
\end{equation*}
with $\lambda_j$'s and $\phi_j(\cdot)$'s defined in (\ref{eqn:mercer_decomp}), and the inner product deduced from
\begin{equation*}\begin{aligned}
\Big\|\sum\nolimits_{j=1}^{\infty} b_j \lambda_j^{s / 2} \phi_j
\Big\|_{[\mathcal{H}]^s}=
\Big(\sum\nolimits_{j=1}^{\infty} b_j^2
\Big)^{1 / 2}.
\end{aligned}
\end{equation*}

It is easy to show that $[\mathcal{H}]^s $ is also a separable Hilbert space with orthonormal basis $ \{\lambda_{j}^{s/2} \phi_{j}\}_{j}$. 
Generally speaking, functions in $[\mathcal{H}]^s$ become smoother as $s$ increases (see, e.g., the example of Sobolev spaces in \cite{edmunds1996function, zhang2023optimality}). 
The two most interesting interpolation spaces are $[\calH]^{0} \subseteq L^{2}$ and $[\calH]^{1} = \calH$. 
%Moreover, since $K$ is an inner product kernel defined on $\mathbb{S}^{d}$, we have $[\calH]^{0} = L^{2}$ if and only if $a_j >0$ for any $j=0, 1, \cdots$, where the absolute constants $\{a_j\}_{j \geq 0}$ are defined in Assumption \ref{assu:coef_of_inner_prod_kernel} (see, e.g., \cite{steinwart2012mercer, fischer2020_SobolevNorm}).

In kernel regression studies, it is typically assumed that $f_{\star}$ falls into the RKHS $\calH$ (e.g., \cite{Caponnetto2006OptimalRF, caponnetto2007optimal, Yao2007OnES, raskutti2014early, Liang_Just_2019}). 
However, subsequent research has suggested that the RKHS $\cal H$ might be too restrictive, prompting interest in the performance of kernel regression in the misspecified case with $s \in (0, 1)$ (\cite{fischer2020_SobolevNorm, zhang2023optimality, zhang2023optimality_2, zhang2024optimal}). 
Recently, several studies on large-dimensional kernel regression have considered the extreme case where $s=0$ (e.g., \cite{ghorbani2021linearized, mei2021learning, mei2022generalization, misiakiewicz_learning_2021}).
To fully capture the performance of large-dimensional kernel regression and provide a unified explanation for previous work, we assume that the regression function falls into the ball in $[\calH]^{s}$ with radius $\sqrt{R}$~:
%, where $s$ and $R$ are positive absolute constants given in Definition \ref{def:abs_constants}.
\begin{assumption}\label{assump_function_calss}
There exist two positive absolute constants $s$ and $R$, such that we have
    \begin{equation}\label{eqn:function_calss}
    f_{\star} \in \sqrt{R}[\calB]^{s} := \left\{f \in [\calH]^{s} \mid \|f\|_{[\calH]^{s}} \leq \sqrt{R}\right\}.
\end{equation}
\end{assumption}

\section{Main Results} %: Pinsker constant of interpolation spaces in large dimensions}
\label{sec:main_result}

We present our main results, demonstrating that the minimax rate of the excess risk for the function class $\sqrt{R}[\calB]^{s}$ is asymptotically equivalent to the {\it Pinsker constant} $\calC^{\star}$ times a corresponding convergence rate $d^{-\zeta}$.

\begin{theorem}\label{thm:main_pinsker_constant}
    Let $\mathcal{P}$ consist of all the distributions $\rho_{f_{\star}}$ on $\mathcal{X} \times \mathcal{Y}$ given by (\ref{equation:true_model}) such that Assumption 
\ref{assump_asymptotic}, 
\ref{assu:coef_of_inner_prod_kernel},  and \ref{assump_function_calss} hold for some $\alpha, \gamma>0$. Then, 
    when $d \geq \mathfrak{C}$ (a sufficiently large constant only depending on the absolute constants given in Definition \ref{def:abs_constants}), we have
    \begin{equation*}
        \inf_{\hat{f}} \sup_{
        \rho_{f_{\star}} \in \mathcal{P}
        %f_{\star} \in \sqrt{R}[\calB]^{s}
        } \mathbb{E}_{(X, Y) \overset{\mathcal{D}}{\sim} \rho_{f_{\star}}^{\otimes n}}\left[
    \|\hat{f}-f_{\star}\|_{L^2}^2
    \right]
    =
    \calC^{\star} d^{-\zeta}(1+o_{d}(1)),
    \end{equation*}
    {\colorblue where $\hat{f}$ is any estimator of $f_{\star}$, measurable with respect to the observed data set $(X, Y)$.}
    Further, define \(p := \left\lfloor \frac{\gamma}{s+1} \right\rfloor\), then we have:
\begin{itemize}
    \item[(i)] When $p(s+1) \leq \gamma < p(s+1)+s$, we have $\zeta = \gamma - p$, and
    $$
    \calC^{\star} :=
    \frac{ \sigma^2   }{\alpha  p!  + \sigma^2 / (R a_p^{s} (p!)^{s})\mathbf{1}\{\gamma=p(s+1)\}}
    $$
    \iffalse
    $$
    \calC^{\star} :=
    \begin{cases}
       \frac{ \sigma^2 }{\alpha  p!  + \sigma^2 / (R a_p^{s} (p!)^{s})}  & \quad \text{if }  \gamma = p(s+1)  \\
       \frac{ \sigma^2 }{\alpha  p! }  & \quad \text{if } \gamma > p(s+1)
    \end{cases}
    $$
    \fi

    \item[(ii)] When $p(s+1)+s \leq \gamma < (p+1)(s+1)$, we have $\zeta = (p+1)s$, and 
    $$
    \calC^{\star} :=
        R  {a_{p+1}^{s}} ((p+1)!)^{s} 
        +
        \frac{ \sigma^2 }{\alpha  p! } \mathbf{1}\{\gamma = p(s+1)+s\}.
    $$
    \iffalse
    $$
    \calC^{\star} :=
    \begin{cases}
        R  {a_{p+1}^{s}} ((p+1)!)^{s} 
        +
        \frac{ \sigma^2 }{\alpha  p! } & \quad \text{if }  \gamma = p(s+1)+s  \\
        R  {a_{p+1}^{s}} ((p+1)!)^{s} & \quad \text{if } \gamma > p(s+1)+s.
    \end{cases}
    $$
    \fi

    \iffalse
    \item[(i)] When $\gamma = p(s+1)$, we have
    $$
    \calC^{\star} := 
        \frac{ \sigma^2 }{\alpha  p!  + \sigma^2 / (R a_p^{s} (p!)^{s})}
        \quad \text{ and } \quad
        \zeta = ps;
    $$
    \item[(ii)] When $p(s+1)<\gamma < p(s+1)+s$, we have
    $$
    \calC^{\star} := 
        \frac{ \sigma^2 }{\alpha  p! }
        \quad \text{ and } \quad
        \zeta = \gamma - p;
    $$
    \item[(iii)] When $\gamma = p(s+1)+s$, we have
    $$
    \calC^{\star} := 
        R  {a_{p+1}^{s}} ((p+1)!)^{s} 
        +
        \frac{ \sigma^2 }{\alpha  p! }
        \quad \text{ and } \quad
        \zeta = (p+1)s;
    $$
    \item[(iv)] When $p(s+1)+s < \gamma < (p+1)(s+1)$, we have
    $$
    \calC^{\star} := 
        R  {a_{p+1}^{s}} ((p+1)!)^{s} 
        \quad \text{ and } \quad
        \zeta = (p+1)s.
    $$
    \fi
\end{itemize}
\end{theorem}

The proof of Theorem \ref{thm:main_pinsker_constant} is organized as follows: 
In Section \ref{sec:pinsker_quantity}, we define a quantity $\calD^{\star}$ and demonstrate that $\calD^{\star} \sim \calC^{\star} d^{-\zeta}$.
In Section \ref{sec:upper_sketch}, we provide a sketch showing that the minimax excess risk in Theorem \ref{thm:main_pinsker_constant} has an upper bound $\calD^{\star}(1+o_d(1))$, and we defer the full proof to Appendix \ref{sec:upper}.
Finally, the proof for the corresponding lower bound in Theorem \ref{thm:main_pinsker_constant}, being relatively routine, is deferred to Appendix \ref{sec:lower}.

\iffalse
\begin{remark}
   Theorem \ref{thm:main_pinsker_constant} assures that as long as \( d \geq \mathfrak{C} \), where \( \mathfrak{C} \) is a sufficiently large constant depending only on the absolute constants \( \gamma, a_0, \ldots, a_{\lfloor \gamma \rfloor + 3} \), the Pinsker bound remains valid for kernel regression. %We note that \( \mathfrak{C} \) does not depend on \( s \). Thus, for any sufficiently large \( d \), the function class \( \sqrt{R}[\mathcal{B}]^{s} \) is sufficiently large if \( s \) is sufficiently small. 
\end{remark}
\fi

%        For fixed $s>0$ and $R>0$, there may be concerns that the function class $\sqrt{R}[\calB]^{s}$ becomes too small as \(d\) increases.        However,        In other words, for large-dimensional data where $d$ exceeds a certain threshold independent of $s$, we can still obtain the corresponding Pinsker bound even for kernel regression problems with large function spaces (i.e., with small $s$).

%Theorem \ref{thm:main_pinsker_constant}  provides the exact asymptotic behavior of the minimax risk, i. e. not only the "optimal rate of convergence for estimators" $d^{-\zeta}$, but also the "optimal constant" $\calC^{\star}$. To help readers understand Theorem \ref{thm:main_pinsker_constant} better, we provide some interpretations of the results in Theorem \ref{thm:main_pinsker_constant} as the following two parts.

Theorem \ref{thm:main_pinsker_constant} delineates the precise asymptotic behavior of the minimax risk.
It specifies not only the optimal convergence rate $d^{-\zeta}$
  for estimation but also the optimal constant $\calC^{\star}$. To enhance readers' comprehension of Theorem \ref{thm:main_pinsker_constant}, we offer interpretations of its results in the following two parts.

%Moreover, in Section \ref{sec:gauss_white_noise}, we will generalize our results to the Gaussian sequence model with parameter space defined by $\mu_k$'s, and compare Theorem \ref{thm:main_pinsker_constant} with the existing result on Sobolev RKHS given by \cite{nussbaum1999minimax}.

\iffalse
\paragraph*{New phenomena occurred in large-dimensional spectral algorithms}
Notice that Theorem \ref{thm:main_pinsker_constant} has successfully reproduced the previously established minimax rate $d^{-\zeta}$ of the interpolation space (\cite{lu2023optimal, zhang2024optimal, lu2024saturation}).
More importantly, these minimax rates match the rates on the excess risk of spectral algorithms with qualification $\tau \geq s$ (see, e.g., Theorem 3 in \cite{zhang2024optimal} and Theorem 4.1 in \cite{lu2024saturation}).
The above works, together with Theorem \ref{thm:main_pinsker_constant}, imply that several new phenomena occur in many spectral algorithms in large dimensions.
For example, we can find the periodic plateau behavior
that $\calC^{\star}$ and $\zeta$ have constant values when $p(s+1)+s < \gamma < (p+1)(s+1)$.
\fi

\paragraph*{Exact convergence rate of the minimax risk}
%Several recent studies (\cite{lu2023optimal, zhang2024optimal, lu2024saturation}) have discussed the exact convergence rate of the minimax risk for kernel regression in high dimensional settings. These studies conjectured that the correct rate is exactly $d^{-\zeta}$. Theorem 3.1 rigorously confirms this conjecture.
Several recent studies (\cite{lu2023optimal, zhang2024optimal, lu2024saturation}) have obtained nearly exact convergence rates, i.e., up to some logarithmic term, of the minimax risk for kernel regression in large-dimensional settings. These studies suggested that the correct rate is \(d^{-\zeta}\). Theorem \ref{thm:main_pinsker_constant} rigorously confirms this conjecture.

%The exact convergence rate of the minimax risk of kernel regression problems has been discussed in several recent studies (\cite{lu2023optimal, zhang2024optimal, lu2024saturation}). However, as indicated in (\ref{eqn:minimax_lower_early_our}) and (\ref{eqn:minimax_upper_early_our}), these studies could only conjecture that the correct rate is exactly $d^{-\zeta}$. Thanks to Theorem \ref{thm:main_pinsker_constant}, now we can rigorously confirm this conjecture.\\

\begin{figure}[ht]
% \vskip 0.05in
\centering
\subfigure[]{\includegraphics[width=0.48\columnwidth]{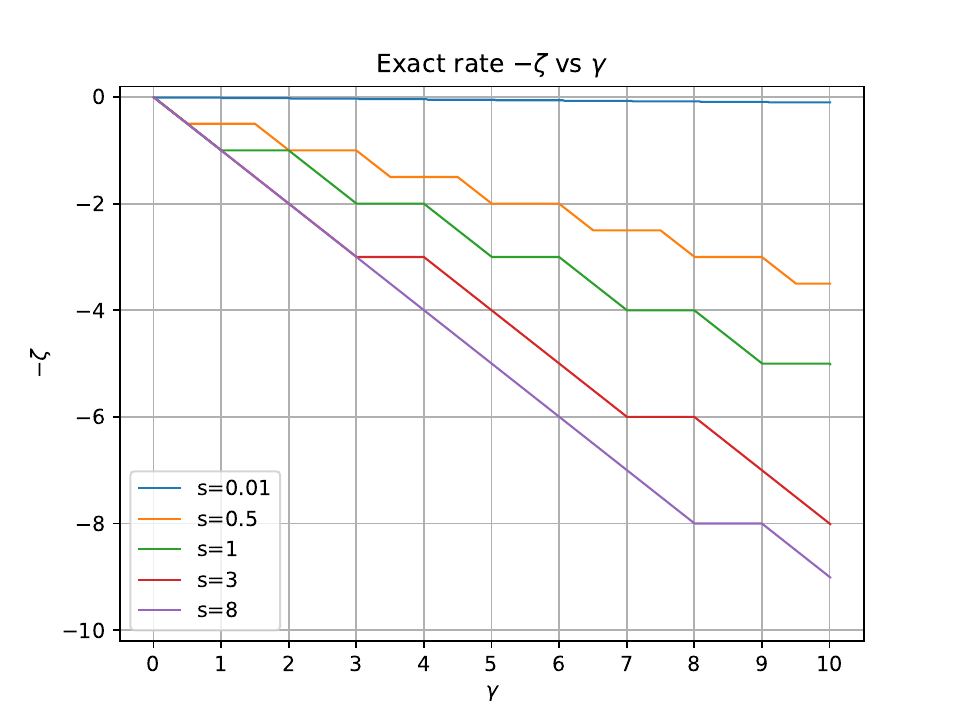}\label{figure_1_exact_rate}}
\subfigure[]{\includegraphics[width=0.48\columnwidth]{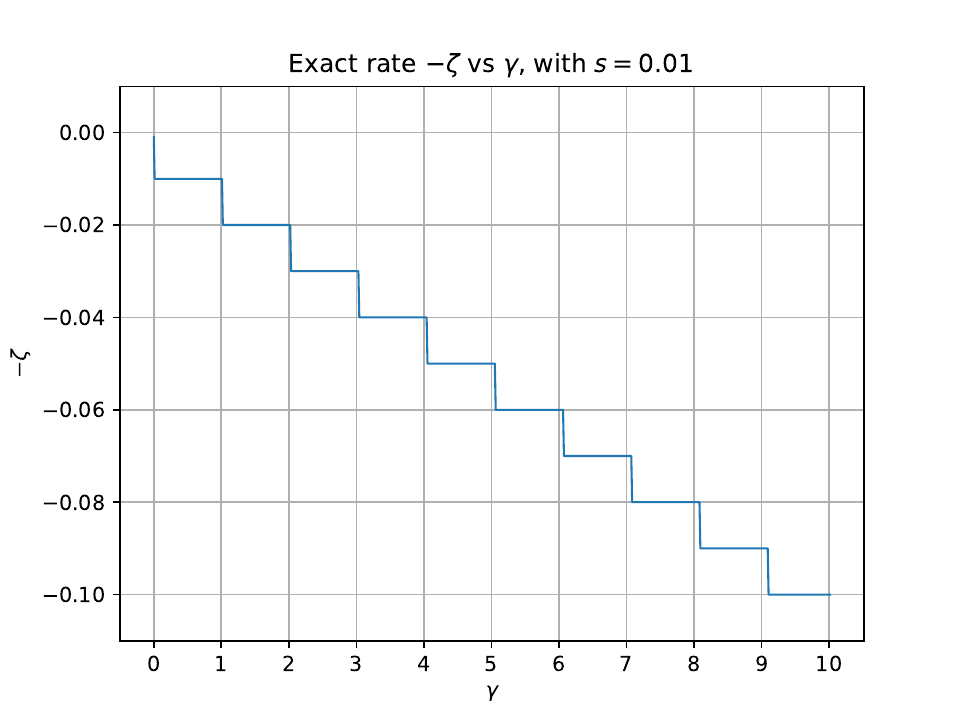}\label{fig_1_rate_0_0_1}}

\caption{
(a) A graphical representation of the exact rate of the minimax risk for kernel regression obtained in Theorem \ref{thm:main_pinsker_constant} with $s=0.01, 0.5, 1, 3$, and $8$.
(b) The exact rate when $s=0.01$.
}
\label{fig:fig_p}
% \vskip 0.05in
\end{figure}

Figure \ref{fig:fig_p} illustrates the curve of the exact rate $\zeta$ with respect to $\gamma$.
Theorem \ref{thm:main_pinsker_constant} and Figure \ref{fig:fig_p} reveal that periodical plateaus, where the rates $\zeta$ remain constant over a range of $\gamma$, occur for any $s>0$.
This phenomenon is termed periodic plateau behavior. As discussed in previous work \cite{lu2023optimal, zhang2024optimal, lu2024saturation}, the periodic plateau behavior suggests that improving the rate of minimax risk for kernel regression requires increasing the sample size above a certain threshold.

Although all plateaus demonstrated above are of length $1$, their proportion in each period (that is, $\gamma \in [p(s+1), (p+1)(s+1))$) gradually decreases as $s$ increases, which is approximately $\frac{1}{s+1}$.
%Specifically, as $s$ approaches infinity and $0 < \gamma < s$, the rate of minimax risk approaches $d^{-\gamma} = \Theta_d(n^{-1})$.
%This aligns with the intuition that kernel regression with sufficiently smooth regression functions (i.e., with large $s$) can achieve a faster minimax rate.

\begin{figure}[ht]
% \vskip 0.05in
\centering
\subfigure[]{\includegraphics[width=0.48\columnwidth]{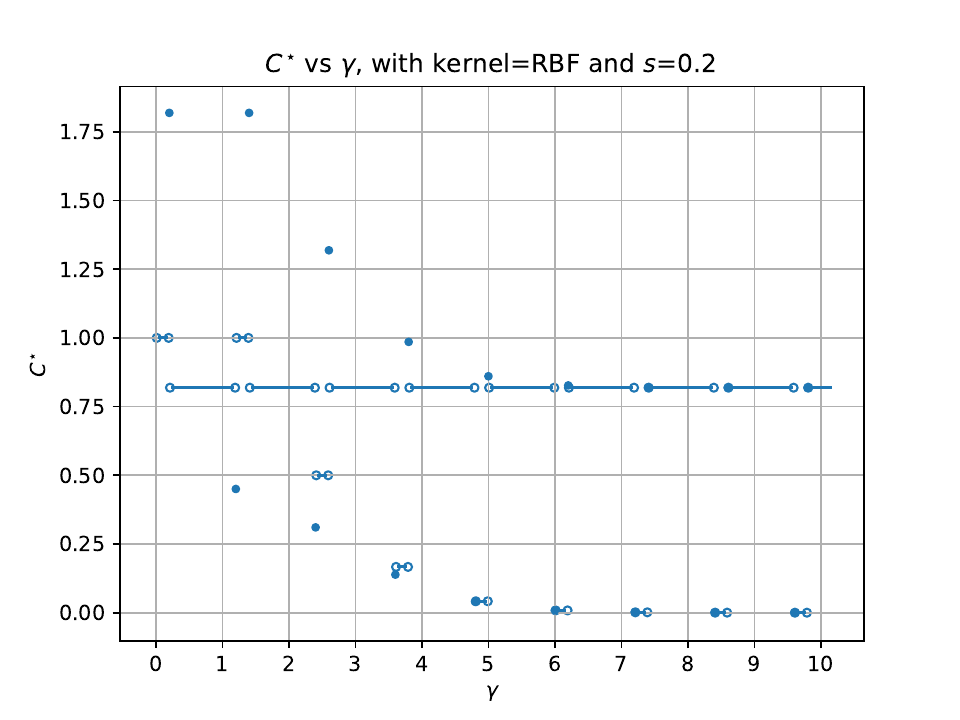}\label{fig_2_1}}
\subfigure[]{\includegraphics[width=0.48\columnwidth]{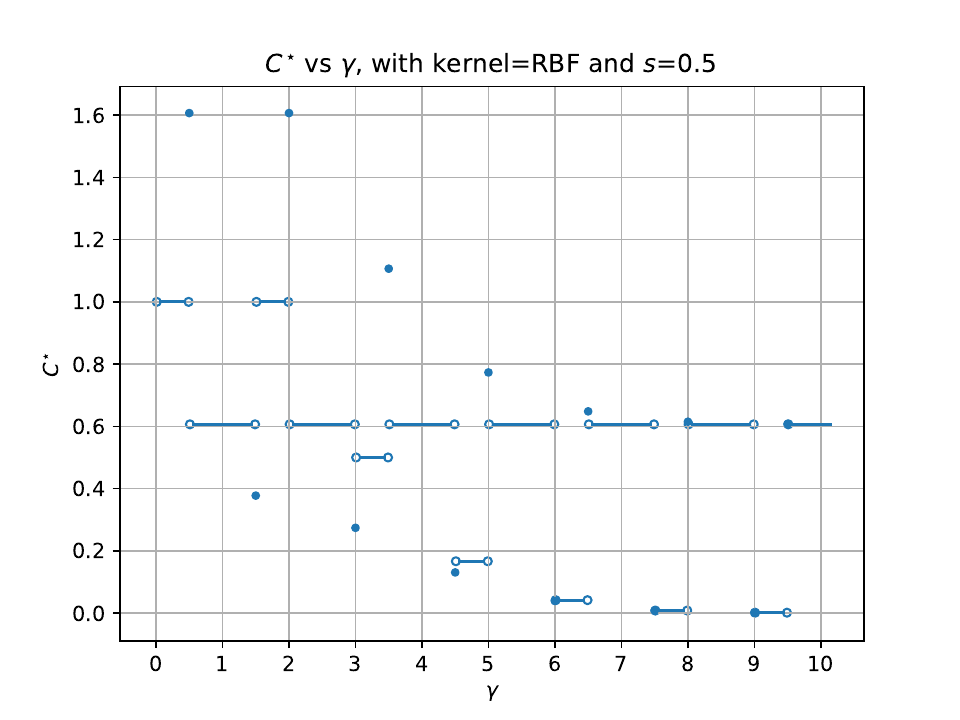}\label{fig_2_2}}

\subfigure[]{\includegraphics[width=0.48\columnwidth]{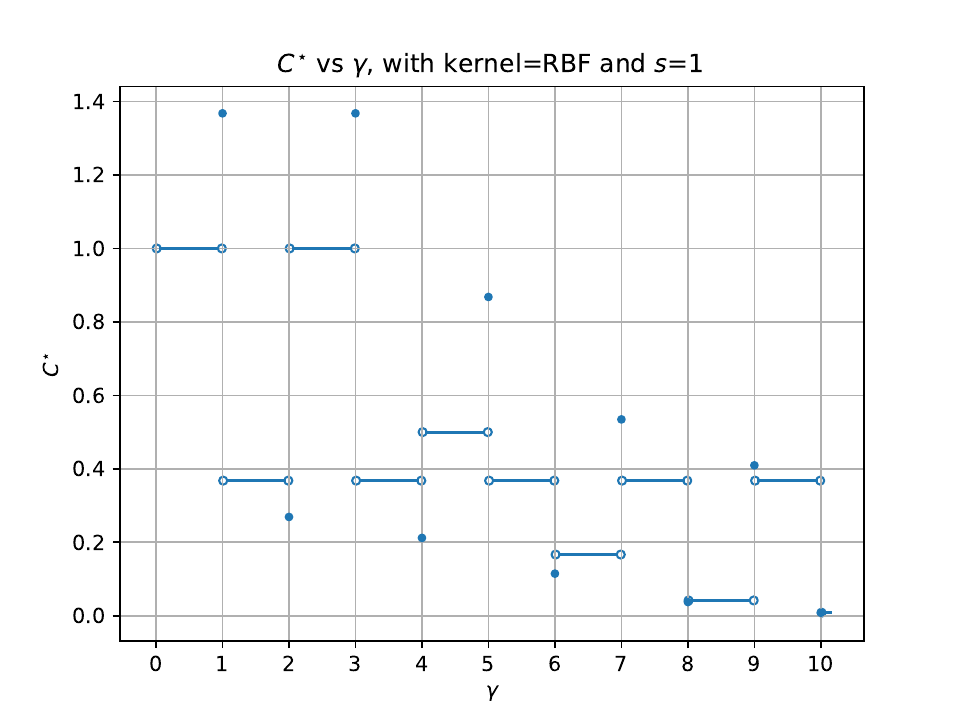}\label{fig_2_3}}
\subfigure[]{\includegraphics[width=0.48\columnwidth]{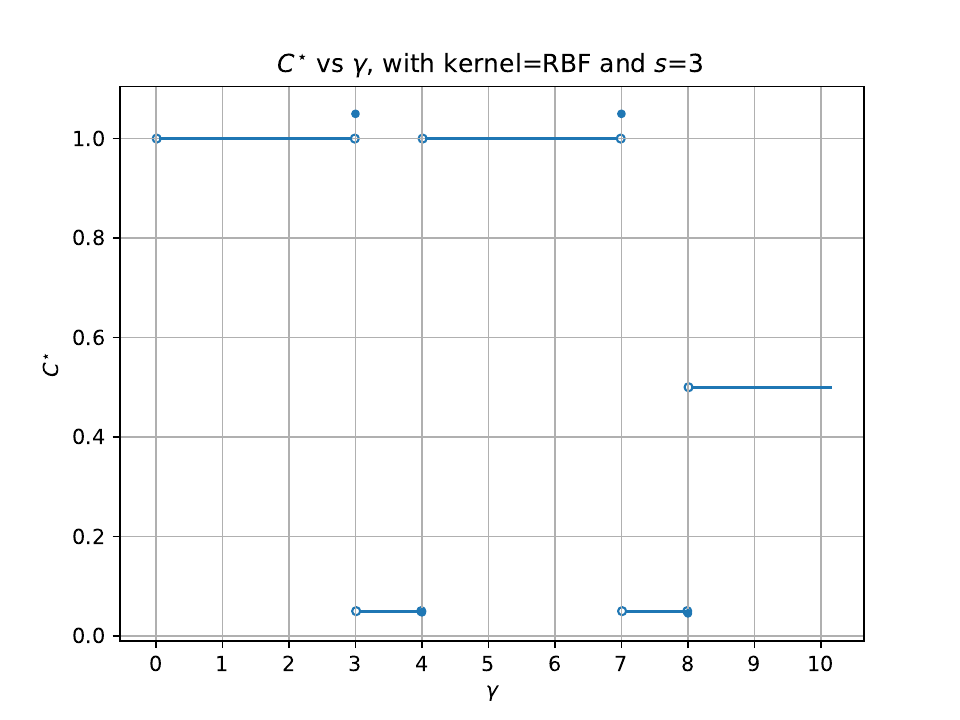}\label{fig_2_4}}

\caption{A graphical representation of the Pinsker constant of minimax risk of kernel regression problems obtained in Theorem \ref{thm:main_pinsker_constant}. We take $\alpha = R = \sigma=1$, and the kernel is the RBF kernel $K(x, x^{\prime}) := \exp(-\|x-x^{\prime}\|^2/2)$ (hence we have $a_p = 1/(e p!)$). In four subfigures, we choose $s=0.2, 0.5, 1$, and $3$.}
\label{figure_2_pinsker_constant}
% \vskip 0.05in
\end{figure}

\paragraph*{Pinsker constant}
Figure \ref{figure_2_pinsker_constant} illustrates the curve of Pinsker constant $\calC^{\star}$ with respect to $\gamma$, and we plot all the jump discontinuities of the Pinsker constant with solid dots. 

Pinsker’s constant represents a significant advancement in non-parametric estimation theory by enabling the comparison of estimators based on constants rather than just convergence rates (see, e.g., \cite{nussbaum1985spline, nussbaum1999minimax}). 
In parametric theory, these constants are expressed as "Fisher’s bound for asymptotic variances" with a corresponding rate of \(n^{-1}\) (\cite{bahadur1964fisher}). 
%Similarly, the Pinsker constant \(\calC^{\star}\) and the corresponding rate \(d^{-\zeta}\) given in Theorem \ref{thm:main_pinsker_constant} can be considered a non-parametric version of Fisher’s bound. 

%One may noticed an interesting scenario: when \(\gamma < s\), the Pinsker bound for the kernel regression problem (as described in Theorem \ref{thm:main_pinsker_constant}) is exactly $\sigma^2/n$, matching the Fisher’s bound for parametric models. 
One may have noticed an interesting scenario: when \( \gamma < s \), the Pinsker bound for the kernel regression problem, as described in Theorem \ref{thm:main_pinsker_constant}, is exactly \( \sigma^2 / n \).%, matching the Fisher's bound for parametric models. 
 To better understand that, notice that we have \(\|\mathrm{P}_{>0} f_{\star}\|_{L^2}^2 = o_d(1/n)\), indicating that the regression function can be approximated as a constant function. Therefore, the minimax risk for kernel regression is $\sigma^2/n + o_d(1/n)$.

%To better understand that, notice that we have \(\| \mathrm{P}_{>0} f_{\star} \|_{L^2}^2 = o_d(1/n)\), indicating that the regression function can be approximated as a constant function. Therefore, the minimax risk for kernel regression is \( \sigma^2 / n + o_d(1/n) \).

%However, it is important to note that the behavior of the problem differs depending on the value of \(\gamma\) within this range. Specifically, we can distinguish two cases according to Appendix \ref{subsec_para_case}:

\iffalse
\begin{itemize}
    \item[(a)] When $\gamma \leq s/2$ and $N=1$, the kernel regression model (\ref{equation:true_model}) degenerates to a parametric model, hence the minimax risk is exactly $\sigma^2/n$;

    \item[(b)] When $s/ 2 \leq \gamma < s$ and $N > 1$, we have \(\|\mathrm{P}_{>0} f_{\star}\|_{L^2}^2 = o_d(1/n)\), indicating that the regression function can be approximated as a constant function. Therefore, the minimax risk is $\sigma^2/n + o_d(1/n)$.
\end{itemize}
\fi

\iffalse
For example, 
\begin{itemize}
    \item[(a)] When $\gamma$ increases from $p(s+1)$ to $\gamma < p(s+1)+s$, $\calC^{\star}$ increases from $\frac{ \sigma^2 }{\alpha  p!  + \sigma^2 / (R a_p^{s} (p!)^{s})}$ to $\frac{ \sigma^2 }{\alpha  p! }$;
    \item[(b)] When $\gamma$ increases from $p(s+1)+s$ to $\gamma < (p+1)(s+1)$, $\calC^{\star}$ decreases from $R  {a_{p+1}^{s}} ((p+1)!)^{s}+\frac{ \sigma^2 }{\alpha  p! }$ to $R {a_{p+1}^{s}} ((p+1)!)^{s}$.
\end{itemize}
To explain (a), when $\gamma \in [p(s+1), p(s+1)+s)$, we can rewrite
$$
\calC^{\star} d^{-\zeta} \sim \frac{\sigma^2 N(d, p)}{n+\sigma^2 N(d, p) /\left(R \mu_p^s\right)}.
$$
and it is obvious that the Pinsker constant increases when $\gamma > p(s+1)$ because $\sigma^2 N(d, p) /\left(R \mu_p^s\right) =o_d(n)$.
A similar argument applies to other gaps in the Pinsker constant.

Interestingly, we notice that (b) can be explained more intuitively since the rate $\zeta$ remains constant for any $\gamma \geq p(s+1)+s$. When $\gamma = p(s+1)+s$, the sample size is not sufficient to fully capture the signal corresponding to $\mu_{p}$'s. %(see Remark \ref{remark_ell_j_values} for details). 
Therefore, the Pinsker constant for $\gamma=p(s+1)+s$ is larger than that for $\gamma>p(s+1)+s$.
\fi

    {\colorblue
    Notice that the Pinsker constant \( \mathcal{C}^{\star} \) decreases when \( \gamma \) increases from \( p(s+1) + s \) to \( \gamma \in (p(s+1) + s, (p+1)(s+1)) \). This is due to the fact that, for this range of \( \gamma \), the asymptotic form of the Pinsker bound is dominated by two terms (see Appendix \ref{lem:quantities_2} (ii)):
    \[
    \mathcal{D}^{\star} \sim R a_{p+1}^s ((p+1)!)^s d^{-(p+1)s} + \frac{\sigma^2}{\alpha p!} d^{p-\gamma}.
    \]
    When \( \gamma > p(s+1) + s \), one term on RHS becomes much larger than the other on RHS, leading to a reduction in the Pinsker constant. Interestingly, this can be explained more intuitively by noting that the rate \( \zeta \) remains constant for any \( \gamma \geq p(s+1) + s \). When \( \gamma = p(s+1) + s \), the sample size is insufficient to fully capture the signal corresponding to \( \mu_p \)'s. Hence, the Pinsker constant for \( \gamma = p(s+1) + s \) is larger than that for \( \gamma > p(s+1) + s \).
}

Lastly, we continue the discussion of the periodic plateau behavior. Recall that when $p(s+1)+s < \gamma < (p+1)(s+1)$, the exact rate $\zeta$ remains constant.
Likewise, we notice that the value of the Pinsker constant remains unchanged within each of these ranges. 
In other words, even if we merely want to reduce the Pinkser constant of the minimax risk, we might have to increase the sample size above a certain threshold.

\vspace{3mm}

\section{Calculation of $\calD^{\star} \sim \calC^{\star}d^{-\zeta}$}\label{sec:pinsker_quantity}

%Our technique of determining the Pinsker constant of interpolation spaces is partly based on the original idea of determining the Pinsker constant of the Gaussian sequence model in \cite{pinsker1980optimal} (one can also refer to \cite{tsybakov2008introduction}).Therefore, in this section, our first goal is to define a quantity $\calD^{\star}$ based on the dimension $d$ and all absolute constants defined in Definition \ref{def:abs_constants}. Then, we will show that $\calD^{\star} \sim \calC^{\star} d^{-\zeta}$, where $\calC^{\star}$ is the Pinsker constant defined in Theorem \ref{thm:main_pinsker_constant}.

Our technique for determining the Pinsker constant of interpolation spaces is partly inspired by the original method for determining the Pinsker constant of the Gaussian sequence model, as presented in Pinsker's seminal work \cite{pinsker1980optimal}. For further insights, one can refer to \cite{tsybakov2008introduction}. In this section, our initial objective is to define a quantity $\calD^{\star}$, which depends on the dimension 
$d$ and all the absolute constants outlined in Definition \ref{def:abs_constants}. Subsequently, we will demonstrate that 
$\calD^{\star}\sim \calC^{\star}d^{-\zeta}$, where $\calC^{\star}$
  is the Pinsker constant introduced in Theorem \ref{thm:main_pinsker_constant}.

Let's first define some quantities that are closely related to the Pinsker constant.

\begin{definition}\label{def:kappa_and_N}
    Denote $\kappa^{\star}$ as the unique solution (if it exists) to the following equation:
\begin{equation}\label{def:kappa}
\frac{\sigma^2}{n\kappa}\sum_{j=1}^{\infty}\lambda_j^{-s/2}(1-\kappa \lambda_j^{-s/2})_+=R,
\end{equation}
where $\lambda_j$'s are the eigenvalues of the kernel defined in Equation \eqref{eqn:maximum_of_kernel}. Moreover, let
\[
N:=\max\left\{j \geq 1~:~ \frac{\sigma^2}{n}\sum_{m=1}^{j}\lambda_m^{-s/2}(\lambda_j^{-s/2}-\lambda_m^{-s/2})<R\right\} \leq \infty.
\]
\end{definition}

%The following proposition ensures that $\kappa^{\star}$ exists and $N$ is finite.

Notice that when $s>0$, $\{\lambda_j^{-s/2}\}_{j=1}^{\infty}$ is a non-decreasing sequence and $\lambda_j^{-s/2} \to \infty$ as $j \to \infty$.
The following proposition restates the results of Lemma 3.1 and equation (3.19) from \cite{tsybakov2008introduction}, confirming the existence and uniqueness of $\kappa^{\star}$ and the finiteness of $N$.

\begin{proposition}\label{prop:exist_of_kappa_and_N}
    There exists a unique solution of (\ref{def:kappa}) given by
\begin{equation}\label{eqn:order_of_kappa}
\kappa^{\star} = \frac{\sigma^2\sum\nolimits_{j=1}^{N}\lambda_j^{-s/2}}{n R +\sigma^2\sum\nolimits_{j=1}^{N} \lambda_j^{-s}}.
\end{equation}
Furthermore, it is established that
\begin{equation}\label{eqn:order_of_N}
N=\max\left\{j:\lambda_j^{s/2}>\kappa^{\star}\right\} < \infty.
\end{equation}
\end{proposition}

Thanks to Proposition \ref{prop:exist_of_kappa_and_N}, we can now define $\calD^{\star}$ in terms of $\kappa^{\star}$ and $N$.

\begin{definition}\label{def:pinsker_quantity}
    For any $j \geq 1$, define $\ell_{j}$ as follows: 
    $$
    \ell_j := (1-\kappa^{\star} \lambda_j^{-s/2})_{+}
    $$
Furthermore, define
\begin{equation*}
    \calD^{\star} := \frac{\sigma^2}{n} \sum_{j=1}^{N} \ell_j,
\end{equation*}
where $\kappa^{\star}$ and $N$ are given in Definition \ref{def:kappa_and_N}.
\end{definition}

To demonstrate that 
$\calD^{\star} \sim \calC^{\star} d^{-\zeta}$, it is necessary to determine the asymptotic values of the leading eigenvalues 
$\lambda_j$'s, or equivalently, the asymptotic values of the 
leading eigenvalues $\mu_k$'s. The following lemma establishes the asymptotic equivalence of the leading eigenvalues
\(\mu_0, \ldots, \mu_{p+3}\)
  and their corresponding multiplicities, as defined in Equation (\ref{spherical_decomposition_of_inner_main}).

\begin{lemma}\label{lemma:inner_edr}
    Suppose Assumption 
\ref{assump_asymptotic} and 
\ref{assu:coef_of_inner_prod_kernel}  hold for some $\alpha, \gamma>0$.
    Let $p=\lfloor\frac{\gamma}{s+1} \rfloor$. 
    Then,
\begin{itemize}
    \item For any $k  = 0, 1, \ldots, p+3$, we have
\begin{equation*}
\begin{aligned}
 \mu_{k} \sim a_{k} k! d^{-k} \quad \text{ and } \quad
N(d, k) \sim \frac{d^{k}}{k!}.
\end{aligned}
\end{equation*}

\item There exists a constant \( \mathfrak{C}_1 \) only depending on the absolute constants \( \gamma, a_0, \ldots, a_{\lfloor \gamma \rfloor + 3} \) given in Definition \ref{def:abs_constants} such that for any \( d \geq \mathfrak{C}_1 \), we have 
\begin{gather*}
  0.9 \cdot a_{k} k! d^{-k} \leq  \mu_{k} \leq  1.1 \cdot a_{k} k! d^{-k} \quad \text{and} \quad 0.9 \cdot \frac{d^{k}}{k!} \leq
N(d, k) \leq 1.1 \cdot \frac{d^{k}}{k!},\\
  \mu_0 > \mu_1 > \cdots > \mu_{p+1} > \mu_{p+2} > \max_{j \geq p+3} \mu_{j}.
\end{gather*}
Consequently, if we denote $v_{-1}=0$ and $v_{\text{k}} = \sum_{k^{\prime}=0}^{k}N(d, k^{\prime})$,
%we have $\lambda_1 = \mu_0$,   $\phi_{1}=Y_{0,1}$, and 
then for any $0 \leq k \leq p+2$, we have:
    \begin{align*}
      \lambda_{v_{k-1}+1} = \lambda_{v_{k-1}+2}= \cdots=\lambda_{v_{k}} = \mu_k, ~ \{ 
         \phi_{v_{k-1}+1}, \phi_{v_{k-1}+2} \cdots, \phi_{v_{k}}\} =\{Y_{k,1},....,Y_{k,N(d,k)}\}.
    \end{align*}
\end{itemize}
\end{lemma}

\iffalse
{\color{red}
    It is noteworthy that the above result can be derived from Rodrigues' formula and Addition Formula. For a detailed exposition, refer to Appendix \ref{appendix_add_and_rog}. From Lemma \ref{lemma:inner_edr}, it follows that $\mu_k N(d, k) \sim a_k$, $k \leq p+3$.
    %From Lemma \ref{lemma:inner_edr} we have $\mu_k N(d, k) \sim a_k$, $k \leq p+3$. Notice that one can also obtain the above result by using Rodrigues’ formula and Addition Formula (see Appendix \ref{appendix_add_and_rog} for details).
}

\fi

Recall that the eigenvalues $\lambda_j$'s in (\ref{eqn:mercer_decomp}) are of non-increasing order, while the eigenvalues $\mu_k$'s in (\ref{spherical_decomposition_of_inner_main}) are not necessarily non-increasing. 
%Hence we can not directly apply Lemma \ref{lemma:inner_edr} on Proposition \ref{prop:exist_of_kappa_and_N} to obtain the value of $N$. 
Fortunately, from Lemma \ref{lemma:inner_edr} we can ensure the monotonicity of the leading eigenvalues \(\mu_0, \ldots, \mu_{p+3}\), and hence we can calculate the value of $N$, stated as the following lemma.

\begin{lemma}\label{lemma:calcul_N}
    Suppose the same conditions as Lemma \ref{lemma:inner_edr}. 
Then, there exists a constant \( \mathfrak{C} \) only depending on the absolute constants given in Definition \ref{def:abs_constants}, such that for any \( d \geq \mathfrak{C} \),we have $\mu_{p+2}^{s/2} \leq  \kappa^{\star} <\mu_p^{s/2}  $. Hence $\ell_{j}=0$ for any $j\geq v_{p+1}+1$ and
    \[N=v_{q}=\sum_{k=0}^{q}N(d,k)\]
    where the value of \(q\) is either equal to \(p\) or \(p+1\), depending on \(\alpha\), \(\gamma\), and the absolute constants in Definition \textit{\ref{def:abs_constants}}.
  Moreover, when $\gamma <p(s+1) + s/2$, we have $q=p$; when $\gamma > p(s+1) + s/2$, we have $q=p+1$. 
\end{lemma}

\begin{remark}\label{remark_unique_structure_of_eigvals}
We would like to point out that the periodic behavior of \( \zeta \) with respect to \( \gamma \) in Theorem \ref{thm:main_pinsker_constant} is closely related to the spectral properties of inner product kernels for uniform data distributed on a large-dimensional sphere.
 %   In Lemma \ref{lemma:inner_edr}, we have showed that we have $\mu_{k} = \Theta_d(d^{-k})$ and $N(d, k)  = \Theta_d(d^{k})$ for $k \leq p+3$.  Such a strong block structure of the spectrum forces $N$ to equal $N(d, q)$ for $q=p$ or $p+1$ (as shown in Lemma \ref{lemma:calcul_N}), and eventually let the rate of $\calD^{\star}$ in Definition \ref{def:pinsker_quantity} increase periodically with respect to $\gamma$.
    In Lemma \ref{lemma:inner_edr}, we have shown that \( \mu_k = \Theta_d(d^{-k}) \) and \( N(d, k) = \Theta_d(d^{k}) \) for \( k \leq p+3 \). 
    The strong block structure in the spectrum, as described, implies that \( N \) must equal \( v_q \) for \( q = p \) or \( p+1 \), as is demonstrated in Lemma \ref{lemma:calcul_N}. 
    This, in turn, results in a periodic decrease in the rate of \( \mathcal{D}^{\star} \) in Definition \ref{def:pinsker_quantity} with respect to \( \gamma \).
\end{remark}

Now we can calculate the Pinsker constant $\calC^{\star}$.

\begin{corollary}\label{coroll:quantities}
    Suppose Assumptions \ref{assump_asymptotic} and \ref{assu:coef_of_inner_prod_kernel} hold for some \(\alpha, \gamma > 0\). Then, when \(d \geq \mathfrak{C}\), where \(\mathfrak{C}\) is the constant defined in Lemma \ref{lemma:calcul_N}, we have 
    \[
    \mathcal{D}^{\star} \sim \mathcal{C}^{\star} d^{-\zeta},
    \]
    where \(\mathcal{D}^{\star}\) is given in Definition \ref{def:pinsker_quantity}, and \(\mathcal{C}^{\star}\) and \(\zeta\) are given in Theorem \ref{thm:main_pinsker_constant}.
\end{corollary}

\section{The matching upper bound}\label{sec:upper_sketch}

In this section, we provide a proof sketch showing that 
\begin{equation}\label{eqn:upper_bound_copy}
        \inf_{\hat{f}} \sup_{
        \rho_{f_{\star}} \in \mathcal{P}
        %f_{\star} \in \sqrt{R}[\calB]^{s}
        } \mathbb{E}_{(X, Y) \overset{\mathcal{D}}{\sim} \rho_{f_{\star}}^{\otimes n}}\left[
    \|\hat{f}-f_{\star}\|_{L^2}^2
    \right]
     \leq
    \mathcal{D}^{\star}(1+o_d(1)).
    \end{equation}
The detailed proof is deferred to Appendix \ref{sec:upper}.
For simplicity, we denote $\mathbb{E} = \mathbb{E}_{(X, Y) \overset{\mathcal{D}}{\sim} \rho_{f_{\star}}^{\otimes n}}$, where the distributions $\rho_{f_{\star}}$ on $\mathcal{X} \times \mathcal{Y}$ is given by (\ref{equation:true_model}), satisfying Assumption 
\ref{assump_asymptotic}, 
\ref{assu:coef_of_inner_prod_kernel},  and \ref{assump_function_calss} for some $\alpha, \gamma>0$.

For any $f_{\star}(\cdot) = \sum_{j=1}^{\infty} \theta_j \phi_j(\cdot) \in \sqrt{R}[\calB]^{s}$, denote $g_{\star}(x) = \sum_{j = 2}^{\infty} \theta_j \phi_j(x)$ where $\phi_j$'s are the eigenfunctions defined in (\ref{eqn:mercer_decomp}).
Let $\bar{z}_j := \frac{1}{n} \sum_{i=1}^{n} y_i \phi_j(x_i)$. We introduce the following linear filter estimator:
\begin{align*}
    \hat{f}_{\ell}(x) :=
    (\ell_1\mathbf{1}\{p=0\}+\mathbf{1}\{p > 0\}) \bar{z}_1 +\hat{g}_{\ell}(x) \quad
    \mbox{ where } \quad \hat{g}_{\ell}(x) = \sum_{j=2}^{N} \ell_j \bar{z}_j \phi_j(x),
\end{align*}
where $p = \lfloor\frac{\gamma}{s+1}\rfloor \geq 0$ is defined as in Theorem \ref{thm:main_pinsker_constant}.

For any $d \geq {\mathfrak{C}}$, where $\frak{C}$ is the sufficiently large constant defined in Lemma \ref{lemma:calcul_N}, we have $\phi_1=Y_{0, 1} \equiv 1$, hence $\mathbb{E}_{x}(g_{\star}(x))=\mathbb{E}_{x}(\hat{g}_{\ell}(x))=0$. It is clear that we have:
%we have the following decomposition:
\begin{equation*}
    \begin{aligned}
    &~ \inf_{\hat{f}} \sup_{
        \rho_{f_{\star}} \in \mathcal{P}
        %f_{\star} \in \sqrt{R}[\calB]^{s}
        } \mathbb{E}_{(X, Y) \overset{\mathcal{D}}{\sim} \rho_{f_{\star}}^{\otimes n}}\left[
    \|\hat{f}-f_{\star}\|_{L^2}^2
    \right]
        \leq  \sup_{
        \rho_{f_{\star}} \in \mathcal{P}
        %f_{\star} \in \sqrt{R}[\calB]^{s}
        }
        \mathbb{E}\left[
    \|\hat{f}_{\ell}-f_{\star}\|_{L^2}^2
    \right].
%    =&~
%    \left\{\begin{matrix}
%\underbrace{\sup_{
%        \rho_{f_{\star}} \in \mathcal{P}
%}
%\mathbb{E}\left[
%    \|\hat{g}_{\ell}-g_{\star}\|_{L^2}^2
%    \right]}_{\bold{I}}
%    +
%  \underbrace{  \sup_{
%        \rho_{f_{\star}} \in \mathcal{P}
%        %f_{\star} \in \sqrt{R}[\calB]^{s}
%        }\mathbb{E}\left(n^{-1} \sum_{i=1}^{n} y_i - \theta_1\right)^2}_{\bold{II}}, & ~ p>0\\ 
%\underbrace{\sup_{
%        \rho_{f_{\star}} \in \mathcal{P}}\mathbb{E}\left[
%    \|\hat{f}_{\ell}-f_{\star}\|_{L^2}^2
%    \right]}_{\bold{III}}, & ~ p=0
%\end{matrix}\right.\\
%    := &~
%    \mathbf{1}\{p>0\}\left(\mathbf{I} + \mathbf{II}\right) + \mathbf{1}\{p=0\}\mathbf{III}.
    \end{aligned}
    \end{equation*}

We first introduce the following  theorem, proof of which is deferred to Appendix \ref{append_zero_mean_case_1}.

\begin{theorem}\label{thm:upper_1_copy}
    Suppose the same conditions as Theorem \ref{thm:main_pinsker_constant}.
%    Further, suppose that $N =\sum_{k=0}^{p}N(d, k)$. 
    %When $d \geq \mathfrak{C}$, a sufficiently large constant only depending on the constants defined in Definition \ref{def:abs_constants}.
    Then, for any $\varepsilon > 0$, there exist a constant $D_{\varepsilon}$ only depending on $\varepsilon$ and  ${\mathfrak{C}}$ defined in Lemma \ref{lemma:calcul_N}, such that for any $d > D_{\varepsilon}$, and
    for any regression function $f_{\star} \in \sqrt{R}[\calB]^{s}$ satisfying one of the following conditions: (i) $\mathbb{E}_x f_{\star}(x) = 0$ or (ii) $p=0$, we have
    \begin{equation*}
       \mathbb{E}\left[
    \|\hat{f}_{\ell}\mathbf{1}\{p=0\}+ \hat{g}_{\ell}\mathbf{1}\{p>0\}-f_{\star}\|_{L^2}^2
    \right]
    \leq
    \mathcal{D}^{\star}(1+\varepsilon).
    \end{equation*}
    \iffalse
    Suppose the same conditions as Theorem \ref{thm:main_pinsker_constant}.
    Then, when \(d \geq \mathfrak{C}\) (the constant defined in Lemma \ref{lemma:calcul_N}) and $N =\sum_{k=0}^{p}N(d, k)$, we have the following statements:
    \begin{itemize}
        \item[(i)] If $p>0$, then we have
        $
        \sup_{
        \rho_{f_{\star}} \in \mathcal{P}} \mathbb{E}\left[
    \|\hat{g}_{\ell}-g_{\star}\|_{L^2}^2
    \right] 
    \leq
    \mathcal{D}^{\star}(1+o_d(1))
        $;

        \item[(ii)] If $p=0$, then we have
        $
        \sup_{
        \rho_{f_{\star}} \in \mathcal{P}} \mathbb{E}\left[
    \|\hat{f}_{\ell}-f_{\star}\|_{L^2}^2
    \right] 
    \leq
    \mathcal{D}^{\star}(1+o_d(1))
        $.
    \end{itemize}
    \fi
\end{theorem}

\iffalse

\begin{theorem}\label{thm:upper_2_copy}
    Suppose the same conditions as Theorem \ref{thm:main_pinsker_constant}. Further, suppose that $N =\sum_{k=0}^{p+1}N(d, k)$. 
    Then, for any $\varepsilon > 0$, there exist a constant $D_{\varepsilon}$ only depending on $\varepsilon$ and  ${\mathfrak{C}}$ defined in Lemma \ref{lemma:calcul_N}, such that for any $d > D_{\varepsilon}$, and
    for any regression function $f_{\star} \in \sqrt{R}[\calB]^{s}$ satisfying one of the following conditions: (i) $\mathbb{E}_x f_{\star}(x) = 0$ or (ii) $p=0$, we have
    \begin{equation*}
       \mathbb{E}\left[
    \|\hat{f}_{\ell}\mathbf{1}\{p=0\}+ \hat{g}_{\ell}\mathbf{1}\{p>0\}-f_{\star}\|_{L^2}^2
    \right]
    \leq
    \mathcal{D}^{\star}(1+\varepsilon).
    \end{equation*}
\end{theorem}
\fi

Now, let’s prove (\ref{eqn:upper_bound_copy}). Notice that when $p=0$, Theorem \ref{thm:upper_1_copy}
implies that
$$
\sup_{
        \rho_{f_{\star}} \in \mathcal{P}
        %f_{\star} \in \sqrt{R}[\calB]^{s}
        } 
    \mathbb{E}\left[
    \|\hat{f}_{\ell}-f_{\star}\|_{L^2}^2
    \right]
    \leq
    \mathcal{D}^{\star}(1+o_d(1)).
$$
As for the case where \(p > 0\), we have the following decomposition:
\begin{equation}
\begin{aligned}
    &\bbE\left[ \|\widehat{f}_{\ell}-f_{\star}\|^{2}_{L^{2}}\right]
    \leq
     \underbrace{\mathbb{E}\left(n^{-1} \sum\nolimits_{i=1}^{n} y_i - \theta_1\right)^2}_{\bold{I}}+\underbrace{
\mathbb{E}\left[
    \|\hat{g}_{\ell}-g_{\star}\|_{L^2}^2
    \right]}_{\bold{II}}.
\end{aligned}
\end{equation}
Since 
$\mathbb{E}(y_i \mid x_i)=\theta_1 + g_{\star}(x_i)$ and $\text{Var}(y_i \mid x_i)\leq \sigma^2$, for any $\varepsilon > 0$, there exists a constant $D_{\varepsilon, 1}$, depending only on $\varepsilon$ and  ${\mathfrak{C}}$ as defined in Lemma \ref{lemma:calcul_N}, such that for any $d > D_{\varepsilon, 1}$, and
    for any regression function $f_{\star}$ belonging to $\sqrt{R}[\calB]^{s}$, we have the following bound (see Theorem \ref{thm:upper_final} for a full proof):
\begin{equation*}
    \begin{aligned}
        \mathbf{I} \leq \frac{\sigma^2}{n} + \frac{\mu_1^{s}}{n} R \leq \mathcal{D}^{\star}\varepsilon.
    \end{aligned}
    \end{equation*}
Furthermore, from Theorem \ref{thm:upper_1_copy}, for any $d > D_{\varepsilon}$, and
    for any regression function $f_{\star} \in \sqrt{R}[\calB]^{s}$, we have
$$
\bold{II} \leq \mathcal{D}^{\star}(1+\varepsilon),
$$
hence when $d \geq \mathfrak{C}$, by the definition of $o_d(1)$, we have
$$
\sup_{
        \rho_{f_{\star}} \in \mathcal{P}
        %f_{\star} \in \sqrt{R}[\calB]^{s}
        } 
    \mathbb{E}\left[
    \|\hat{f}_{\ell}-f_{\star}\|_{L^2}^2
    \right] \leq \mathcal{D}^{\star}(1+o_d(1)).
$$

{\colorblue
\begin{remark}
Obtaining the upper bound in (\ref{eqn:upper_bound_copy}) is a challenging task due to several technical difficulties:
 \begin{itemize}
    \item In the Gaussian sequence model one observes
    $
    z_j=\int_0^1 \phi_j(t) d Y(t)=\theta_j+ \xi_j^{\text{normal}}$ with $\xi_j^{\text{normal}} {\sim}_{\text{i.i.d.}} \calN (0, \sigma^2)$,
    allowing a straightforward linear filter analysis (\cite{pinsker1980optimal}). In our kernel-regression framework only empirical estimators (refer to Eq.(\ref{def:bar_z_j})),
    $
    \bar{z}_j=\frac{1}{n} \sum_{i=1}^{n} y_i \phi_j(x_i)=\theta_{j} + \sum_{j^{\prime}=1}^{\infty} \theta_{j^{\prime}} \Delta_n(j, j^{\prime})+\xi_j,
    $
    are available. This replacement introduces an error term $\sum_{j^{\prime}=1}^{\infty} \theta_{j^{\prime}} \Delta_n(j, j^{\prime})$ and destroys the i.i.d. Gaussian structure of $\xi_j$, thereby significantly complicating the analysis.
    
    \item In fixed-dimensional Sobolev spaces with equidistant inputs on $[0, 1]^{d}$, the basis functions satisfy the so-called strong cancellation property, ensuring that $\Delta_n(j, j^{\prime}) \equiv 0$ (\cite{nussbaum1985spline, brown1996asymptotic, carter2006continous, reiss2008asymptotic}). 
    In contrast, spherical harmonics do not. In Appendix \ref{append_bound_E_1}, we developed new tools to control the interaction terms $\Delta_n(j, j^{\prime})$.

\end{itemize}
\end{remark}
}

\section{{\colorblue Equalness of Pinsker bounds for kernel regression model and Gaussian sequence model}}\label{sec:gauss_white_noise}

{\colorblue

In this section, we will obtain the Pinsker bound for an equivalent Gaussian sequence model, with eigenvalues $\lambda_j$ defined in (\ref{eqn:mercer_decomp}).
We will then show that this Pinsker bound is equal to the Pinsker bound for kernel regression model in Theorem 3.1.

Consider countably many observations
\begin{equation}\label{equation:GSM_model}
    z_j=\theta_j+ \varepsilon\xi_j, j=1,2, \cdots,
\end{equation}
where $\xi_j$ are i.i.d. $\calN\left(0, 1\right)$ and the sequence $\theta=\left(\theta_j\right)$ is in the following parameter space
$$
\Theta_{R}=\left\{\theta: \sum_{j=1}^{\infty} \lambda_j^{-s/2} \theta_j^2 \leq R\right\},
$$
where $\lambda_j$'s are the eigenvalues of the inner product kernel $K$ defined in (\ref{eqn:mercer_decomp}). 

Pinsker's result (\cite{pinsker1980optimal}) proposed to use the linear filtering estimator $\hat{\theta}^{c} = (c_j z_j)_{j \geq 1}$ to estimate $\theta$, where $c = (c_j)_{j \geq 1}$ is a sequence in $\ell^2$ such that $0 \leq c_j \leq 1$ for all $j$. The following results can be obtained by combining results in Lemma 3.2 in \cite{tsybakov2008introduction} and Corollary \ref{coroll:quantities}.

\begin{proposition}
[Restate Lemma 3.2 in \cite{tsybakov2008introduction}]\label{prop_tsyb_upper}
Let $\varepsilon^2 = \sigma^2 / n$.
    Suppose Assumption \ref{assump_asymptotic} and \ref{assu:coef_of_inner_prod_kernel} hold for some $\alpha>0$. Then we have
    $$
    \inf_{\hat{\theta}} \sup _{\theta \in \Theta_{R}} \mathbb{E}_\theta\|\hat{\theta}-\theta\|_{\ell^2}^2 
 \leq \sup _{\theta \in \Theta_{R}} \mathbb{E}_\theta\|\hat{\theta}^{\ell}-\theta\|_{\ell^2}^2
    %=\inf _c \sup _{\theta \in \Theta_{R}} \mathbb{E}_\theta\|\hat{\theta}^{c}-\theta\|_{\ell^2}^2 
    = \calD^{\star} \sim
    \calC^{\star} d^{-\zeta},
    $$
where 
$\hat{\theta}$ is any estimator of $\theta$, measurable with respect to the observed data set $\{z_j\}_{j=1}^{\infty}$,
$\ell = (\ell_j)$, $\ell_1, \cdots, \ell_N$ are given in Definition \ref{def:pinsker_quantity}, $\ell_j = 0$ for all $j >N$, and $\calC^{\star}$ and $\zeta$ are given in Theorem \ref{thm:main_pinsker_constant}.
\end{proposition}

Then, we can obtain the Pinsker bound for the above Gaussian sequence model based on Proposition \ref{prop_tsyb_upper} and Subsection 3.3.2 in \cite{tsybakov2008introduction}.

\begin{corollary}\label{coroll_asymp_minimax}
    Let $\varepsilon^2 = \sigma^2 / n$.
    Suppose Assumption \ref{assump_asymptotic} and \ref{assu:coef_of_inner_prod_kernel} hold for some $\alpha>0$. Then we have
    $$
    \inf_{\hat{\theta}} \sup _{\theta \in \Theta_{R}} \mathbb{E}_\theta\|\hat{\theta}-\theta\|_{\ell^2}^2 
    \sim
    \calC^{\star} d^{-\zeta},
    $$
where 
$\hat{\theta}$ is any estimator of $\theta$, measurable with respect to the observed data set $\{z_j\}_{j=1}^{\infty}$,
and $\calC^{\star}$ and $\zeta$ are given in Theorem \ref{thm:main_pinsker_constant}.
\end{corollary}

\begin{remark}
    For readers' convenience, we provide a quick proof of Corollary \ref{coroll_asymp_minimax} as follows.
    The upper bound is given by Proposition \ref{prop_tsyb_upper}.
    The lower bound can be obtained in the following way: (1) when $\gamma>s/2$, we can use the proof in Section 3.3.2 in \cite{tsybakov2008introduction} to get desired results, with (3.48) in \cite{tsybakov2008introduction} replaced by Appendix \ref{subsec_bound_of_resi}; (2) when $\gamma \leq s/2$, we can use the proof of Theorem \ref{thm_lower_para} instead.
\end{remark}

It is well known that Le Cam’s equivalence can, in many cases, reduce nonparametric problems to equivalent sequence models (\cite{le2000asymptotics, le2012asymptotic}).
However, we can not attain the Pinsker constant of large-dimensional kernel regression from Corollary \ref{coroll_asymp_minimax}. We would like to discuss existing literature and some of the challenges we encountered along the way.
\begin{itemize}
    \item[(i)] For fixed $d$, \cite{pinsker1980optimal} derived Pinsker bound for sequence model, and \cite{brown1996asymptotic, carter2006continous, reiss2008asymptotic} developed the Le Cam equivalence between kernel regression model over $[\mathcal{H}]^s$ ($s>1$) and sequence model. As a result, two models have same Pinsker bounds when the Le Cam equivalence holds. 
    However, the Le Cam equivalence fails for \( s \leq 1 \).
    In fact, \cite{efromovich1996asymptotic, brown1998asymptotic} 
    gave counterexamples
    that the Le Cam equivalence fails for \( s = 1/2 \) and for the boundary case \( s = 1 \) in the case of equidistant designs in \( [0,1]^{d} \).
    As a result, the Pinsker bound for kernel regression over \( \{[\mathcal{H}]^{s}, 0 < s \leq 1\} \) has not been established in the literature. 

    \item[(ii)] For large $d$ where $n\asymp d^\gamma$, whether Le Cam equivalence holds (even for $s>1$) is an open problem.
In fact, we derived our results without establishing the large-dimensional Le Cam equivalence. 
Consequently, we leveraged harmonic analysis on spheres and performed large-dimensional calculations involving eigenvalues to address this issue.

\end{itemize}

Nonetheless, notice that $\sqrt{R}[\calB]^{s}$ can be parametrized by the parameter space $\Theta_{R}$.
Hence, when Assumption \ref{assump_asymptotic} and \ref{assu:coef_of_inner_prod_kernel} hold,
Theorem \ref{thm:main_pinsker_constant} and Corollary \ref{coroll_asymp_minimax} build equalness between Pinsker bounds for kernel regression model (\ref{equation:true_model}) and Gaussian sequence model (\ref{equation:GSM_model}). We hope it offers heuristic evidence of a deeper connection between the two models, possibly even a new Le Cam equivalence.
}

\section{ Discussion}\label{sec:comparison}

This paper determined the exact asymptotic behavior of the minimax risk for kernel regression in large-dimensional settings. 
Specifically, we consider the nonparametric regression problem $y=f_{\star}(x)+\epsilon$, where the sample size $n \sim \alpha d^{\gamma}$ and $f_{\star} \in [\calH]^s$,  an interpolation space associated with an inner product kernel $K$ defined on the sphere $\bbS^{d}$.
As stated in Theorem \ref{thm:main_pinsker_constant}, the exact minimax risk bound is given by
\begin{align*}
    &~
    \inf_{\hat{f}} \sup_{
        \rho_{f_{\star}} \in \mathcal{P}
        %f_{\star} \in \sqrt{R}[\calB]^{s}
        } \mathbb{E}_{(X, Y) \overset{\mathcal{D}}{\sim} \rho_{f_{\star}}^{\otimes n}}\left[
    \|\hat{f}-f_{\star}\|_{L^2}^2
    \right] \sim \calC^{\star} d^{-\zeta}\\
    \sim 
    &~
    R  {a_{p+1}^{s}} ((p+1)!)^{s} d^{-(p+1)s}
        +
        \frac{ \sigma^2   }{\alpha  p!  + \sigma^2 / (R a_p^{s} (p!)^{s})\mathbf{1}\{\gamma=(s+1)p\}} d^{p-\gamma},
\end{align*}
where {\colorblue $\hat{f}$ is any estimator of $f_{\star}$, measurable with respect to the observed data set $(X, Y)$,} and $f_{\star}$ is in $\sqrt{R}[\calB]^{s} = \{f \in [\calH]^{s} \mid \|f\|_{[\calH]^{s}} \leq \sqrt{R}\}$, and all absolute constants above are given in Definition \ref{def:abs_constants}.

\iffalse
We also generalize our result and determine the exact asymptotic behavior of the minimax risk for the Gaussian sequence model $z_j=\theta_j+ \varepsilon\xi_j, j=1,2, \cdots$.
When the sequence $\theta=\left(\theta_j\right)$ is in the parameter space
$
\Theta=\{\theta: \sum_{j=1}^{\infty} \lambda_j^{-s/2} \theta_j^2 \leq R\}
$ with $\lambda_j$'s being the eigenvalues of the inner product kernel $K$, we show in Corollary \ref{coroll_asymp_minimax} that
$$
\inf_{\hat{\theta}} \sup _{\theta \in \Theta} \mathbb{E}_\theta\|\hat{\theta}-\theta\|_{\ell^2}^2 
    \sim
    \calC^{\star} d^{-\zeta}.
$$
\fi

%similar to the discussion in Subsection \ref{subsec:poly_approx_barrier}, we believe that certain spectral algorithms can achieve the Pinsker bound $\calC^{\star} d^{-\zeta}$ for any $s \geq 0$.

%In this section, let's review several existing works. We then provide some possible further research based on our results and the above works.
%\subsection{Polynomial approximation barrier}\label{subsec:poly_approx_barrier}

%
%Several studies have focused on the polynomial approximation barrier phenomenon 

%It is quite interesting to synthesize these studies to obtain several intriguing insights.
%Recently, a significant amount of research has been conducted on kernel regression in large-dimensional settings
%They found that, for square-integrable regression functions, kernel ridge regression (KRR) and kernel gradient flow are consistent if and only if the regression function is a polynomial of low degree.
It is quite interesting to compare our results with the extensive research conducted on kernel regression in large-dimensional settings (e.g., \cite{ghorbani2021linearized, Donhauser_how_2021, mei2022generalization, xiao2022precise, misiakiewicz_spectrum_2022, hu2022sharp}). Specifically, we restate Theorem 4 from \cite{ghorbani2021linearized} in the following proposition:
\begin{proposition}\label{restate_thm_4_in_linear}
%Given a dimension $d$, l
Let $f_{\star} \in L^2
%(\mathbb{S}^{d}, \rho_{\calX})
$ be a fixed regression function.
Suppose there exists an integer $\ell \in \{0, 1, \cdots\}$, and a constant $0<\delta < 1$, such that $n = \Theta_d(d^{\ell+1-\delta})$.
%Suppose further that Assumption 3 in \cite{ghorbani2021linearized} holds for the kernel $K$.
 Denote $\hat{f}_{\lambda}^{\krr}$ as the estimator of KRR and $R_{\text{KRR}}\left(f_{\star}, X, \lambda\right) := \mathbb{E}[\|\hat{f}_{\lambda}^{\krr}  - f_{\star}\|^2_{L^2}\mid X]$ as the conditional excess risk of KRR.\\
Under certain conditions,
    for any $\varepsilon>0$, and any regularization parameter $0<\lambda<\lambda^*$ ($\lambda^*$ is defined as (20) in \cite{ghorbani2021linearized}), there exists a constant $\mathfrak{C}_1$, such that if $d \geq \mathfrak{C}_1$, then with probability $1-o_d(1)$ we have
\begin{equation*}
    \left|R_{\text{KRR}}\left(f_{\star}, X, \lambda\right)-\left\|\mathrm{P}_{>\ell} f_{\star}\right\|_{L^2}^2\right| \leq \varepsilon\left(\left\|f_{\star}\right\|_{L^2}^2+\sigma^2\right).
\end{equation*}
\end{proposition}

%In summary, \cite{ghorbani2021linearized} revealed that when $s=0$ and $\gamma$ is not an integer, the excess risks of certain spectral algorithms remain close to the constant $\left\|\mathrm{P}_{>\lfloor\gamma\rfloor} f_{\star}\right\|_{L^2}^2$.

We observe that if the works of \cite{ghorbani2021linearized} and subsequent research could further obtain a union bound for \( R_{\text{KRR}}(f_{\star}, X, \lambda) \) over all functions \( f_{\star}\) in  \(\sqrt{R}[\mathcal{B}]^{0} \subseteq L^{2} \), then
\[
\sup_{f_{\star} \in \sqrt{R}[\mathcal{B}]^{0}} R_{\text{KRR}}(f_{\star}, X, \lambda) 
=
\sup_{f_{\star} \in \sqrt{R}[\mathcal{B}]^{0}} \|\mathrm{P}_{>\lfloor\gamma\rfloor} f_{\star}\|_{L^2}^2 (1+o_d(1))
= R(1+o_d(1)).
\]
This is intriguing because, by letting \( s \rightarrow 0 \) in our Pinsker's bound, we find
\[
\mathcal{C}^{\star} d^{-\zeta} = \lim_{s \to 0} R a_{\gamma+1}^{s} ((\gamma+1)!)^{s} d^{-s(\gamma+1)} = R.
\]
In other words, the conclusions of \cite{ghorbani2021linearized} and subsequent works align with our findings, particularly in the limit as \( s \) approaches zero.

On the other hand, when \( s > 0 \), Proposition \ref{restate_thm_4_in_linear} is not precise enough to provide an exact {\colorblue minimax} rate, even if the above union bound is obtained. Notice that we have
\[
\sup_{f_{\star} \in \sqrt{R}[\mathcal{B}]^{s}}\left\|\mathrm{P}_{>\lfloor\gamma\rfloor} f_{\star}\right\|_{L^2}^2  =
\mu_{\lfloor\gamma\rfloor+1}^{s}R = \Theta_d(d^{-s(\lfloor\gamma\rfloor+1)}),
\]
on the contrary, from Theorem \ref{thm:main_pinsker_constant} we know that the minimax rate is \( \Theta_d(d^{-\min\{\gamma - p, s(p+1)\}}) \) with \( p = \left\lfloor \frac{\gamma}{s+1} \right\rfloor \leq \lfloor\gamma\rfloor \).

{\colorblue
Two recent studies (\cite{cheng2022dimension, misiakiewicz2024non}) established concentration bounds for 
(i) the conditional excess risk of kernel ridge regression (KRR) in kernel regression model and (ii) the excess risk of ridge regression (RR) on the Gaussian sequence model. 
Specifically, they consider the following two settings:
\begin{itemize}
    \item[(i)]  They consider the kernel regression model (\ref{equation:true_model})
    with a regression function $f_{\star} = \sum_j \theta_j \phi_j \in L^2$.
    Specifying a kernel $K$ with eigenvalues $\lambda_j$'s, they then consider the KRR estimator with regularization parameter $\lambda$. The conditional excess risk is defined as
    $R_{\text{KRR}}\left(f_{\star}, X, \lambda\right)$;
    \item[(ii)] They also consider the Gaussian sequence model with 
    a specific variance of the Gaussian noise.
    %$\varepsilon^2 = \omega^2/n$, where $\omega$ is given in (20) of \cite{cheng2022dimension}.
    Let $\lambda_{\star} = \lambda_{\star}(\lambda)$ be given as in (7) of \cite{cheng2022dimension} and $R_{\text{RR}}(\lambda_{\star})$ be the excess risk of the RR estimator with regularization level $\lambda_{\star}$.
    
\end{itemize}

Under certain assumptions on the kernel $K$ and the regression function $f_{\star}$, \cite{cheng2022dimension, misiakiewicz2024non} proved that $|R_{\text{KRR}}\left(f_{\star}, X, \lambda\right) - R_{\text{RR}}(\lambda_{\star})|=o_d(R_{\text{RR}}(\lambda_{\star}))$ with high probability, as stated in the following propositions.

\begin{proposition}[Restate Theorem 1 in \cite{cheng2022dimension}]\label{prop_restate_cheng}
    Given a dimension $d$, let $f_{\star} \in \calH$ be a fixed regression function.
    Suppose that $\mathbb{E}\phi_j = 0$, $j=1, \cdots$. 
    Further suppose that there exists a constant $C>0$, such that for any $1$-Lipschitz convex function $\varphi: \mathbb{R}^{\infty} \rightarrow \mathbb{R}$, and for every $t>0$, we have
$$
\mathbb{P}\left(\left|\varphi\left(z_i\right)-\mathbb{E} \varphi\left(z_i\right)\right| \geq t\right) \leq 2 \exp \left(-t^2 / C^2\right),
$$
where $z_i = (\phi_1(x_i), \phi_2(x_i),\cdots)^{\top}$, $i \leq n$.
    Then under certain conditions, with probability $1-o_d(1)$, we have
    $$
    \left|R_{\text{KRR}}\left(f_{\star}, X, \lambda\right) - R_{\text{RR}}(\lambda_{\star})\right| = o_{d}(R_{\text{RR}}(\lambda_{\star})).
    $$
\end{proposition}

\begin{proposition}[Restate Theorem 2 in \cite{misiakiewicz2024non}]\label{prop_restate_misia}
    Given a dimension $d$, let $f_{\star} \in L^2$ be a fixed regression function.
    Suppose Assumption \ref{assump_asymptotic} and \ref{assu:coef_of_inner_prod_kernel} hold for some $\alpha, \gamma > 0$. 
    Denote $\ell = \lfloor\gamma\rfloor$.
    Suppose there exists a constant $C$, such that $\left\|\mathrm{P}_{>\ell} f_{\star}\right\|_{L^2} \geq\left\|f_{\star}\right\|_{L^2} / C$, and for any integer $q \geq 2$,  we have $\left\|f_{\star}\right\|_{L^q} \leq (C q)^{(\ell+1) / 2}\left\|f_{\star}\right\|_{L^2}$. Then under certain conditions, with probability $1-o_d(1)$, we have
    $$
    \left|R_{\text{KRR}}\left(f_{\star}, X, \lambda\right) - R_{\text{RR}}(\lambda_{\star})\right| = O_{d}\left(\log^{3(\ell+2)}(d) \cdot \left(\sqrt{\frac{d^{\ell-1}}{n}}+\sqrt{\frac{n}{d^{\ell+1}}}\right)\right)R_{\text{RR}}(\lambda_{\star}).
    $$
\end{proposition}

These results imply that the exact order of excess risk of the KRR is possibly same as the the exact order of excess risk of ridge estimator in sequence model (when $d \to \infty$). In particular, when ridge estimator in sequence model is minimax optimal, KRR is also minimax optimal.
However, they are insufficient for us to directly derive our Pinsker bound from sequence models:
    \begin{itemize}
        \item The saturation effect demonstrates that for %certain values of 
        \( s > 1 \), KRR cannot achieve the minimax rate (\cite{ zhang2024optimal, lu2024saturation}).
        
        \item Even when KRR achieves the minimax rate, our results \cite{li2024generalization} suggest that for a class of analytic spectral algorithms ( including the gradient flow, gradient descent, KRR etc.) cannot attain the constant optimality on excess risk. Hence, we can not determine the Pinsker constant through KRR.
        
        \item Their assumptions are incompatible with ours. For example, inner product kernels defined on the sphere do not satisfy the conditions in Proposition 7.2 since $\mathbb{E}Y_{0, 1}=1$. 
Similarly, functions in $\sqrt{R}[\calB]^{s}$ 
%that are not smooth or 
with non-zero $L^2$ norms do not satisfy the conditions in Proposition 7.3 since $\left\|\mathrm{P}_{>\ell} f_{\star}\right\|_{L^2} \to 0$.
    \end{itemize}
}

Finally, Theorem \ref{thm:main_pinsker_constant} strongly suggests that related nonparametric estimation problems with similar structures are worth considering, such as density estimation \cite{efroimovich1981estimation, golubev1994nonparametric}, 
%regression models with non-Gaussian noise or random designs \cite{golubev1990risk, efromovich1996nonparametric}, 
Besov bodies and wavelet estimation \cite{donoho1994minimax, donoho1990minimax}, and analogs of Theorem \ref{thm:main_pinsker_constant} when the square loss is substituted by other types of losses \cite{korostelev1993exact, tsybakov1997asymptotically}. 
Moreover, since our results heavily rely on the rotation-invariant property of the inner product kernels on the sphere (see, e.g., Remark \ref{remark_unique_structure_of_eigvals}), we believe that determining Pinsker bounds for other types of kernels on general domains in \(\mathbb{R}^{d}\) remains a more challenging question for future work.

%\section{Conclusions}\label{sec:con}

\section*{Acknowledgments}
Lin’s research was supported in part by the National Natural Science Foundation of China (Grant 92370122, Grant 11971257). This work has been partially supported by the New Cornerstone Science Foundation.
The authors would like to thank the anonymous referees, the Associate Editor, and the Editor for their constructive comments that improved the quality of this paper.

%Bibliography
\bibliographystyle{unsrt}  
\bibliography{reference.bib}

\begin{thebibliography}{10}

\bibitem{pinsker1980optimal}
Mark~Semenovich Pinsker.
\newblock Optimal filtering of square-integrable signals in gaussian noise.
\newblock {\em Problemy Peredachi Informatsii}, 16(2):52--68, 1980.

\bibitem{nussbaum1985spline}
Michael Nussbaum.
\newblock Spline smoothing in regression models and asymptotic efficiency in l2.
\newblock {\em The Annals of Statistics}, pages 984--997, 1985.

\bibitem{brown1996asymptotic}
Lawrence~D Brown and Mark~G Low.
\newblock Asymptotic equivalence of nonparametric regression and white noise.
\newblock {\em The Annals of Statistics}, 24(6):2384--2398, 1996.

\bibitem{le2012asymptotic}
Lucien Le~Cam.
\newblock {\em Asymptotic methods in statistical decision theory}.
\newblock Springer Science \& Business Media, 2012.

\bibitem{le2000asymptotics}
Lucien~Marie Le~Cam and Grace~Lo Yang.
\newblock {\em Asymptotics in statistics: some basic concepts}.
\newblock Springer Science \& Business Media, 2000.

\bibitem{efroimovich1981estimation}
S~Yu Efroimovich and Mark~Semenovich Pinsker.
\newblock Estimation of square-integrable density on the basis of a sequence of observations.
\newblock {\em Problemy Peredachi Informatsii}, 17(3):50--68, 1981.

\bibitem{golubev1994nonparametric}
Georgii~Ksenofontovich Golubev.
\newblock Nonparametric estimation of smooth spectral densities of gaussian stationary sequences.
\newblock {\em Theory of Probability \& Its Applications}, 38(4):630--639, 1994.

\bibitem{golubev1990risk}
Grigori~K. Golubev and Michael Nussbaum.
\newblock A risk bound in sobolev class regression.
\newblock {\em The Annals of Statistics}, 18(2):758--778, 1990.

\bibitem{efromovich1996nonparametric}
Sam Efromovich.
\newblock On nonparametric regression for iid observations in a general setting.
\newblock {\em The Annals of Statistics}, 24(3):1126--1144, 1996.

\bibitem{donoho1994minimax}
David~L Donoho and Iain~M Johnstone.
\newblock Minimax risk over $l_p$-balls for $l_q$-error.
\newblock {\em Probability Theory and Related Fields}, 99:277--303, 1994.

\bibitem{donoho1990minimax}
David~L Donoho, Richard~C Liu, and Brenda MacGibbon.
\newblock Minimax risk over hyperrectangles, and implications.
\newblock {\em The Annals of Statistics}, pages 1416--1437, 1990.

\bibitem{korostelev1993exact}
A.~P. Korostelev.
\newblock An asymptotically minimax regression estimator in the uniform norm up to exact constant.
\newblock {\em Theory of Probability \& Its Applications}, 38(4):737--743, 1994.

\bibitem{tsybakov1997asymptotically}
Aleksandr~Borisovich Tsybakov.
\newblock Asymptotically efficient signal estimation in l\_2 under general loss functions.
\newblock {\em Problemy Peredachi Informatsii}, 33(1):94--106, 1997.

\bibitem{nussbaum1999minimax}
Michael Nussbaum.
\newblock Minimax risk: Pinsker bound.
\newblock {\em Encyclopedia of Statistical Sciences}, 3:451--460, 1999.

\bibitem{le1953some}
Lucien Le~Cam.
\newblock On some asymptotic properties of maximum likelihood estimates and related bayes' estimates.
\newblock {\em University of California Publications in Statistics}, 1:277--330, 1953.

\bibitem{brown1997superefficiency}
Lawrence~D Brown, Mark~G Low, and Linda~H Zhao.
\newblock Superefficiency in nonparametric function estimation.
\newblock {\em The Annals of Statistics}, 25(6):2607--2625, 1997.

\bibitem{van1997superefficiency}
Aad~W van~der Vaart.
\newblock Superefficiency.
\newblock In {\em Festschrift for Lucien Le Cam: Research Papers in Probability and Statistics}, pages 397--410. Springer, 1997.

\bibitem{cai2005nonparametric}
T.~Tony Cai and Mark~G. Low.
\newblock {Nonparametric estimation over shrinking neighborhoods: Superefficiency and adaptation}.
\newblock {\em The Annals of Statistics}, 33(1):184 -- 213, 2005.

\bibitem{Jacot_NTK_2018}
Arthur Jacot, Franck Gabriel, and Cl{\'e}ment Hongler.
\newblock Neural tangent kernel: Convergence and generalization in neural networks.
\newblock {\em Advances in Neural Information Processing Systems}, 31, 2018.

\bibitem{Arora_on_2019}
Sanjeev Arora, Simon~S Du, Wei Hu, Zhiyuan Li, Russ~R Salakhutdinov, and Ruosong Wang.
\newblock On exact computation with an infinitely wide neural net.
\newblock {\em Advances in Neural Information Processing Systems}, 32, 2019.

\bibitem{Du_gradient_2019_b}
Simon Du, Jason Lee, Haochuan Li, Liwei Wang, and Xiyu Zhai.
\newblock Gradient descent finds global minima of deep neural networks.
\newblock In {\em International Conference on Machine Learning}, pages 1675--1685. PMLR, 2019.

\bibitem{Du_gradient_2019_a}
Simon~S Du, Xiyu Zhai, Barnabas Poczos, and Aarti Singh.
\newblock Gradient descent provably optimizes over-parameterized neural networks.
\newblock {\em arXiv preprint arXiv:1810.02054}, 2018.

\bibitem{jianfa2022generalization}
Jianfa Lai, Manyun Xu, Rui Chen, and Qian Lin.
\newblock Generalization ability of wide neural networks on $\mathbb{R}$.
\newblock {\em arXiv preprint arXiv:2302.05933}, 2023.

\bibitem{li2023statistical}
Yicheng Li, Zixiong Yu, Guhan Chen, and Qian Lin.
\newblock On the eigenvalue decay rates of a class of neural-network related kernel functions defined on general domains.
\newblock {\em Journal of Machine Learning Research}, 25(82):1--47, 2024.

\bibitem{Caponnetto2006OptimalRF}
Andrea Caponnetto.
\newblock Optimal rates for regularization operators in learning theory.
\newblock Technical Report CBCL Paper \#264/AI Technical Report \#062, Massachusetts Institute of Technology, September 2006.

\bibitem{caponnetto2007optimal}
Andrea Caponnetto and Ernesto De~Vito.
\newblock Optimal rates for the regularized least-squares algorithm.
\newblock {\em Foundations of Computational Mathematics}, 7(3):331--368, 2007.

\bibitem{raskutti2014early}
Garvesh Raskutti, Martin~J. Wainwright, and Bin Yu.
\newblock Early stopping and non-parametric regression: An optimal data-dependent stopping rule.
\newblock {\em Journal of Machine Learning Research}, 15(11):335--366, 2014.

\bibitem{Lin_Optimal_2020}
Junhong Lin, Alessandro Rudi, Lorenzo Rosasco, and Volkan Cevher.
\newblock Optimal rates for spectral algorithms with least-squares regression over hilbert spaces.
\newblock {\em Applied and Computational Harmonic Analysis}, 48(3):868--890, may 2020.

\bibitem{zhang2023optimality}
Haobo Zhang, Yicheng Li, Weihao Lu, and Qian Lin.
\newblock On the optimality of misspecified kernel ridge regression.
\newblock In {\em International Conference on Machine Learning}, pages 41331--41353. PMLR, 2023.

\bibitem{pmlr-v99-rakhlin19a}
Alexander Rakhlin and Xiyu Zhai.
\newblock Consistency of interpolation with laplace kernels is a high-dimensional phenomenon.
\newblock In {\em Conference on Learning Theory}, pages 2595--2623. PMLR, 2019.

\bibitem{beaglehole2022kernel}
Daniel Beaglehole, Mikhail Belkin, and Parthe Pandit.
\newblock On the inconsistency of kernel ridgeless regression in fixed dimensions.
\newblock {\em SIAM Journal on Mathematics of Data Science}, 5(4):854--872, 2023.

\bibitem{buchholz2022_KernelInterpolation}
Simon Buchholz.
\newblock Kernel interpolation in sobolev spaces is not consistent in low dimensions.
\newblock In {\em Conference on Learning Theory}, pages 3410--3440. PMLR, 2022.

\bibitem{li2023kernel}
Yicheng Li, Haobo Zhang, and Qian Lin.
\newblock Kernel interpolation generalizes poorly.
\newblock {\em Biometrika}, 111(2):715--722, 2024.

\bibitem{Bordelon_Spectrum_2020}
Blake Bordelon, Abdulkadir Canatar, and Cengiz Pehlevan.
\newblock Spectrum dependent learning curves in kernel regression and wide neural networks.
\newblock In {\em International Conference on Machine Learning}, pages 1024--1034. PMLR, 2020.

\bibitem{Cui2021GeneralizationER}
Hugo Cui, Bruno Loureiro, Florent Krzakala, and Lenka Zdeborov{\'a}.
\newblock Generalization error rates in kernel regression: The crossover from the noiseless to noisy regime.
\newblock {\em Advances in Neural Information Processing Systems}, 34:10131--10143, 2021.

\bibitem{jin2021learning}
Hui Jin, Pradeep~Kr Banerjee, and Guido Mont{\'u}far.
\newblock Learning curves for gaussian process regression with power-law priors and targets.
\newblock {\em arXiv preprint arXiv:2110.12231}, 2021.

\bibitem{li2023asymptotic}
Yicheng Li, Haobo Zhang, and Qian Lin.
\newblock On the asymptotic learning curves of kernel ridge regression under power-law decay.
\newblock {\em Advances in Neural Information Processing Systems}, 36:49341--49364, 2024.

\bibitem{li2024generalization}
Yicheng Li, Weiye Gan, Zuoqiang Shi, and Qian Lin.
\newblock Generalization error curves for analytic spectral algorithms under power-law decay.
\newblock {\em arXiv preprint arXiv:2401.01599}, 2024.

\bibitem{Liang_Just_2019}
Tengyuan Liang and Alexander Rakhlin.
\newblock {Just interpolate: Kernel “Ridgeless” regression can generalize}.
\newblock {\em The Annals of Statistics}, 48(3):1329 -- 1347, 2020.

\bibitem{ghorbani2021linearized}
Behrooz Ghorbani, Song Mei, Theodor Misiakiewicz, and Andrea Montanari.
\newblock {Linearized two-layers neural networks in high dimension}.
\newblock {\em The Annals of Statistics}, 49(2):1029 -- 1054, 2021.

\bibitem{mei2021learning}
Song Mei, Theodor Misiakiewicz, and Andrea Montanari.
\newblock Learning with invariances in random features and kernel models.
\newblock In {\em Conference on Learning Theory}, pages 3351--3418. PMLR, 2021.

\bibitem{Ghosh_three_2021}
Nikhil Ghosh, Song Mei, and Bin Yu.
\newblock The three stages of learning dynamics in high-dimensional kernel methods.
\newblock {\em arXiv preprint arXiv:2111.07167}, 2021.

\bibitem{mei2022generalization}
Song Mei, Theodor Misiakiewicz, and Andrea Montanari.
\newblock Generalization error of random feature and kernel methods: Hypercontractivity and kernel matrix concentration.
\newblock {\em Applied and Computational Harmonic Analysis}, 59:3--84, 2022.

\bibitem{misiakiewicz_learning_2021}
Theodor Misiakiewicz and Song Mei.
\newblock Learning with convolution and pooling operations in kernel methods.
\newblock {\em Advances in Neural Information Processing Systems}, 35:29014--29025, 2022.

\bibitem{aerni2023strong}
Michael Aerni, Marco Milanta, Konstantin Donhauser, and Fanny Yang.
\newblock Strong inductive biases provably prevent harmless interpolation.
\newblock {\em arXiv preprint arXiv:2301.07605}, 2023.

\bibitem{barzilai2023generalization}
Daniel Barzilai and Ohad Shamir.
\newblock Generalization in kernel regression under realistic assumptions.
\newblock {\em arXiv preprint arXiv:2312.15995}, 2023.

\bibitem{liang2020multiple}
Tengyuan Liang, Alexander Rakhlin, and Xiyu Zhai.
\newblock On the multiple descent of minimum-norm interpolants and restricted lower isometry of kernels.
\newblock In {\em Conference on Learning Theory}, pages 2683--2711. PMLR, 2020.

\bibitem{zhang2024phase}
Haobo Zhang, Weihao Lu, and Qian Lin.
\newblock The phase diagram of kernel interpolation in large dimensions.
\newblock {\em arXiv preprint arXiv:2404.12597}, 2024.

\bibitem{lu2023optimal}
Weihao Lu, Haobo Zhang, Yicheng Li, Manyun Xu, and Qian Lin.
\newblock Optimal rate of kernel regression in large dimensions.
\newblock {\em arXiv preprint arXiv:2309.04268}, 2023.

\bibitem{zhang2024optimal}
Haobo Zhang, Yicheng Li, Weihao Lu, and Qian Lin.
\newblock Optimal rates of kernel ridge regression under source condition in large dimensions.
\newblock {\em arXiv preprint arXiv:2401.01270}, 2024.

\bibitem{lu2024saturation}
Weihao Lu, Haobo Zhang, Yicheng Li, and Qian Lin.
\newblock On the saturation effects of spectral algorithms in large dimensions.
\newblock 2024.

\bibitem{Karoui_spectrum_2010}
Noureddine~El Karoui.
\newblock {The spectrum of kernel random matrices}.
\newblock {\em The Annals of Statistics}, 38(1):1 -- 50, 2010.

\bibitem{bartlett2021deep}
Peter~L Bartlett, Andrea Montanari, and Alexander Rakhlin.
\newblock Deep learning: a statistical viewpoint.
\newblock {\em Acta numerica}, 30:87--201, 2021.

\bibitem{pandit2024universality}
Parthe Pandit, Zhichao Wang, and Yizhe Zhu.
\newblock Universality of kernel random matrices and kernel regression in the quadratic regime.
\newblock {\em arXiv preprint arXiv:2408.01062}, 2024.

\bibitem{misiakiewicz2024non}
Theodor Misiakiewicz and Basil Saeed.
\newblock A non-asymptotic theory of kernel ridge regression: deterministic equivalents, test error, and gcv estimator.
\newblock {\em arXiv preprint arXiv:2403.08938}, 2024.

\bibitem{wang2023overparameterized}
Zhichao Wang and Yizhe Zhu.
\newblock Overparameterized random feature regression with nearly orthogonal data.
\newblock In {\em International Conference on Artificial Intelligence and Statistics}, pages 8463--8493. PMLR, 2023.

\bibitem{wang2024deformed}
Zhichao Wang and Yizhe Zhu.
\newblock Deformed semicircle law and concentration of nonlinear random matrices for ultra-wide neural networks.
\newblock {\em The Annals of Applied Probability}, 34(2):1896--1947, 2024.

\bibitem{Montanari_interpolation_2021}
Andrea Montanari and Yiqiao Zhong.
\newblock {The interpolation phase transition in neural networks: Memorization and generalization under lazy training}.
\newblock {\em The Annals of Statistics}, 50(5):2816 -- 2847, 2022.

\bibitem{Donhauser_how_2021}
Konstantin Donhauser, Mingqi Wu, and Fanny Yang.
\newblock How rotational invariance of common kernels prevents generalization in high dimensions.
\newblock In {\em International Conference on Machine Learning}, pages 2804--2814. PMLR, 2021.

\bibitem{xiao2022precise}
Lechao Xiao, Hong Hu, Theodor Misiakiewicz, Yue~M Lu, and Jeffrey Pennington.
\newblock Precise learning curves and higher-order scaling limits for dot product kernel regression.
\newblock {\em Journal of Statistical Mechanics: Theory and Experiment}, 2023(11):114005, 2023.

\bibitem{misiakiewicz_spectrum_2022}
Theodor Misiakiewicz.
\newblock Spectrum of inner-product kernel matrices in the polynomial regime and multiple descent phenomenon in kernel ridge regression.
\newblock {\em arXiv preprint arXiv:2204.10425}, 2022.

\bibitem{hu2022sharp}
Hong Hu and Yue~M Lu.
\newblock Sharp asymptotics of kernel ridge regression beyond the linear regime.
\newblock {\em arXiv preprint arXiv:2205.06798}, 2022.

\bibitem{steinwart2012mercer}
Ingo Steinwart and Clint Scovel.
\newblock Mercer’s theorem on general domains: On the interaction between measures, kernels, and rkhss.
\newblock {\em Constructive Approximation}, 35:363--417, 2012.

\bibitem{fischer2020_SobolevNorm}
Simon Fischer and Ingo Steinwart.
\newblock Sobolev norm learning rates for regularized least-squares algorithms.
\newblock {\em Journal of Machine Learning Research}, 21(205):1--38, 2020.

\bibitem{gneiting2013strictly}
Tilmann Gneiting.
\newblock {Strictly and non-strictly positive definite functions on spheres}.
\newblock {\em Bernoulli}, 19(4):1327 -- 1349, 2013.

\bibitem{Bietti_on_2019}
Alberto Bietti and Julien Mairal.
\newblock On the inductive bias of neural tangent kernels.
\newblock {\em Advances in Neural Information Processing Systems}, 32, 2019.

\bibitem{gallier2009notes}
Jean Gallier.
\newblock Notes on spherical harmonics and linear representations of lie groups.
\newblock Preprint, [OL], 2009.

\bibitem{https://doi.org/10.1002/cpa.22008}
Song Mei and Andrea Montanari.
\newblock The generalization error of random features regression: Precise asymptotics and the double descent curve.
\newblock {\em Communications on Pure and Applied Mathematics}, 75(4):667--766, 2022.

\bibitem{edmunds1996function}
David~Eric Edmunds and Hans Triebel.
\newblock {\em Function Spaces, Entropy Numbers, Differential Operators}.
\newblock Cambridge University Press, Cambridge, 1996.

\bibitem{Yao2007OnES}
Y.~Yao, Lorenzo Rosasco, and Andrea Caponnetto.
\newblock On early stopping in gradient descent learning.
\newblock {\em Constructive Approximation}, 26:289--315, 2007.

\bibitem{zhang2023optimality_2}
Haobo Zhang, Yicheng Li, and Qian Lin.
\newblock On the optimality of misspecified spectral algorithms.
\newblock {\em Journal of Machine Learning Research}, 25(188):1--50, 2024.

\bibitem{bahadur1964fisher}
R~Raj Bahadur.
\newblock On fisher's bound for asymptotic variances.
\newblock {\em The Annals of Mathematical Statistics}, 35(4):1545--1552, 1964.

\bibitem{tsybakov2008introduction}
A.B. Tsybakov.
\newblock {\em Introduction to Nonparametric Estimation}.
\newblock Springer Series in Statistics. Springer New York, 2008.

\bibitem{carter2006continous}
Andrew~V. Carter.
\newblock {A continuous Gaussian approximation to a nonparametric regression in two dimensions}.
\newblock {\em Bernoulli}, 12(1):143 -- 156, 2006.

\bibitem{reiss2008asymptotic}
Markus Rei{\ss}.
\newblock Asymptotic equivalence for nonparametric regression with multivariate and random design.
\newblock {\em The Annals of Statistics}, pages 1957--1982, 2008.

\bibitem{efromovich1996asymptotic}
Sam Efromovich and Alex Samarov.
\newblock Asymptotic equivalence of nonparametric regression and white noise model has its limits.
\newblock {\em Statistics \& probability letters}, 28(2):143--145, 1996.

\bibitem{brown1998asymptotic}
Lawrence~D Brown and Cun-Hui Zhang.
\newblock Asymptotic nonequivalence of nonparametric experiments when the smoothness index is 1/2.
\newblock {\em Annals of statistics}, pages 279--287, 1998.

\bibitem{cheng2022dimension}
Chen Cheng and Andrea Montanari.
\newblock Dimension free ridge regression.
\newblock {\em arXiv preprint arXiv:2210.08571}, 2022.

\bibitem{azevedo2015eigenvalues}
Douglas Azevedo and Valdir~A Menegatto.
\newblock Eigenvalues of dot-product kernels on the sphere.
\newblock {\em Proceeding Series of the Brazilian Society of Computational and Applied Mathematics}, 3(1), 2015.

\bibitem{jing2024advanced}
Bing-Yi Jing.
\newblock {\em Advanced Probability Theory}.
\newblock 2012.
\newblock Unpublished textbook.

\bibitem{shao2003mathematical}
Jun Shao.
\newblock {\em Mathematical statistics}.
\newblock Springer Science \& Business Media, 2003.

\end{thebibliography}

\appendix

\newpage

{\colorblue
\section{Notation Table}

Various statistical quantities are used in our proof to determine the Pinsker bound. Most of these notations are borrowed from \cite{tsybakov2008introduction}, ensuring consistency with established literature.

For readers' convenience, we provide the following Notation Table, listing all quantities used in the proof, their meaning, and the pages where they first appear.

\begin{table}[ht]
  \centering
  \caption{Notation Table}
  \label{tab:notation}
  \renewcommand{\arraystretch}{1.5}
  \begin{tabular}{c p{6cm} c}
    \toprule
    Symbol & Description & First Occurrence Page \\
    \midrule
    $\calP$               & a set of distributions on $\mathcal{X} \times \mathcal{Y}$               & \pageref{thm:main_pinsker_constant}  \\
    $\calC^{\star}$               & Pinsker constant               & \pageref{thm:main_pinsker_constant}  \\
    $\zeta$               & minimax rate               & \pageref{thm:main_pinsker_constant}  \\
    $\kappa^{\star}$               & defined in Definition \ref{def:kappa_and_N}               & \pageref{def:kappa_and_N}  \\
    $N$               & defined in Definition \ref{def:kappa_and_N}               & \pageref{def:kappa_and_N}  \\
    $\ell_j$               & defined in Definition \ref{def:pinsker_quantity}               & \pageref{def:pinsker_quantity}  \\
    $\calD^{\star}$               & Pinsker bound               & \pageref{def:pinsker_quantity}  \\
    $\hat{f}_{\ell}$               & linear filter estimator               & \pageref{eqn:upper_bound_copy}  \\
    $\Delta_n(j, j^{\prime})$               & defined in (\ref{def:bar_z_j})               & \pageref{def:bar_z_j}  \\
    $\Theta_N$               & a subset of $\mathbb{R}^{N}$               & \pageref{def:theta_spca_lower_bound}  \\
    $\calF_N$               & a function space associated with $\Theta_N$               & \pageref{def:theta_spca_lower_bound}  \\
    $\tilde{\mathcal{P}}$               & a subset of $\calP$               & \pageref{def:theta_spca_lower_bound}  \\
    $v_j^2$               & defined in (\ref{def:v_j_and_s_j})              & \pageref{def:v_j_and_s_j}  \\
    $s_j^{2}$               & defined in (\ref{def:v_j_and_s_j})              & \pageref{def:v_j_and_s_j}  \\
    $\mu_s(\cdot)$               & the p.d.f. of \(\mathcal{N}(0, s^2)\)              & \pageref{def:v_j_and_s_j}  \\
    $\mu(\cdot)$               & the p.d.f. of \(\calN \left(\mathbf 0, \text{diag} \left(s_1^2, \ldots, s_N^2\right) \right)\)              & \pageref{def:v_j_and_s_j}  \\
    \bottomrule
  \end{tabular}
  \renewcommand{\arraystretch}{1.0}
\end{table}

}

\newpage
\section{Proof of results in Section \ref{sec:pinsker_quantity}}

\subsection{Proof of Lemma \ref{lemma:inner_edr}}

\begin{proof}
The equation (22) in \cite{ghorbani2021linearized} holds for data uniformly distributed on $\sqrt{d} \mathbb{S}^d$, However, the spectrum estimates in \cite{ghorbani2021linearized} are invariant with respect to this scaling. Hence,  for any $k \geq 0$, we have 
\begin{equation}\label{eqn:mu_order_in_proof}
\mu_{k} = d^{-k}(\Phi^{(k)}(0)+o_d(d^{-1}))= d^{-k}(a_{k} k!+o_d(d^{-1})).
\end{equation}

For any $0 \leq k \leq p+3$, it is clear that
\begin{equation}\label{eqn:224_n_def}
    \begin{aligned}
    N(d, k)
     &= \frac{2k+d-1}{k(k+d-1)} \cdot \frac{(k+d-1)!}{(d-1)!(k-1)!}=\frac{d^{k}}{k!}(1+O_{d}(d^{-1})).
    \end{aligned}
\end{equation}
Now we begin to proof the second part of Lemma \ref{lemma:inner_edr}.
Notice that, for any $k \geq 0$, from \cite{azevedo2015eigenvalues}, we have
\begin{equation}\label{eqn:monotone_of_eigen}
    \begin{aligned}
        \frac{\mu_{k+2}}{\mu_{k}}
        &=
        \frac{1}{4}\cdot \frac{\sum_{s=0}^\infty a_{2s+k+2} \frac{(2s+k+2)!}{(2s)!} \frac{\Gamma(s+\frac{1}{2})}{\Gamma(s+k+2+\frac{d+1}{2})}}
        {\sum_{s=0}^\infty a_{2s+k} \frac{(2s+k)!}{(2s)!} \frac{\Gamma(s+\frac{1}{2})}{\Gamma(s+k+\frac{d+1}{2})}}=
        \frac{1}{4}\cdot \frac{\sum_{s=1}^\infty a_{2s+k} \frac{(2s+k)!}{(2s-2)!} \frac{\Gamma(s-\frac{1}{2})}{\Gamma(s+k+1+\frac{d+1}{2})}}
        {\sum_{s=0}^\infty a_{2s+k} \frac{(2s+k)!}{(2s)!} \frac{\Gamma(s+\frac{1}{2})}{\Gamma(s+k+\frac{d+1}{2})}}\\
        &=
        \frac{\sum_{s=1}^\infty a_{2s+k} \frac{(2s+k)!}{(2s)!} \frac{\Gamma(s+\frac{1}{2})}{\Gamma(s+k+\frac{d+1}{2})}
        \cdot \frac{s}{s+k+\frac{d+1}{2}}}
        {\sum_{s=0}^\infty a_{2s+k} \frac{(2s+k)!}{(2s)!} \frac{\Gamma(s+\frac{1}{2})}{\Gamma(s+k+\frac{d+1}{2})}} \overset{\text{Assumption } \ref{assu:coef_of_inner_prod_kernel}}{\leq}1.
    \end{aligned}
\end{equation}
Furthermore, since $a_{p+2}>0$, similar to (\ref{eqn:monotone_of_eigen}), we have $\mu_{p+2}> \mu_{p+4}$. Therefore, from (\ref{eqn:mu_order_in_proof}), (\ref{eqn:monotone_of_eigen}), and the definition of $p=\lfloor \gamma / (s+1)\rfloor \leq \lfloor \gamma \rfloor$, 
    there exists a constant $\mathfrak{C}_1$ (only depends on $\gamma, a_0, \cdots, a_{\lfloor \gamma \rfloor+3}$), such that
    for any $d \geq \mathfrak{C}_1$, we have
    \begin{gather*}
  0.9 \cdot a_{k} k! d^{-k} \leq  \mu_{k} \leq  1.1 \cdot a_{k} k! d^{-k} \quad \text{and} \quad 0.9 \cdot \frac{d^{k}}{k!} \leq
N(d, k) \leq 1.1 \cdot \frac{d^{k}}{k!},\\
  \mu_0 > \mu_1 > \cdots > \mu_{p+1} > \mu_{p+2} > \max_{j \geq p+3} \mu_{j}.
\end{gather*}
Consequently, from (\ref{spherical_decomposition_of_inner_main}), for any $0 \leq k \leq p+2$, we have:
\begin{align*}
      \lambda_{v_{k-1}+1} = \lambda_{v_{k-1}+2}= \cdots=\lambda_{v_{k}} = \mu_k, ~ \{ 
         \phi_{v_{k-1}+1}, \phi_{v_{k-1}+2} \cdots, \phi_{v_{k}}\} =\{Y_{k,1},....,Y_{k,N(d,k)}\},
    \end{align*}
    finishing the proof.
\end{proof}

\subsection{Proof of Lemma \ref{lemma:calcul_N}}

\begin{proof}
    From Lemma \ref{lemma:inner_edr}, 
    there exists a constant $\mathfrak{C}_1$, depending only on the absolute constants $\gamma, a_0, \cdots, a_{\lfloor \gamma \rfloor+3}$, such that
    for any $d \geq \mathfrak{C}_1$, we have
    \begin{equation}\label{eqn:eqn_lemma_4_5_prep_1}
        \mu_0 > \mu_1 > \cdots > \mu_{p+2} > \max_{j \geq p+3}\mu_{j}.
    \end{equation}

%We proceed by contradiction to show that $ \mu_{p+2}^{s/2}\leq \kappa^{\star} <\mu_p^{s/2} $. 
To proceed, we will demonstrate that any of the following four cases leads to a contradiction: (i) $\mu_p^{s/2} \leq \kappa^{\star}$, (ii) $\mu_{p+2}^{s/2} > \kappa^{\star}$, (iii) $\gamma < p(s+1) + s/2$ and $N=\sum_{k=0}^{p+1}N(d, k)$, or (iv) $\gamma > p(s+1) + s/2$ and $N=\sum_{k=0}^{p}N(d, k)$.
These will establish that $\mu_{p+2}^{s/2}\leq \kappa^{\star} <\mu_p^{s/2}$, implying that $\ell_{j}=0$ for any $j\geq v_{p+1}= \sum_{k=0}^{p+1}N(d,k)+1$ and $\ell_{j}\neq 0$ for any $j\leq v_{p}=\sum_{k=0}^{p}N(d,k)$. Therefore:
$$
    N = \sum_{k=0}^{p}N(d, k) \quad \text{ or } \quad N=\sum_{k=0}^{p+1}N(d, k);
$$
Moreover, when $\gamma <p(s+1) + s/2$, we have $q=p$; when $\gamma > p(s+1) + s/2$, we have $q=p+1$.

\noindent{\bf Case (i):} If $\mu_p^{s/2} \leq \kappa^{\star}$, then $\ell_{j}=0$ for any $j\geq v_{p-1}=\sum_{k=0}^{p-1}N(d,k)+1$. Therefore, 
\begin{equation}\label{eqn:lemma_calcul_N_i}
    \begin{aligned}
        R \overset{(\ref{def:kappa})}{=}
        &
        \frac{\sigma^2}{n \kappa^{\star}} \sum_{k=0}^{p-1} N(d, k) \mu_k^{-s/2} \left( 1 - \kappa^{\star} \mu_k^{-s/2} \right)_{+}
        \leq 
        \frac{\sigma^2}{n \mu_p^{s/2}} \sum_{k=0}^{p-1} \mu_k^{-s/2} N(d, k)\\
        \sim
        &
        \frac{\sigma^2}{\alpha d^{\gamma} \left( a_p \right)^{s/2} (p!)^{s/2} d^{-sp/2}} \sum_{k=0}^{p-1} (a_k)^{-\frac{s}{2}} d^{\frac{sk}{2}} \frac{d^k}{(k!)^{s/2+1} }\\
        \sim 
        &
        \frac{\sigma^2}{\alpha (a_p)^{s/2} (a_{p-1})^{s/2}(p!)^{s/2} ((p-1)!)^{s/2+1} } d^{-\gamma+p(s+1)-s/2-1},
    \end{aligned}
\end{equation}
where the approximation in the second line follows from Assumption \ref{assump_asymptotic} and Lemma \ref{lemma:inner_edr}. 
Since $R$ is an absolute positive constant and $\gamma \geq p(s+1)$, when $d \geq \mathfrak{C}_2$ (a sufficiently large constant only depending on the absolute constants defined in Definition \ref{def:abs_constants}), we get a contradiction.

\noindent{\bf Case (ii)} If $\mu_{p+2}^{s/2} > \kappa^{\star}$, then for $d \geq \mathfrak{C}_1$, 
Lemma \ref{lemma:inner_edr} implies  $\kappa^{\star} \mu_{p+1}^{-s/2} < [\mu_{p+2} / \mu_{p+1}]^{s/2}<1$.
Therefore,
\begin{equation}\label{eqn:lemma_calcul_N_ii}
    \begin{aligned}
        R \overset{(\ref{def:kappa})}{\geq}
        &~
        \frac{\sigma^2}{n \kappa^{\star}} \sum_{k=0}^{p+1} N(d, k) \mu_k^{-s/2} \left( 1 - \kappa^{\star} \mu_k^{-s/2} \right)_{+}\\
        > 
        &~
        \frac{\sigma^2}{n (\mu_{p+2})^{s/2}} \sum_{k=0}^{p+1} \mu_k^{-s/2} N(d, k) - \frac{\sigma^2}{n} \sum_{k=0}^{p+1} \mu_k^{-s} N(d, k)\\
        \sim 
        &~
        \frac{\sigma^2}{\alpha d^{\gamma} \left( a_{p+2} \right)^{s/2} ((p+2)!)^{s/2} d^{-s(p+2)/2}} \sum_{k=0}^{p+1} (a_k)^{-\frac{s}{2}} d^{\frac{sk}{2}} \frac{d^k}{(k!)^{s/2+1} } +O_d(d^{-\gamma + (p+1)(s+1)})\\
        \sim 
        &~\frac{\sigma^2}{\alpha (a_{p+2})^{s/2} (a_{p+1})^{s/2} ((p+2)!)^{s/2}((p+1)!)^{s/2+1} } d^{-\gamma+(p+1)(s+1)+s/2}.
    \end{aligned}
\end{equation}
Since $R$ is an absolute positive constant and $\gamma < (p+1)(s+1)$, when $d \geq \mathfrak{C}_3$ (a sufficiently large constant only depending on the absolute constants defined in Definition \ref{def:abs_constants}), we also get a contradiction.

\noindent{\bf Case (iii)} If $\gamma < p(s+1) + s/2$ and $N=\sum_{k=0}^{p+1}N(d, k)$, then by the definition of $N$ we have $1 - \kappa^{\star}\mu_{p+1}^{-s/2}>0$. However, from (\ref{eqn:order_of_kappa}) we find
\begin{equation*}
\begin{aligned}
    &~1 - \kappa^{\star}\mu_{p+1}^{-s/2} = 1-\frac{\sigma^2\mu_{p+1}^{-s/2} \sum\limits_{k=0}^{{p+1}} \mu_k^{-s/2} N(d, k)}{n R +\sigma^2\sum\limits_{k=0}^{{p+1}} \mu_k^{-s} N(d, k)}\\
    =&~
    \frac{n R +\sigma^2\sum\limits_{k=0}^{{p}} \left(\mu_k^{-s} - \mu_{p+1}^{-s/2}\mu_k^{-s/2}\right) N(d, k)}{n R +\sigma^2\sum\limits_{k=0}^{{p+1}} \mu_k^{-s} N(d, k)}
    \sim 
    \frac{n R   -\sigma^2 \mu_{p+1}^{-s/2}\mu_p^{-s/2} N(d, p)}{n R +\sigma^2 \mu_{p+1}^{-s} N(d, p+1)}\\
    %=&~
    %\frac{n R + \sigma^2 \left(\mu_p^{-s} - \mu_{p+1}^{-s/2}\mu_p^{-s/2}\right) N(d, p)+\sigma^2\sum\limits_{k=0}^{{p-1}} \left(\mu_k^{-s} - \mu_{p+1}^{-s/2}\mu_k^{-s/2}\right) N(d, k)}{n R +\sigma^2\sum\limits_{k=0}^{{p+1}} \mu_k^{-s} N(d, k)}\\
    \sim &~
    \frac{\alpha R d^{\gamma}  -  \frac{\sigma^2 }{a_{p}^{s/2} a_{p+1}^{s/2} (p!)^{s/2+1} (({p+1})!)^{s/2} } d^{(s+1)p+s/2} }{\alpha R d^{\gamma} +   \frac{\sigma^2 }{ a_{p+1}^{s} (({p+1})!)^{s+1} }d^{({p+1})s + {p+1}}}.
\end{aligned}
\end{equation*}
Therefore, when $d \geq \mathfrak{C}_4$ (a sufficiently large constant only depending on the absolute constants defined in Definition \ref{def:abs_constants}), we get a contradiction that $1 - \kappa^{\star}\mu_{p+1}^{-s/2}<0$.

\noindent{\bf Case (iv)} If $\gamma > p(s+1) + s/2$ and $N=\sum_{k=0}^{p}N(d, k)$, then by the definition of of $N$ we have $1 - \kappa^{\star}\mu_{p+1}^{-s/2} \leq 0$. 
However, similar to (iii), for $d \geq \mathfrak{C}_5$ (a sufficiently large constant only depending on the absolute constants defined in Definition \ref{def:abs_constants}), from (\ref{eqn:order_of_kappa}) we get a contradiction that $1 - \kappa^{\star}\mu_{p+1}^{-s/2}>0$.

Combining the results from cases (i) through (iv), we define $\mathfrak{C} = \max\{\mathfrak{C}_1, \mathfrak{C}_2, \mathfrak{C}_3, \mathfrak{C}_4, \mathfrak{C}_5\}$. With this definition, we obtain the desired results.
\end{proof}

\subsection{Proof of Corollary \ref{coroll:quantities}}
When \( d \geq \mathfrak{C} \),
Lemma \ref{lemma:calcul_N} implies $N=\sum_{k=0}^{q}N(d,k)$ for $q=p$ or $q=p+1$. Hence, we only need to show that $\mathcal{D}^{\star} \sim \mathcal{C}^{\star} d^{-\zeta}$ in the following two situations: 
%(1) $(s+1)p\leq \gamma<(s+1)p+s/2$, and (2) $(s+1)p+s/2 \leq \gamma < (s+1)(p+1)$.

\subsubsection{}\label{lem:quantities_1} %Case 1: $(s+1)p\leq \gamma<(s+1)p+s/2$} 
When $q=p$ and $N = \sum_{k=0}^{p}N(d, k)$, by Lemma \ref{lemma:calcul_N}, we know that $\gamma \leq  p(s+1)+s/2$.  
We will prove the Corollary \ref{coroll:quantities} in the following two steps.

\begin{itemize}

\item[(i)]  $\kappa^{\star} \sim \frac{\sigma^2 a_p^{s/2} (p!)^{s/2} }{\alpha R a_p^{s} (p!)^{s+1}   +   \sigma^2 \mathbf{1}\{\gamma  = ps + p\} } d^{ps/2 + p - \gamma} $;\\[1em]

 From (\ref{eqn:order_of_kappa}) we have
\begin{equation}\label{eqn:quantities_1_bound_ii}
\begin{aligned}
    \kappa^{\star} = &~\frac{\sigma^2\sum\limits_{k=0}^{p} \mu_k^{-s/2} N(d, k)}{n R +\sigma^2\sum\limits_{k=0}^{p} \mu_k^{-s} N(d, k)}
    \sim
    \frac{  \frac{\sigma^2 }{a_p^{s/2}  (p!)^{s/2+1} } d^{ps/2 + p} }{\alpha R d^{\gamma} +   \frac{\sigma^2 }{ a_p^{s} (p!)^{s+1} }d^{ps + p}}\\
    \sim &~\begin{cases}
    \frac{  \sigma^2 a_p^{s/2} (p!)^{s/2} }{\alpha R a_p^{s} (p!)^{s+1}   +   \sigma^2} d^{-ps/2}  & \mbox{ if } \gamma=ps+p\\
    \frac{\sigma^2 }{\alpha R a_p^{s/2}  (p!)^{s/2+1} } d^{ps/2 + p - \gamma}  & \mbox{ if } \gamma>ps+p
    \end{cases}.
\end{aligned}
\end{equation}\\

\item[(ii)] $
\calD^* \sim \frac{ \sigma^2   }{\alpha  p!  + \sigma^2 / (R a_p^{s} (p!)^{s})\mathbf{1}\{\gamma=(s+1)p\}} d^{p-\gamma}.
$ \\[1em]

When $\gamma  = ps + p$, from Lemma \ref{lemma:inner_edr} and (\ref{eqn:quantities_1_bound_ii}), we have 
\begin{equation*}
    \begin{aligned}
        \calD^{\star}
        = &~
        \frac{\sigma^2}{n} \sum_{k=0}^{p} N(d, k) (1 - \kappa^{\star} \mu_k^{-s/2})_{+}\\
        \sim  &~
        \frac{\sigma^2}{n} N(d, p) 
        \frac{  \alpha R a_p^{s} (p!)^{s+1}  }{\alpha R a_p^{s} (p!)^{s+1}  + \sigma^2}
        \sim 
        \frac{ \sigma^2   }{\alpha  p!  + \sigma^2 / (R a_p^{s} (p!)^{s})} d^{p-\gamma};
    \end{aligned}
\end{equation*}
When $\gamma  > ps + p$, from Lemma \ref{lemma:inner_edr} and (\ref{eqn:quantities_1_bound_ii}), we have 
\begin{equation*}
    \begin{aligned}
        \calD^{\star}
        \sim
        \frac{\sigma^2}{n} \sum_{k=0}^{p} N(d, k) (1 - \kappa^{\star} \mu_p^{-s/2})_{+}
        \sim
        \frac{\sigma^2}{n} N(d, p)
        \sim
        \frac{ \sigma^2  }{\alpha  p!  } d^{p-\gamma},
    \end{aligned}
\end{equation*}
and we get the desired results.
\end{itemize}

\subsubsection{}\label{lem:quantities_2}
%{Proof of Corollary \ref{coroll:quantities} when $(s+1)p+s/2 \leq \gamma < (s+1)(p+1)$}
When $q=p+1$ and $N = \sum_{k=0}^{p+1}N(d, k)$, by Lemma \ref{lemma:calcul_N}, we know that $\gamma\geq p(s+1)+s/2$.
We will prove the Corollary \ref{coroll:quantities} in the following two steps.

\begin{itemize}

\item[(i)] $\kappa^{\star} = \Theta_d( d^{-({p+1})s/2})$.\\[1em]

%Recall that from Lemma \ref{lemma:calcul_N}, we have $\gamma \geq p(s+1)+s/2$. When $\gamma > p(s+1)+s/2$, 
If $p(s+1)+s/2<\gamma<(p+1)(s+1)$, then (\ref{eqn:order_of_kappa}) implies
\begin{equation}\label{eqn:quantities_2_bound_ii_origin}
\begin{aligned}
    &~1 - \kappa^{\star}\mu_{p+1}^{-s/2} = 1-\frac{\sigma^2\mu_{p+1}^{-s/2} \sum\limits_{k=0}^{{p+1}} \mu_k^{-s/2} N(d, k)}{n R +\sigma^2\sum\limits_{k=0}^{{p+1}} \mu_k^{-s} N(d, k)}\\
    =&~
    \frac{n R +\sigma^2\sum\limits_{k=0}^{{p}} \left(\mu_k^{-s} - \mu_{p+1}^{-s/2}\mu_k^{-s/2}\right) N(d, k)}{n R +\sigma^2\sum\limits_{k=0}^{{p+1}} \mu_k^{-s} N(d, k)}
    \sim 
    \frac{n R   -\sigma^2 \mu_{p+1}^{-s/2}\mu_p^{-s/2} N(d, p)}{n R +\sigma^2 \mu_{p+1}^{-s} N(d, p+1)}\\
    %=&~
    %\frac{n R + \sigma^2 \left(\mu_p^{-s} - \mu_{p+1}^{-s/2}\mu_p^{-s/2}\right) N(d, p)+\sigma^2\sum\limits_{k=0}^{{p-1}} \left(\mu_k^{-s} - \mu_{p+1}^{-s/2}\mu_k^{-s/2}\right) N(d, k)}{n R +\sigma^2\sum\limits_{k=0}^{{p+1}} \mu_k^{-s} N(d, k)}\\
    \sim &~
    \frac{\alpha R d^{\gamma}  -  \frac{\sigma^2 }{a_{p}^{s/2} a_{p+1}^{s/2} (p!)^{s/2+1} (({p+1})!)^{s/2} } d^{(s+1)p+s/2} }{\alpha R d^{\gamma} +   \frac{\sigma^2 }{ a_{p+1}^{s} (({p+1})!)^{s+1} }d^{({p+1})s + {p+1}}}\\
    \sim &~
    \frac{\alpha R d^{\gamma}}{ \frac{\sigma^2 }{ a_{p+1}^{s} (({p+1})!)^{s+1} }d^{({p+1})s + {p+1}}}
    =~\frac{\alpha R a_{p+1}^{s}((p+1)!)^{s+1}}{\sigma^{2}}d^{\gamma-(p+1)(s+1)},
\end{aligned}
\end{equation}
where the last line follows from $p(s+1)+s/2<\gamma<(p+1)(s+1)$.
Hence we have
\begin{equation*}
\begin{aligned}
\kappa^{\star} \sim \mu_{p+1}^{s/2} \sim a_{p+1}^{s/2} (({p+1})!)^{s/2} d^{-(p+1)s/2}.
\end{aligned}
\end{equation*}
If $\gamma = p(s+1)+s/2$, then $q=p+1$ implies $1 - \kappa^{\star}\mu_{p+1}^{-s/2}>0$. Hence, similar to (\ref{eqn:quantities_2_bound_ii}), we can show that $0<1 - \kappa^{\star}\mu_{p+1}^{-s/2} = O_d(d^{\gamma-(p+1)(s+1)})=o_d(1)$. Therefore, we have $\kappa^{\star} = \Theta_d( d^{-({p+1})s/2})$.

Combining all above, for any $\gamma \geq p(s+1)+s/2$, we have
\begin{equation}\label{eqn:quantities_2_bound_ii}
\begin{aligned}
0< 1 - \kappa^{\star}\mu_{p+1}^{-s/2} = O_d(d^{\gamma-(p+1)(s+1)}),
\end{aligned}
\end{equation}
and
\begin{equation}\label{eqn:quantities_2_bound_ii_coro}
\begin{aligned}
\kappa^{\star} = \Theta_d( d^{-({p+1})s/2})
\end{aligned}
\end{equation}

\item[(ii)]  $
\calD^{\star} \sim R {a_{p+1}^{s}} ((p+1)!)^{s} d^{-(p+1)s}
        +
        \frac{ \sigma^2  }{\alpha  p!  } d^{p-\gamma}
$.\\[1em]
If $p(s+1)+s/2<\gamma<(p+1)(s+1)$, then from Lemma \ref{lemma:inner_edr}, (\ref{eqn:quantities_2_bound_ii_origin}), and (\ref{eqn:quantities_2_bound_ii_coro}), we have 
\begin{equation*}
    \begin{aligned}
        \calD^{\star}
        &\sim
        \frac{\sigma^2}{n} \sum_{k=0}^{p+1} N(d, k) (1 - \kappa^{\star} \mu_k^{-s/2})_{+}\\
        &\sim \frac{\sigma^2}{n}  N(d, p+1) (1 - \kappa^{\star} \mu_{p+1}^{-s/2})_{+}+\frac{\sigma^2}{n}  N(d, p) (1 - \kappa^{\star} \mu_p^{-s/2})_{+}\\
        &\sim
        \frac{\sigma^2}{n} N(d, p+1) \frac{\alpha R {a_{p+1}^{s}  (({p+1})!)^{s+1} }}{\sigma^2} d^{\gamma-(p+1)(s+1)} + \frac{\sigma^2}{n} N(d, p)\\
        &\sim
         R {a_{p+1}^{s}} (({p+1})!)^{s} d^{-(p+1)s}
        +
        \frac{ \sigma^2  }{\alpha  p! } d^{p-\gamma}.
    \end{aligned}
\end{equation*}
If $\gamma = p(s+1)+s/2$, then similarly, we have
\begin{equation*}
    \begin{aligned}
        \calD^{\star}
        &\sim
        \frac{\sigma^2}{n} \sum_{k=0}^{p+1} N(d, k) (1 - \kappa^{\star} \mu_k^{-s/2})_{+}\\
        &\sim \frac{\sigma^2}{n}  N(d, p+1) (1 - \kappa^{\star} \mu_{p+1}^{-s/2})_{+}+\frac{\sigma^2}{n}  N(d, p) (1 - \kappa^{\star} \mu_p^{-s/2})_{+}\\
        &\sim
        \frac{\sigma^2}{n} N(d, p+1) O_d(d^{-s/2-1}) + \frac{\sigma^2}{n} N(d, p)
        \sim
        \frac{ \sigma^2  }{\alpha  p! } d^{p-\gamma}.
    \end{aligned}
\end{equation*}

\end{itemize}
\iffalse
Now, consider the case $\gamma = (s+1)p+s/2$.
Denote $N = \sum_{k=0}^{q}N(d, k)$ and $\beta = \alpha R -  \frac{\sigma^2 }{a_{p}^{s/2} a_{p+1}^{s/2} (p!)^{s/2+1} (({p+1})!)^{s/2} }$. If $\beta \leq 0$, then similar to (\ref{eqn:quantities_2_bound_ii}), we have $1 - \kappa^{\star}\mu_{p+1}^{-s/2} =o_d(d^{\gamma-(p+1)(s+1)})$, hence
$$
\kappa^{\star} = \Omega_d(d^{-(p+1)s/2}).
$$
On the other hand, if $\beta > 0$, then similar to (\ref{eqn:quantities_2_bound_ii}), we have
$$
1 - \kappa^{\star}\mu_{p+1}^{-s/2} \sim 
    \frac{\alpha R {a_{p+1}^{s}  (({p+1})!)^{s+1} }}{\sigma^2} d^{\gamma-(p+1)(s+1)}
$$
To sum up, for any value of $\beta$, similar to the proof of Lemma \ref{lem:quantities_1} and the above discussion, we have
\fi

%%%%%%%%%%%%%%%%%%%%%%
Before we conclude this section, we present a proposition that will be useful in establishing the lower bound on the minimax risk.

\begin{proposition}\label{prop_calcul_max_ell}
Suppose Assumptions \ref{assump_asymptotic} and \ref{assu:coef_of_inner_prod_kernel} hold for some \(\alpha, \gamma > 0\). Further, suppose $\gamma \geq s$. Then, when \(d \geq \mathfrak{C}\), where \(\mathfrak{C}\) is the constant defined in Lemma \ref{lemma:calcul_N}, we have 
    $$
    \max_{1\leq j\leq N} \frac{\ell_j}{n\lambda_j^{s/2}\kappa^{\star}} = O_d(d^{-\min\{1, \gamma-s/2\}}).
    $$
\end{proposition}
    \begin{proof}
When \( d \geq \mathfrak{C} \),
Lemma \ref{lemma:calcul_N} implies that $N=\sum_{k=0}^{q}N(d,k)$ for $q=p \geq 1$ or $q=p+1\geq 2$, and that $q=1$ when $p=0$.
We therefore need to prove two main cases:
\begin{itemize}
    \item[(i)] If \( q = p \geq 1 \), then $\max_{1\leq j\leq N} \frac{\ell_j}{n\lambda_j^{s/2}\kappa^{\star}} = O_d(d^{-p})$;
    \item[(ii)] If \( q = p+1 \), then $\max_{1\leq j\leq N} \frac{\ell_j}{n\lambda_j^{s/2}\kappa^{\star}} =  O_d(d^{-\gamma+ps+s/2} + d^{-p-1})$.
\end{itemize}
Then when $\gamma \geq s$, these will establish that
\begin{align*}
    &~ \max_{1\leq j\leq N} \frac{\ell_j}{n\lambda_j^{s/2}\kappa^{\star}}\\
    =
    &~
    O_d(d^{-p})\mathbf{1}\{q=p \geq 1\}
    + O_d(d^{-\gamma+ps+s/2} + d^{-p-1})\mathbf{1}\{q=p+1 \geq 2\}\\
    &~+ O_d(d^{-\gamma+s/2} + d^{-1})\mathbf{1}\{p=0\}\\
    =&~
    O_d(d^{-1}) + O_d(d^{-1-s/2} + d^{-2}) + O_d(d^{-\gamma+s/2} + d^{-1})\\
    =&~ O_d(d^{-\min\{1, \gamma-s/2\}}).
\end{align*}

\noindent {\bf Case (i): } If \( q = p \geq 1 \), since $\ell_j = (1-\kappa^{\star} \lambda_j^{-s/2})_{+}$, we have
\begin{align*}
    \max_{1\leq j\leq N} \frac{\ell_j}{n\lambda_j^{s/2}\kappa^{\star}}
    =
    \max_{k \leq p}
    \frac{1-\kappa^{\star} \mu_k^{-s/2}}{n\mu_k^{s/2}\kappa^{\star}}
    \leq
    \max_{k \leq p}
    \frac{1}{n\mu_k^{s/2}\kappa^{\star}}
    =
    \frac{1}{n\mu_p^{s/2}\kappa^{\star}}.
\end{align*}
From the bounds in (\ref{eqn:quantities_1_bound_ii}) and Lemma \ref{lemma:inner_edr}, we have $n\mu_p^{s/2}\kappa^{\star} = \Theta_d(d^p)$.
\iffalse
\begin{equation}\label{eqn:quantities_1_bound_iii}
    %\kappa^{\star} \mu_p^{-s/2} = \frac{  \sigma^2 \mathbf{1}\{\gamma  = ps + p\} }{\alpha R a_p^{s} (p!)^{s+1}  + \sigma^2}  + o_d(1) \mathbf{1}\{\gamma  > ps + p\} ~ \text{ and } ~ 
    n\mu_p^{s/2}\kappa^{\star} = \Theta_d(d^p),
\end{equation}
\fi
Thus,
$$
\max_{1\leq j\leq N} \frac{\ell_j}{n\lambda_j^{s/2}\kappa^{\star}} = O_d(d^{-p}).
$$\\

\noindent {\bf Case (ii): } If \( q = p +1 \), using (\ref{eqn:quantities_2_bound_ii}) and $\ell_j = (1-\kappa^{\star} \lambda_j^{-s/2})_{+}$, we have
\begin{align*}
    \max_{1\leq j\leq N} \frac{\ell_j}{n\lambda_j^{s/2}\kappa^{\star}}
    =
    \max_{k \leq p+1}
    \frac{1-\kappa^{\star} \mu_k^{-s/2}}{n\mu_k^{s/2}\kappa^{\star}}
    \leq
    \max\left\{
    \frac{1}{n\mu_p^{s/2}\kappa^{\star}}, 
    \frac{O_d(d^{\gamma-(p+1)(s+1)})}{n\mu_{p+1}^{s/2}\kappa^{\star}}
    \right\}.
\end{align*}
From (\ref{eqn:quantities_2_bound_ii_coro}) and Lemma \ref{lemma:inner_edr}, we have
\begin{equation}\label{eqn:quantities_2_bound_iii}
    %\kappa^{\star} \mu_p^{-s/2} =  o_d(1) 
    n\mu_{p}^{s/2}\kappa^{\star} = \Omega_d(d^{\gamma-ps-s/2}) \quad \text{ and } \quad n\mu_{p+1}^{s/2}\kappa^{\star} = \Omega_d(d^{\gamma - (p+1)s}),
\end{equation}
thus,
$$
\max_{1\leq j\leq N} \frac{\ell_j}{n\lambda_j^{s/2}\kappa^{\star}} = O_d(d^{-\gamma+ps+s/2} + d^{-p-1}).
$$
\end{proof}

\section{Proof of upper bound in Theorem \ref{thm:main_pinsker_constant}}\label{sec:upper}

In this section, our goal is to show that
\begin{equation}\label{eqn:upper_bound}
        \inf_{\hat{f}} \sup_{
        \rho_{f_{\star}} \in \mathcal{P}
        %f_{\star} \in \sqrt{R}[\calB]^{s}
        } \mathbb{E}_{(X, Y) \overset{\mathcal{D}}{\sim} \rho_{f_{\star}}^{\otimes n}}\left[
    \|\hat{f}-f_{\star}\|_{L^2}^2
    \right]
     \leq
    \mathcal{D}^{\star}(1+o_d(1)).
    \end{equation}

For notation simplicity, we denote $\mathbb{E} = \mathbb{E}_{(X, Y) \overset{\mathcal{D}}{\sim} \rho_{f_{\star}}^{\otimes n}}$, where the distributions $\rho_{f_{\star}}$ on $\mathcal{X} \times \mathcal{Y}$ is given by (\ref{equation:true_model}) such that Assumption 
\ref{assump_asymptotic}, 
\ref{assu:coef_of_inner_prod_kernel},  and \ref{assump_function_calss} hold for some $\alpha, \gamma>0$.

\subsection{Regression function with zero expectation}\label{append_zero_mean_case_1}

In this subsection, we consider regression functions in $\sqrt{R}[\calB]^{s}$ and have zero expectation, that is, we assume that
\begin{equation}\label{assumption_zero_mean}
    f_{\star}(\cdot) = \sum_j \theta_j \phi_j(\cdot) \in \sqrt{R}[\calB]^{s} \quad \text{ and } \quad \theta_{1}=\mathbb{E}_x f_{\star}(x) := \int f_{\star}(x) \rho_{\calX}(x) ~\mathsf{d}x = 0.
\end{equation}

For any $j \leq N$, denote
\begin{equation}\label{def:bar_z_j}
    \begin{aligned}
    \bar{z}_j &= \frac{1}{n} \sum_{i=1}^{n} y_i \phi_j(x_i)
    =
    \frac{1}{n}\sum_{i=1}^{n}  f_{\star}(x_i) \phi_j(x_i)
    +
    \frac{1}{n}\sum_{i=1}^{n}  \epsilon_i \phi_j(x_i)\\
    &=
    \sum_{j^{\prime}=1}^{\infty} \theta_{j^{\prime}} \left(\frac{1}{n}\sum_{i=1}^{n} \phi_{j^{\prime}}(x_i) \phi_j(x_i) \right)
    +
    \frac{1}{n}\sum_{i=1}^{n}  \epsilon_i \phi_j(x_i)\\
    &:=
    \theta_{j} + \sum_{j^{\prime}=1}^{\infty} \theta_{j^{\prime}} \Delta_n(j, j^{\prime})
    +
    \xi_j,
\end{aligned}
\end{equation}
where $\Delta_n(j, j^{\prime}) = \frac{1}{n}\sum_{i=1}^{n} \phi_{j^{\prime}}(x_i) \phi_j(x_i) - \delta_{j, j^{\prime}}$ and $\xi_j = \frac{1}{n}\sum_{i=1}^{n}  \epsilon_i \phi_j(x_i)$.

Let's construct an estimator of the regression function as
\begin{align*}
    \hat{f}_{\ell, 0}(x) := \ell_1 \bar{z}_1 \mathbf{1}\{p=0\} +
    \sum_{j=2}^{N} \ell_j \bar{z}_j \phi_j(x).
\end{align*}

Recall that from Lemma \ref{lemma:calcul_N} we have $N =\sum_{k=0}^{q}N(d, k)$ for $q=p$ or $q=p+1$. The following Theorem proves (\ref{eqn:upper_bound}) when $q=p$.

\begin{theorem}[Restate Theorem \ref{thm:upper_1_copy} when $q=p$]\label{thm:upper_1}
    Suppose the same conditions as Theorem \ref{thm:main_pinsker_constant}.
    Further, suppose that $N =\sum_{k=0}^{p}N(d, k)$. 
    %When $d \geq \mathfrak{C}$, a sufficiently large constant only depending on the constants defined in Definition \ref{def:abs_constants}.
    Then, for any $\varepsilon > 0$, there exist a constant $D_{\varepsilon}$ only depending on $\varepsilon$ and  ${\mathfrak{C}}$ defined in Lemma \ref{lemma:calcul_N}, such that for any $d > D_{\varepsilon}$, and
    for any regression function $f_{\star} \in \sqrt{R}[\calB]^{s}$ satisfying one of the following conditions: (i) $\mathbb{E}_x f_{\star}(x) = 0$ or (ii) $p=0$, we have
    \begin{equation*}
       \mathbb{E}\left[
    \|\hat{f}_{\ell, 0}-f_{\star}\|_{L^2}^2
    \right]
    \leq
    \mathcal{D}^{\star}(1+\varepsilon).
    \end{equation*}
\end{theorem}
\begin{proof}
    %Recall that for any $j \leq N$, we have $\ell_j = (1-\kappa^{\star} \lambda_j^{-s/2})$ given in Definition \ref{def:pinsker_quantity}. 
    If $p>0$, from Lemma \ref{lemma:calcul_N}, when $d \geq {\mathfrak{C}}$ (a sufficiently large constant defined in Lemma \ref{lemma:calcul_N}), we have $\phi_1=Y_{0, 1} \equiv 1$, hence $0 = \mathbb{E}_x f_{\star}(X) = \theta_1$. Therefore, for any $p \geq 0$, we have
    $$
        (\ell_1 \bar{z}_1 \mathbf{1}\{p=0\} - \theta_1 \mathbf{1}\{p=0\})^2 + \sum_{j=2}^\infty (\ell_j \bar{z}_j -\theta_j)^2 \leq \sum_{j=1}^\infty (\ell_j \bar{z}_j -\theta_j)^2.
    $$
    Moreover, $\xi_j \mid x_1, \cdots, x_n$
    are mutually independent zero-mean variables with variance no greater than
    $\frac{\sigma^2}{n^2}\sum_{i=1}^{n}  \phi_j^2(x_i)$. Hence, we have
    \begin{equation}\label{eqn:main_decomposotion_in_upper_bound}
\begin{aligned}
    &~
    \mathbb{E}\left[
    \|\hat{f}_{\ell, 0}-f_{\star}\|_{L^2}^2
    \mid x_1, \cdots, x_n\right]
    \leq
    \mathbb{E}\left[
    \sum_{j=1}^\infty (\ell_j \bar{z}_j -\theta_j)^2
    \mid x_1, \cdots, x_n\right]\\
    =&~
    \sum_{j=1}^\infty \mathbb{E}\left.\left[\left((\ell_j-1)\theta_j
    + \ell_j \sum_{j^{\prime}=1}^{\infty} \theta_{j^{\prime}} \Delta_n(j, j^{\prime})
    +
    \ell_j\xi_j\right)^2
    \right| x_1, \cdots, x_n\right]\\
    \leq &~
    \underbrace{\left[\sum_{j=1}^\infty (1-\ell_j)^2\theta_j^2 + \frac{\sigma^2}{n} \sum_{j=1}^\infty\ell_j^2\right]}_{\calD_{0}^{*}}
    +
    \underbrace{\sum_{j=1}^\infty\ell_j^2 \left(\sum_{j^{\prime}=1}^{\infty} \theta_{j^{\prime}} \Delta_n(j, j^{\prime})\right)^2}_{\mathbf{E}_{1}}\\
    &~+
    \underbrace{2 \sum_{j=1}^\infty(\ell_j-1)\theta_j \ell_j \sum_{j^{\prime}=1}^{\infty} \theta_{j^{\prime}} \Delta_n(j, j^{\prime})}_{\mathbf{E}_{2}}
    +
    \underbrace{\frac{\sigma^2}{n} \sum_{j=1}^\infty\ell_j^2 \left[ 
    \frac{1}{n}\sum_{i=1}^{n}\phi_j^2(x_i) - 1
    \right]}_{\mathbf{E}_{3}},
\end{aligned}
\end{equation}
where the second equation can be proven by applying the monotone convergence theorem to the sequence $\{\sum_{j=1}^{k} (\ell_j \bar{z}_j -\theta_j)^2\}_{k \geq 1}$.

We bound the above terms separately.

\subsubsection{Term $\mathcal{D}_{0}^{\star}$ }
From Lemma 3.2 in \cite{tsybakov2008introduction} we have 
\begin{align}\label{eqn:D0}
    \mathcal{D}_{0}^{\star} \leq \mathcal{D}^{\star}.
\end{align}
\begin{remark}
    For readers' convenience, we copy the proof for $\mathcal{D}_{0}^{\star} \leq \mathcal{D}^{\star}$ in \cite{tsybakov2008introduction} as follows. We have
\begin{equation*}
    \begin{aligned}
        \mathcal{D}_{0}^{\star} &= \sum_{j=1}^{\infty} \left( (1 - \ell_j)^2 \theta_j^2 + \frac{\sigma^2}{n} \ell_j^2 \right) 
        = \frac{\sigma^2}{n} \sum_{j=1}^{\infty} \ell_j^2 + \sum_{j=1}^{\infty} (1 - \ell_j)^2 \lambda_j^{s} \lambda_j^{-s} \theta_j^2 \\
        &\leq \frac{\sigma^2}{n} \sum_{j=1}^{\infty} \ell_j^2 + R \sup_{j \geq 1} \left[ (1 - \ell_j)^2 \lambda_j^{s} \right] \\
        &\leq \frac{\sigma^2}{n} \sum_{j=1}^{\infty} \ell_j^2 + R (\kappa^{\star})^2 \quad \text{(since $1 - \kappa^{\star} \lambda_j^{-s/2} \leq \ell_j \leq 1$)} \\
        &\equiv \frac{\sigma^2}{n} \sum_{j=1}^{\infty} \ell_j^2 + \frac{\sigma^2}{n} \kappa^{\star} \sum_{j=1}^{\infty} \lambda_j^{-s/2} \ell_j \quad \text{(by (\ref{def:kappa}))} \\
        &= \frac{\sigma^2}{n} \sum_{j=1}^{\infty} \ell_j (\ell_j + \kappa^{\star} \lambda_j^{-s/2}) 
        = \frac{\sigma^2}{n} \sum_{j=1}^{N} \ell_j (\ell_j + \kappa^{\star} \lambda_j^{-s/2}) = \frac{\sigma^2}{n} \sum_{j=1}^{N} \ell_j = \mathcal{D}^{\star}.
    \end{aligned}
\end{equation*}
\end{remark}

\subsubsection{Term $\mathbf{E}_1$}\label{append_bound_E_1}
Since $\ell_{j}=0$ for any $j>N$ and $\ell_{j}\leq 1$ for any $1\leq j\leq N$,  
we have
\begin{equation}\label{eqn:E1}
    \begin{aligned}
        \mathbf{E}_1 
    &=
    \sum_{j=1}^\infty\ell_j^2 \left(\sum_{j^{\prime}=1}^{\infty} \theta_{j^{\prime}} \Delta_n(j, j^{\prime})\right)^2\\
    &\leq
    \sum_{j=1}^{N} \left(\sum_{j^{\prime}=1}^{\infty} \theta_{j^{\prime}} \Delta_n(j, j^{\prime})\right)^2\\
    &\leq
    \underbrace{2\sum_{j=1}^{N} \left(\sum_{j^{\prime}=1}^{N} \theta_{j^{\prime}} \Delta_n(j, j^{\prime})\right)^2}_{\mathbf{E}_{11}}
    +
    \underbrace{2\sum_{j=1}^{N}\left(\sum_{j^{\prime}=N+1}^{\infty} \theta_{j^{\prime}} \Delta_n(j, j^{\prime})\right)^2}_{\mathbf{E}_{12}}.
    %&=2\|\Delta_n(:N, :N) \theta(:N)\|_2^2 +2\sum_{j=1}^{N} (\sum_{j^{\prime}=N+1}^{\infty} \theta_{j^{\prime}} \Delta_n(j, j^{\prime}))^2\\
    %&\leq 2\|\Delta_n(:N, :N)\|_2^2 \|f\|_{L^2}^2 + 2\sum_{j=1}^{N} (\sum_{j^{\prime}=N+1}^{\infty} \theta_{j^{\prime}} \Delta_n(j, j^{\prime}))^2\\
    \end{aligned}
\end{equation}

For the first term, we have
\begin{equation}\label{eqn:E11}
\begin{aligned}
\mathbb{E}\mathbf{E}_{11}
&=
   2\mathbb{E}\sum_{j=1}^{N} \left(\sum_{j^{\prime}=1}^{N} \theta_{j^{\prime}} \Delta_n(j, j^{\prime})\right)^2\\
&= 2 \sum_{j=1}^{N} \sum_{j'=1}^{N} \theta_{j'}^2 \mathbb{E} \Delta_n(j,j')^2 + 2\sum_{j=1}^{N} \sum_{u \neq v}^{N} \theta_u \theta_v \mathbb{E} \left[\Delta_n(j,u) \Delta_n(j,v)\right].
    \end{aligned}
\end{equation}

For any $j \leq N$, $a \neq j$, and $b \neq j$, we have
\begin{equation}\label{eqn:E1_1}
\begin{aligned}
&~\mathbb{E} \Delta_n(j,a)^2 = \frac{1}{n^2} \mathbb{E} \left( \sum_{i=1}^{n} \phi_j(x_i)  \phi_a(x_i) \right)^2 \\
= &~ \frac{1}{n^2} \sum_{i=1}^{n} \mathbb{E} \left( \phi_j(x_i)^2 \phi_a(x_i)^2 \right) 
 + \frac{1}{n^2} \sum_{i \neq i'} \mathbb{E} \left( \phi_j(x_i) \phi_j(x_{i'}) \phi_a(x_i) \phi_a(x_{i'}) \right) \\
 = &~ \frac{1}{n^2} \sum_{i=1}^{n} \mathbb{E} \left( \phi_j(x_i)^2 \phi_a(x_i)^2 \right) 
 + \frac{1}{n^2} \sum_{i \neq i'} \mathbb{E} \left( \phi_j(x_i) \phi_a(x_i)\right)
 \mathbb{E} \left( \phi_j(x_{i'})  \phi_a(x_{i'}) \right)\\
=&~ \frac{1}{n^2} \sum_{i=1}^{n} \mathbb{E} \left( \phi_j(x_i)^2 \phi_a(x_i)^2 \right);
\end{aligned}\end{equation}
and
\begin{equation}\label{eqn:E1_2}
\begin{aligned}
&~\mathbb{E} \Delta_n(j,j)^2 = \mathbb{E} \left( \frac{1}{n} \sum_{i=1}^{n} \phi_j^2(x_i) - 1 \right)^2 \\
=&~ \left[\frac{1}{n^2} \sum_{i=1}^{n} \mathbb{E} \phi_j^4(x_i) + \frac{1}{n^2} \sum_{i \neq i'} \mathbb{E} \left(\phi_j^2(x_i) \phi_j^2(x_{i'})\right) \right]
 - \frac{2}{n} \sum_{i=1}^{n} \mathbb{E} \phi_j^2(x_i) + 1 \\
=&~ \left[\frac{1}{n^2} \sum_{i=1}^{n} \mathbb{E} \phi_j^4(x_i) + \frac{n-1}{n} \right] - 2 + 1 = \frac{1}{n^2} \sum_{i=1}^{n} \mathbb{E} \phi_j^4(x_i) - \frac{1}{n};
\end{aligned}\end{equation}
and
\begin{equation}\label{eqn:E1_3}
\begin{aligned}
&~ \mathbb{E} \left[\Delta_n(j,a) \Delta_n(j,b)\right] = \frac{1}{n^2} \mathbb{E} \left[\left( \sum_{i=1}^{n} \phi_j(x_i) \phi_a(x_i) \right) \left( \sum_{i'=1}^{n} \phi_j(x_{i'}) \phi_b(x_{i'}) \right)\right] \\
=&~ \frac{1}{n^2} \sum_{i=1}^{n} \mathbb{E} \left( \phi_j(x_i)^2 \phi_a(x_i) \phi_b(x_i) \right) 
 + \frac{1}{n^2} \sum_{i \neq i'} \mathbb{E} \left( \phi_j(x_i) \phi_j(x_{i'}) \phi_a(x_i) \phi_b(x_{i'}) \right) \\
=&~ \frac{1}{n^2} \sum_{i=1}^{n} \mathbb{E} \left( \phi_j^2(x_i) \phi_a(x_i) \phi_b(x_i) \right);
\end{aligned}\end{equation}
and
\begin{equation}\label{eqn:E1_4}
\begin{aligned}
\mathbb{E} \left[\Delta_n(j,j) \Delta_n(j,b)\right] &= \frac{1}{n^2} \mathbb{E} \left[\left( \sum_{i=1}^{n} \phi_j(x_i)^2 - 1 \right) \left( \sum_{i'=1}^{n} \phi_j(x_{i'}) \phi_b(x_{i'}) \right)\right] \\
&= \frac{1}{n^2} \sum_{i=1}^{n} \mathbb{E} \left( \phi_j^3(x_i) \phi_b(x_i) \right).
\end{aligned}\end{equation}

Combining (\ref{eqn:E1_1}) and (\ref{eqn:E1_2}) we have
\begin{equation}\label{eqn:E1_5}
\begin{aligned}
&~\sum_{j=1}^{N} \sum_{j^{\prime}=1}^{N} \theta_{j'}^2 \mathbb{E} \Delta_n(j,j')^2\\
=&~
\frac{1}{n^2}\sum_{j=1}^{N} \sum_{j^{\prime} \neq j}^{N} \theta_{j^{\prime}}^2 \sum_{i=1}^{n} \mathbb{E} \left[\phi_j^2(x_i) \phi_{j^{\prime}}^2(x_i)\right]
+
\frac{1}{n^2}\sum_{j=1}^{N} \theta_{j}^2 \sum_{i=1}^{n} \mathbb{E} \phi_j^4(x_i) -\frac{1}{n}\sum_{j=1}^{N}\theta_j^2\\
=&~
\frac{1}{n^2}\sum_{j=1}^{N} \left(\sum_{j^{\prime} \neq j}^{N} \theta_{j^{\prime}}^2 \sum_{i=1}^{n} \mathbb{E} \left[\phi_j^2(x_i) \phi_{j^{\prime}}^2(x_i)\right]
+
 \theta_{j}^2 \sum_{i=1}^{n} \mathbb{E} \left[\phi_j^2(x_i) \phi_{j}^2(x_i)\right]\right) -\frac{1}{n}\sum_{j=1}^{N}\theta_j^2\\
=&~\frac{1}{n^{2}}\sum_{j=1}^{N}\sum_{j^{\prime}=1}^{N}\theta_{j'}^{2}\sum_{i=1}^{n}\bbE\left[\phi_{j}^{2}(x_{i})\phi_{j^{\prime}}^{2}(x_{i})\right]-\frac{1}{n}\sum_{j=1}^{N}\theta_j^2\\
=&~
\frac{1}{n^2}
\sum_{j^{\prime}=1}^{N}
  \theta_{j^{\prime}}^2 \sum_{i=1}^{n}\mathbb{E} \left[ \left(\sum_{j=1}^{N}\phi_j^2(x_i)\right) \phi_{j^{\prime}}^2(x_i)\right]
  -
 \frac{1}{n}\sum_{j=1}^{N}\theta_j^2\\
 =&~
 \frac{1}{n^2}
\sum_{j^{\prime}=1}^{N}
  \theta_{j^{\prime}}^2 \sum_{i=1}^{n} \mathbb{E} \left[N \phi_{j^{\prime}}^2(x_i)\right]
  -
 \frac{1}{n}\sum_{j=1}^{N}\theta_j^2
 =
 \sum_{j=1}^{N} \theta_j^2 \cdot \frac{N-1}{n},
\end{aligned}\end{equation}
where in the fifth equation we use the Addition Formula $\sum_{j=1}^{N}\phi_j^2(x) = N$, $x \in \mathbb{S}^{d}$ (see, e.g., Proposition 1.18 in \cite{gallier2009notes}).

Combining (\ref{eqn:E1_3}) and (\ref{eqn:E1_4}), for any $u \neq v \geq 1$, we have
\begin{equation}\label{eqn:E1_6}
\begin{aligned}
&~\sum_{j=1}^{N} \mathbb{E} \left[\Delta_n(j,u) \Delta_n(j,v)\right] = 
\sum_{j=1}^{N}
\frac{1}{n^2} \sum_{i=1}^{n} \mathbb{E} \left[ \phi_j^2(x_i) \phi_u(x_i) \phi_v(x_i)\right]\\
=&~
\frac{1}{n^2} \sum_{i=1}^{n} \mathbb{E} \left[\sum_{j=1}^N  \phi_j^2(x_i) \phi_u(x_i) \phi_v(x_i)\right]
=\frac{1}{n^2} \sum_{i=1}^{n} \mathbb{E}\left[ N \phi_u(x_i) \phi_v(x_i)\right] = 0,
\end{aligned}\end{equation}
where in the third equation we use the Addition Formula again. 

Finally, from (\ref{eqn:E11}), (\ref{eqn:E1_5}), and (\ref{eqn:E1_6}) we have
\begin{equation}\label{eqn:E11_final_control}
    \mathbb{E} \mathbf{E}_{11} = 2 \sum_{j=1}^{N} \theta_j^2 \cdot \frac{N-1}{n}.
\end{equation}

Now we begin to calculate the second term in (\ref{eqn:E1}). We first recall an elementary result, and readers can refer to, e.g., page 67 in \cite{jing2024advanced}:
\begin{proposition}[Integration term by term]\label{prop_integration_term_by_term}
    If $\sum_{j^{\prime}=1}^{\infty} \mathbb{E}|Z_{j^{\prime}}| < \infty$, then
\[
\sum\nolimits_{j^{\prime}=1}^{\infty} |Z_{j^{\prime}}| < \infty, \, \text{a.s.}
\]
so that $\sum_{j^{\prime}=1}^{\infty} Z_{j^{\prime}}$ converges a.s., and
\[
\mathbb{E}\left( \sum_{j^{\prime}=1}^{\infty} Z_{j^{\prime}} \right) = \sum_{j^{\prime}=1}^{\infty} \mathbb{E} Z_{j^{\prime}}.
\]
\end{proposition}
\begin{proof}[Proof of Proposition \ref{prop_integration_term_by_term}]
Let $Y_n = \sum_{i=1}^{n} Z_i$ and $Y = \sum_{i=1}^{\infty} Z_i$. Define $X = \sum_{i=1}^{\infty} |Z_i|$. 
Notice that $Y_n$ converges almost surely to $Y$, and $|Y_n| \leq X$ almost surely.
Moreover, by the monotone convergence theorem, we have:
\[
\mathbb{E}X \leq \sum_{i=1}^{\infty} \mathbb{E}|Z_i| < \infty.
\]
Therefore, by the dominated convergence theorem, we obtain:
\[
\mathbb{E}\left( \sum_{j^{\prime}=1}^{\infty} Z_{j^{\prime}} \right) = \mathbb{E}Y
=
\lim_{n \to \infty} \mathbb{E}Y_n
=
\lim_{n \to \infty} \mathbb{E}\sum_{i=1}^{n} Z_i
=
\lim_{n \to \infty} \sum_{i=1}^{n}\mathbb{E} Z_i
=
\sum_{j^{\prime}=1}^{\infty} \mathbb{E} Z_{j^{\prime}},
\]
and this completes the proof.
\end{proof}

Define
$$
\mathbf{E}_{121} =
    \sum_{j=1}^{N} \sum_{j^{\prime}=N+1}^{\infty} \theta_{j^{\prime}}^2  \Delta_n(j,j')^2 
    \quad \text{ and } \quad
\mathbf{E}_{122} = 
    \sum_{j=1}^{N} \sum_{u \neq v \geq N+1}^{\infty} \theta_u \theta_v  \left[\Delta_n(j,u) \Delta_n(j,v)\right], 
$$
and let's use Proposition \ref{prop_integration_term_by_term} to calculate their expectations.

\noindent {\bf Term } $\mathbb{E}\mathbf{E}_{121}$. \quad 
For any $k \leq N$ and $j^{\prime}>N$, let $Z_{j^{\prime}, k} =  \theta_{j^{\prime}}^2  \Delta_n(k,j')^2$. 
%\leq \sum_{j=1}^{N} \theta_{j^{\prime}}^2  \Delta_n(j,j')^2$. 
It is clear that we have
\begin{align*}
    \sum_{j^{\prime}=N+1}^{\infty} \mathbb{E}|Z_{j^{\prime}, k}| 
    =&~
    \sum_{j^{\prime}=N+1}^{\infty} \mathbb{E} \theta_{j^{\prime}}^2  \Delta_n(k,j')^2
    \leq
    \sum_{j^{\prime}=N+1}^{\infty} \mathbb{E}\sum_{j=1}^{N} \theta_{j^{\prime}}^2  \Delta_n(j,j')^2\\
    \overset{(\ref{eqn:E1_1})}{=}
    &~
    \sum_{j^{\prime}=N+1}^{\infty} \theta_{j^{\prime}}^2
    \sum_{j=1}^{N} \frac{1}{n^2} \sum_{i=1}^{n} \mathbb{E} \left( \phi_j(x_i)^2 \phi_{j^{\prime}}(x_i)^2 \right)\\
    \overset{\text{Addition formula}}{=}&~
    2 \sum_{j^{\prime}=N+1}^{\infty} \theta_{j^{\prime}}^2
    \frac{1}{n^2} \sum_{i=1}^{n} \mathbb{E} \left( N \phi_{j^{\prime}}(x_i)^2 \right)
    =
    \frac{2N}{n} \sum_{j^{\prime}=N+1}^{\infty} \theta_{j^{\prime}}^2
    < \infty.
\end{align*}
Therefore, from Proposition \ref{prop_integration_term_by_term} we have
\begin{align*}
    &~\mathbb{E}\mathbf{E}_{121}
    =
    \mathbb{E}\sum_{j=1}^{N} \sum_{j^{\prime}=N+1}^{\infty} \theta_{j^{\prime}}^2  \Delta_n(j,j')^2
    =
    \sum_{j=1}^{N}
    \mathbb{E}\left( \sum_{j^{\prime}=N+1}^{\infty} Z_{j^{\prime}, j} \right) = \sum_{j=1}^{N}\sum_{j^{\prime}=N+1}^{\infty} \mathbb{E} Z_{j^{\prime}, j}\\
    =&~
    \sum_{j^{\prime}=N+1}^{\infty} \sum_{j=1}^{N}\mathbb{E} Z_{j^{\prime}, j}
    =
    \sum_{j^{\prime}=N+1}^{\infty} \mathbb{E}\sum_{j=1}^{N} \theta_{j^{\prime}}^2  \Delta_n(j,j')^2
    =
    \frac{2N}{n} \sum_{j^{\prime}=N+1}^{\infty} \theta_{j^{\prime}}^2.
\end{align*}

\noindent {\bf Term } $\mathbb{E}\mathbf{E}_{122}$. \quad 
For any $k \leq N$ and $u, v \geq N+1$, let $Z_{u, v, k} =  \theta_{u} \theta_{v} \Delta_n(k, u) \Delta_n(k, v)$. We have
\begin{align*}
    &~
    \sum_{u \neq v \geq N+1}^{\infty} \mathbb{E}|Z_{u, v, k}| 
    =
    \sum_{u \neq v \geq N+1}^{\infty} \mathbb{E}|\theta_{u} \theta_{v} \Delta_n(k, u) \Delta_n(k, v)|\\
    \overset{\text{Cauchy–Schwarz inequality}}{\leq} &~
    \sum_{u \neq v \geq N+1}^{\infty}
    \left(\mathbb{E}|\theta_{u}^2 \Delta_n(k, u)^2|\right)^{1/2}
    \left(\mathbb{E}|\theta_{v}^2 \Delta_n(k, v)^2|\right)^{1/2}\\
    \overset{\text{Cauchy–Schwarz inequality}}{\leq} &~
    \left(\sum_{u = N+1}^{\infty}
    \mathbb{E}|\theta_{u}^2 \Delta_n(k, u)^2|\right)^{1/2}
    \left(\sum_{v = N+1}^{\infty}
    \mathbb{E}|\theta_{v}^2 \Delta_n(k, v)^2|\right)^{1/2}\\
    =&~
    \sum_{j^{\prime} = N+1}^{\infty}
    \mathbb{E}|\theta_{j^{\prime}}^2 \Delta_n(k, j^{\prime})^2|=
    \sum_{j^{\prime}=N+1}^{\infty} \mathbb{E}|Z_{j^{\prime}, k}|<\infty;
\end{align*}
Therefore, from Proposition \ref{prop_integration_term_by_term} we have
\begin{align*}
    &~\mathbb{E}\mathbf{E}_{122} = 
    \mathbb{E}\sum_{j=1}^{N} \sum_{u \neq v \geq N+1}^{\infty} \theta_u \theta_v  \left[\Delta_n(j,u) \Delta_n(j,v)\right]
    =
    \sum_{j=1}^{N} \mathbb{E}\left(\sum_{u \neq v \geq N+1}^{\infty} Z_{u, v, j}\right)\\
    =&~
    \sum_{j=1}^{N} \sum_{u \neq v \geq N+1}^{\infty}\mathbb{E}Z_{u, v, j}
    =
    \sum_{u \neq v \geq N+1}^{\infty} \sum_{j=1}^{N}\mathbb{E}Z_{u, v, j}\\
    =&~
    \sum_{u \neq v \geq N+1}^{\infty} \sum_{j=1}^{N}
    \mathbb{E}\left[\theta_u \theta_v  \Delta_n(j,u) \Delta_n(j,v)\right]
    \overset{(\ref{eqn:E1_6})}{=}0.
\end{align*}

Combining all these, we have
\begin{equation}\label{eqn:E12}
\begin{aligned}
&~\mathbb{E}\mathbf{E}_{12}
=
2\mathbb{E}\sum_{j=1}^{N} \left(\sum_{j^{\prime}=N+1}^{\infty} \theta_{j^{\prime}} \Delta_n(j, j^{\prime})\right)^2\\
=&~
2\mathbb{E}\sum_{j=1}^{N} \sum_{j^{\prime}=N+1}^{\infty} \theta_{j^{\prime}}^2  \Delta_n(j,j')^2
+2\mathbb{E}\sum_{j=1}^{N} \sum_{u \neq v \geq N+1}^{\infty} \theta_u \theta_v  \left[\Delta_n(j,u) \Delta_n(j,v)\right]\\
=&~
2\mathbb{E}\mathbf{E}_{121} + 2\mathbb{E}\mathbf{E}_{122}
=
\frac{2N}{n} \sum_{j=N+1}^{\infty} \theta_{j}^2.
\end{aligned}\end{equation}

Now we begin to bound $\mathbb{E}\mathbf{E}_1$ in (\ref{eqn:main_decomposotion_in_upper_bound}). We separate the proof into the following two cases.

\noindent{\bf (i). } We first consider the case when $p>0$. Recall that when $d \geq \mathfrak{C}$, from Lemma \ref{lemma:calcul_N} we have
\begin{align*}
     \mu_0 &> \mu_1 > \cdots > \mu_{p+1} > \max_{j \geq p+2} \mu_{j},
   \end{align*}
hence we have $\lambda_1 = \mu_0$ and $\lambda_2 = \mu_1$.
Notice that we have
$\mathbb{E}_x f_{\star}(x) = \theta_1 = 0$. Therefore, 
\begin{equation}\label{eqn:use_bound_L2}
    \begin{aligned}
        \mathbb{E}\mathbf{E}_1 
    &\leq \mathbb{E}\mathbf{E}_{11}+\mathbb{E}\mathbf{E}_{12}
    \overset{(\ref{eqn:E11_final_control}) \text{ and } (\ref{eqn:E12})}{\leq}
    \frac{2N}{n}\sum_{j=2}^{\infty} \theta_{j}^2
    \leq
    \frac{2N}{n} \cdot \mu_1^{s} \sum_{j=2}^{\infty} \lambda_j^{-s}\theta_{j}^2
    \leq
    \frac{2N}{n} \cdot \mu_1^{s} R,
    \end{aligned}
\end{equation}
where the last inequality comes from the definition of interpolation space $[\calB]^{s}$ in Subsection \ref{subsec:interpolation_space}.

\noindent{\bf (ii). } Next, we consider the case when $p=0$. Notice that we have $N = N(d, 0) = 1$, and hence from (\ref{eqn:E11_final_control}) we have $\mathbb{E} \mathbf{E}_{11}  = 0$. 
From Lemma \ref{lemma:calcul_N}, when $d \geq {\mathfrak{C}}$,
   a sufficiently large constant defined in Lemma \ref{lemma:calcul_N}, we have $\lambda_1=\mu_0$ and $\lambda_2 = \mu_1$. Similar to (\ref{eqn:use_bound_L2}), we have
\begin{equation}\label{eqn:use_bound_Lp_0}
    \begin{aligned}
     &~ \mathbb{E}\mathbf{E}_1 
    \leq
    \mathbb{E}\mathbf{E}_{12}
    \overset{(\ref{eqn:E12})}{=}
    \frac{2N}{n}\sum_{j=2}^{\infty} \theta_{j}^2
    \leq \frac{2}{n} \cdot \mu_1^{s} R.
    \end{aligned}
\end{equation}

\subsubsection{Term $\mathbf{E}_2$}
We have
\begin{align}\label{eqn:main_decomposotion_E_2_in_upper_bound}
     \mathbb{E} \mathbf{E}_2 &\leq \sqrt{\mathcal{D}_{0}^{\star} \cdot 
    \mathbb{E}\mathbf{E}_1}.
\end{align}

\subsubsection{Term $\mathbf{E}_3$} We have 
\begin{align}\label{eqn:main_decomposotion_E_3_in_upper_bound}
    \mathbb{E}\mathbf{E}_3 &= \frac{\sigma^2}{n} \sum_{j=1}^{N}\ell_j^2 \mathbb{E}\left[ 
    \frac{1}{n}\sum_{i=1}^{n}\phi_j^2(x_i) - 1
    \right] = 0.
\end{align}

\subsubsection{Final result}
When $p > 0$, as shown in Corollary \ref{coroll:quantities} (and also in Appendices \ref{lem:quantities_1}), we have
$
\mathcal{D}^{\star} = \Omega_d(d^{p-\gamma})
 \gg \frac{2N}{n} \cdot \mu_1^{s} R$, hence for any $\varepsilon > 0$, there exist a constant $D_{\varepsilon}$ only depending on $\varepsilon$ and  ${\mathfrak{C}}$ defined in Lemma \ref{lemma:calcul_N}, such that for any $d > D_{\varepsilon}$, and
    for any regression function $f_{\star} \in \sqrt{R}[\calB]^{s}$ satisfying $\mathbb{E}_x f_{\star}(x) = 0$, we have
 \begin{equation*}
       \mathbb{E}\left[
    \|\hat{f}_{\ell, 0}-f_{\star}\|_{L^2}^2
    \right]
    \leq
    \mathcal{D}_{0}^{\star}
    +
    \mathbb{E}\mathbf{E}_1
    +
    \mathbb{E}\mathbf{E}_2
    +
    \mathbb{E}\mathbf{E}_3
    \leq
    \mathcal{D}^{\star}(1+\varepsilon).
    \end{equation*}
    
Similarly, when $p=0$, from Corollary \ref{coroll:quantities}, we have
$
\mathcal{D}^{\star} = \Omega_d(d^{-\gamma}) \gg \frac{2}{n} \cdot \mu_1^{s} R
$, hence there exist a constant $D_{\varepsilon}$ only depending on $\varepsilon$ and  ${\mathfrak{C}}$ defined in Lemma \ref{lemma:calcul_N}, such that for any $d > D_{\varepsilon}$, and
    for any regression function $f_{\star} \in \sqrt{R}[\calB]^{s}$, we have
\begin{equation*}
       \mathbb{E}\left[
    \|\hat{f}_{\ell, 0}-f_{\star}\|_{L^2}^2
    \right]
    \leq
    \mathcal{D}_{0}^{\star}
    +
    \mathbb{E}\mathbf{E}_1
    +
    \mathbb{E}\mathbf{E}_2
    +
    \mathbb{E}\mathbf{E}_3
    \leq
    \mathcal{D}^{\star}(1+\varepsilon).
    \end{equation*}
\end{proof}

The following Theorem proves (\ref{eqn:upper_bound}) when $N =\sum_{k=0}^{p+1}N(d, k)$.

\begin{theorem}[Restate Theorem \ref{thm:upper_1_copy} when $q=p+1$]\label{thm:upper_2}
    Suppose the same conditions as Theorem \ref{thm:main_pinsker_constant}. Further, suppose that $N =\sum_{k=0}^{p+1}N(d, k)$. 
    Then, for any $\varepsilon > 0$, there exist a constant $D_{\varepsilon}$ only depending on $\varepsilon$ and  ${\mathfrak{C}}$ defined in Lemma \ref{lemma:calcul_N}, such that for any $d > D_{\varepsilon}$, and
    for any regression function $f_{\star} \in \sqrt{R}[\calB]^{s}$ satisfying one of the following conditions: (i) $\mathbb{E}_x f_{\star}(x) = 0$ or (ii) $p=0$, we have
    \begin{equation*}
       \mathbb{E}\left[
    \|\hat{f}_{\ell, 0}-f_{\star}\|_{L^2}^2
    \right]
    \leq
    \mathcal{D}^{\star}(1+\varepsilon).
    \end{equation*}
\end{theorem}
\begin{proof}
    Recall that from the proof in Theorem \ref{thm:upper_1},  we have the following decomposition:
    \begin{align*}
        \mathbb{E}\left[
    \|\hat{f}_{\ell, 0}-f_{\star}\|_{L^2}^2
    \mid x_1, \cdots, x_n\right]
    &\leq \mathcal{D}_{0}^{\star} + \mathbf{E}_1 + \mathbf{E}_2 +
    \mathbf{E}_3,
    \end{align*}
and from (\ref{eqn:D0}), (\ref{eqn:main_decomposotion_E_2_in_upper_bound}), and (\ref{eqn:main_decomposotion_E_3_in_upper_bound}), 
\iffalse
we have
\begin{equation*}
    \begin{aligned}
        \mathcal{D}_{0}^{\star}\leq \mathcal{D}^{\star}, \quad
        \mathbb{E} \mathbf{E}_2 \leq \sqrt{\mathcal{D}^{\star} \cdot \mathbb{E}\mathbf{E}_1}, \quad
        \mathbb{E}\mathbf{E}_3 =0,
    \end{aligned}
\end{equation*}
therefore, 
\fi
we only need to show that 
$$
 \mathbb{E}\mathbf{E}_1 \leq \mathcal{D}^{\star}\varepsilon.
$$        

Denote $N^{\prime} = \sum_{k=0}^{p}N(d, k)$. We have
\begin{equation}\label{eqn:E1m}
    \begin{aligned}
      &~  \mathbf{E}_1 
    =
    \sum_{j=1}^\infty\ell_j^2 \left(\sum_{j^{\prime}=1}^{\infty} \theta_{j^{\prime}} \Delta_n(j, j^{\prime})\right)^2\\
    =&~
    \sum_{j=1}^{N^{\prime}}\ell_j^2 \left(\sum_{j^{\prime}=1}^{\infty} \theta_{j^{\prime}} \Delta_n(j, j^{\prime})\right)^2 + \sum_{j=N^{\prime}+1}^{N}\ell_j^2 \left(\sum_{j^{\prime}=1}^{\infty} \theta_{j^{\prime}} \Delta_n(j, j^{\prime})\right)^2\\
    \leq &~
    \underbrace{2\sum_{j=1}^{N^{\prime}} \left(\sum_{j^{\prime}=1}^{N} \theta_{j^{\prime}} \Delta_n(j, j^{\prime})\right)^2}_{2\mathbf{E}_{13}}
    +
    \underbrace{2\sum_{j=1}^{N^{\prime}}\left(\sum_{j^{\prime}=N+1}^{\infty} \theta_{j^{\prime}} \Delta_n(j, j^{\prime})\right)^2}_{2\mathbf{E}_{14}}\\
   &~ + 
    \underbrace{2\sum_{j=N^{\prime}+1}^{N} \ell_j^2 \left(\sum_{j^{\prime}=1}^{N} \theta_{j^{\prime}} \Delta_n(j, j^{\prime})\right)^2}_{2\mathbf{E}_{15}}
    +
    \underbrace{2\sum_{j=N^{\prime}+1}^{N} \ell_j^2
    \left(\sum_{j^{\prime}=N+1}^{\infty} \theta_{j^{\prime}} \Delta_n(j, j^{\prime})\right)^2}_{2\mathbf{E}_{16}}.
    \end{aligned}
\end{equation}

For the first term in (\ref{eqn:E1m}), we have
\begin{equation}\label{eqn:E13}
\begin{aligned}
&~\mathbb{E}\mathbf{E}_{13}
=
   \mathbb{E}\sum_{j=1}^{N^{\prime}} \left(\sum_{j^{\prime}=1}^{N} \theta_{j^{\prime}} \Delta_n(j, j^{\prime})\right)^2\\
=&~  \sum_{j=1}^{N^{\prime}} \sum_{j'=1}^{N} \theta_{j'}^2 \mathbb{E} \Delta_n(j,j')^2 + \sum_{j=1}^{N^{\prime}} \sum_{u \neq v}^{N} \theta_u \theta_v \mathbb{E} \left[\Delta_n(j,u) \Delta_n(j,v)\right]\\
\overset{(\ref{eqn:E1_1}), (\ref{eqn:E1_2}), \text{ and } (\ref{eqn:E1_3})}{=}
&~
\sum_{j'=1}^{N} \theta_{j'}^2 \left[\frac{N^{\prime}-\mathbf{1}(j' \leq N^{\prime})}{n}\right]
=
\frac{N^{\prime}}{n}\sum_{j'=1}^{N} \theta_{j'}^2
-
\frac{1}{n}\sum_{j'=1}^{N^{\prime}} \theta_{j'}^2,
    \end{aligned}
\end{equation}
where in the third equation we use the Addition Formula.

For the second term in (\ref{eqn:E1m}), similarly we have
\begin{equation}\label{eqn:E14}
\begin{aligned}
&~\mathbb{E}\mathbf{E}_{14}
=
   \mathbb{E}\sum_{j=1}^{N^{\prime}} \left(\sum_{j^{\prime}=N+1}^{\infty} \theta_{j^{\prime}} \Delta_n(j, j^{\prime})\right)^2\\
   =&~  \mathbb{E}\sum_{j=1}^{N^{\prime}} \sum_{j^{\prime}=N+1}^{\infty} \theta_{j^{\prime}}^2  \Delta_n(j,j')^2 + \mathbb{E}\sum_{j=1}^{N^{\prime}} \sum_{u \neq v \geq N+1}^{\infty} \theta_u \theta_v  \left[\Delta_n(j,u) \Delta_n(j,v)\right]\\
=&~  \sum_{j=1}^{N^{\prime}} \sum_{j^{\prime}=N+1}^{\infty} \theta_{j^{\prime}}^2 \mathbb{E} \Delta_n(j,j')^2 + \sum_{u \neq v \geq N+1}^{\infty} \theta_u \theta_v \sum_{j=1}^{N^{\prime}} \mathbb{E} \left[\Delta_n(j,u) \Delta_n(j,v)\right]\\
\overset{(\ref{eqn:E1_1})  \text{ and } (\ref{eqn:E1_3})}{=}
&~
 \sum_{j^{\prime}=N+1}^{\infty} \theta_{j^{\prime}}^2 \frac{N^{\prime}}{n},
\end{aligned}\end{equation}
where the interchangeable order of infinite summation and expectation in the third equation can be argued similar to $\mathbf{E}_{12}$ in (\ref{eqn:E12}).

For the third term in (\ref{eqn:E1m}), 
notice that from we have $\ell_{N^{\prime}+1} = \cdots = \ell_{N} = 1 - \kappa^{\star}\mu_{p+1}^{-s/2}$, and hence
\begin{equation}\label{eqn:E15}
\begin{aligned}
&~\mathbb{E}\mathbf{E}_{15}
=
   \mathbb{E}\sum_{j=N^{\prime}+1}^{N} \ell_j^2 \left(\sum_{j^{\prime}=1}^{N} \theta_{j^{\prime}} \Delta_n(j, j^{\prime})\right)^2\\
= &~ \ell_{N}^2 \cdot\left[\sum_{j=N^{\prime}+1}^{N} \sum_{j'=1}^{N} \theta_{j'}^2 \mathbb{E} \Delta_n(j,j')^2 + \sum_{j=N^{\prime}+1}^{N} \sum_{u \neq v}^{N} \theta_u \theta_v \mathbb{E} \left[\Delta_n(j,u) \Delta_n(j,v)\right]\right]\\
\overset{(\ref{eqn:E1_1}) - (\ref{eqn:E1_4})}{=} &~
\ell_{N}^2 \cdot \sum_{j'=1}^{N} \theta_{j'}^2 \left[\frac{N - N^{\prime}-\mathbf{1}(N^{\prime} < j' \leq N)}{n}\right]\\
= &~
\ell_{N}^2 \cdot \left[\frac{N-N^{\prime}}{n}\sum_{j'=1}^{N} \theta_{j'}^2
-
\frac{1}{n}\sum_{j'=N^{\prime}+1}^{N} \theta_{j'}^2\right].
    \end{aligned}
\end{equation}

For the fourth term in (\ref{eqn:E1m}), we have
\begin{equation}\label{eqn:E16}
\begin{aligned}
&~\mathbb{E}\mathbf{E}_{16}
=
   \mathbb{E}\sum_{j=N^{\prime}+1}^{N} \ell_j^2
    \left(\sum_{j^{\prime}=N+1}^{\infty} \theta_{j^{\prime}} \Delta_n(j, j^{\prime})\right)^2\\
=&~  \ell_{N}^2 \cdot\left[\sum_{j=N^{\prime}+1}^{N} \sum_{j^{\prime}=N+1}^{\infty} \theta_{j^{\prime}}^2 \mathbb{E} \Delta_n(j,j')^2 + \sum_{j=N^{\prime}+1}^{N} \sum_{u \neq v \geq N+1}^{\infty} \theta_u \theta_v \mathbb{E} \left[\Delta_n(j,u) \Delta_n(j,v)\right]\right]\\
\overset{(\ref{eqn:E1_1})  \text{ and } (\ref{eqn:E1_3})}{=}
&~
 \ell_{N}^2 \cdot \sum_{j^{\prime}=N+1}^{\infty} \theta_{j^{\prime}}^2 \frac{N-N^{\prime}}{n},
\end{aligned}\end{equation}
where the interchangeable order of infinite summation and expectation in the second equation can be argued similar to $\mathbf{E}_{12}$ in (\ref{eqn:E12}).

Now we begin to bound $\mathbf{E}_1$. We separate the proof into the following two cases.

\noindent{\bf (i). } We first consider the case when $p>0$ and $\theta_1=0$.
Notice that from (\ref{eqn:quantities_2_bound_ii}) we have $\ell_{N}^2 =  O_d(d^{2\gamma - 2(s+1)(p+1)})$.
Furthermore, from Corollary \ref{coroll:quantities} we have $\frac{2N^{\prime}}{n} = O_d(\calD^{\star})$ and $\ell_{N}^2 \cdot \frac{N}{n} = o_d(\calD^{\star})$.
Finally, similar to (\ref{eqn:use_bound_L2}), since $\theta_1=0$, we have $\sum_{j=1}^{\infty} \theta_{j}^2 \leq \mu_1^{s} R = o_d(1)$.

Therefore, for any $\varepsilon > 0$, there exist a constant $D_{\varepsilon}$ only depending on $\varepsilon$ and  ${\mathfrak{C}}$ defined in Lemma \ref{lemma:calcul_N}, such that for any $d > D_{\varepsilon}$, and
    for any regression function $f_{\star} \in \sqrt{R}[\calB]^{s}$ satisfying $\mathbb{E}_x f_{\star}(x) = 0$, we have
\begin{equation}
    \begin{aligned}
        \mathbb{E}\mathbf{E}_{1} \leq
        2\mathbb{E}\mathbf{E}_{13}+2\mathbb{E}\mathbf{E}_{14} + 2\mathbb{E}\mathbf{E}_{15}+2\mathbb{E}\mathbf{E}_{16}
        \leq
        \frac{2N^{\prime}}{n} \sum_{j=1}^{\infty} \theta_{j}^2 + \ell_{N}^2 \cdot \frac{N}{n}\sum_{j=1}^{\infty} \theta_{j}^2
        \leq \mathcal{D}^{\star}\varepsilon.
    \end{aligned}
\end{equation}

\noindent{\bf (ii). } Next, we consider the case when $p=0$. Notice that we have $N^{\prime} = N(d, 0) = 1$, and hence from (\ref{eqn:E13}) we have $\mathbb{E} \mathbf{E}_{13}  = \frac{1}{n} \sum\nolimits_{j=2}^{N}\theta_j^2$. 
Similar to above, we have $\frac{2}{n} = O_d(\calD^{\star})$, $\ell_{N}^2 \cdot \frac{N}{n} = o_d(\calD^{\star})$, and $\sum_{j=2}^{\infty} \theta_{j}^2= o_d(1)$. 

Therefore, for any $\varepsilon > 0$, there exist a constant $D_{\varepsilon}$ only depending on $\varepsilon$ and  ${\mathfrak{C}}$ defined in Lemma \ref{lemma:calcul_N}, such that for any $d > D_{\varepsilon}$, and
    for any regression function $f_{\star} \in \sqrt{R}[\calB]^{s}$, we have
\begin{equation}\label{eqn:use_bound_Lp_0_q_great_than_p}
    \begin{aligned}
        \mathbb{E}\mathbf{E}_{1} \leq
        2\mathbb{E}\mathbf{E}_{13}+2\mathbb{E}\mathbf{E}_{14} + 2\mathbb{E}\mathbf{E}_{15}+2\mathbb{E}\mathbf{E}_{16}
        \leq 
        \frac{2}{n} \sum_{j=2}^{\infty} \theta_{j}^2 + \ell_{N} \cdot \frac{N}{n}\sum_{j=1}^{\infty} \theta_{j}^2
        \leq \mathcal{D}^{\star}\varepsilon.
    \end{aligned}
\end{equation}
\end{proof}

\subsection{Proof of (\ref{eqn:upper_bound})}

Now we can give the final result.
Recall that in Section \ref{sec:upper_sketch}, we define the linear filter estimator as:
\begin{align*}
    \hat{f}_{\ell}(x) :=
    (\ell_1\mathbf{1}\{p=0\}+\mathbf{1}\{p > 0\}) \bar{z}_1 +\hat{g}_{\ell}(x) \quad
    \mbox{ where } \quad \hat{g}_{\ell}(x) = \sum_{j=2}^{N} \ell_j \bar{z}_j \phi_j(x).
\end{align*}
where $\ell_j$'s are given in Definition \ref{def:pinsker_quantity} and $\bar{z}_j$'s are given in (\ref{def:bar_z_j}).

\begin{theorem}\label{thm:upper_final}
    Suppose the same conditions as Theorem \ref{thm:main_pinsker_constant}.
    Then, when $d \geq {\mathfrak{C}}$,
   a sufficiently large constant defined in Lemma \ref{lemma:calcul_N}, 
    we have
    \begin{equation*}
        \inf_{\hat{f}} \sup_{
        \rho_{f_{\star}} \in \mathcal{P}
        %f_{\star} \in \sqrt{R}[\calB]^{s}
        } \mathbb{E}\left[
    \|\hat{f}-f_{\star}\|_{L^2}^2
    \right]
     \leq
        \sup_{
        \rho_{f_{\star}} \in \mathcal{P}
        %f_{\star} \in \sqrt{R}[\calB]^{s}
        } 
    \mathbb{E}\left[
    \|\hat{f}_{\ell}-f_{\star}\|_{L^2}^2
    \right]
    \leq
    \mathcal{D}^{\star}(1+o_d(1)).
    \end{equation*}
\end{theorem}
\begin{proof}
From Lemma \ref{lemma:calcul_N}, when $d \geq {\mathfrak{C}}$,
   a sufficiently large constant defined in Lemma \ref{lemma:calcul_N}, we have $\lambda_1=\mu_0$ and $\phi_1(x) = Y_{0, 1}(x)\equiv 1$. 

Notice that when $p=0$, Theorem \ref{thm:upper_1} and \ref{thm:upper_2}
imply that
$$
\sup_{
        \rho_{f_{\star}} \in \mathcal{P}
        %f_{\star} \in \sqrt{R}[\calB]^{s}
        } 
    \mathbb{E}\left[
    \|\hat{f}_{\ell}-f_{\star}\|_{L^2}^2
    \right]
    \leq
    \mathcal{D}^{\star}(1+o_d(1)),
$$
and hence we only need to prove the case when $p \geq 1$.

For any $f_{\star}(\cdot) = \sum_{j=1}^{\infty} \theta_j \phi_j(\cdot) \in \sqrt{R}[\calB]^{s}$, denote $g_{\star}(x) = \sum_{j = 2}^{\infty} \theta_j \phi_j(x)$ where $\phi_j$'s are the eigenfunctions defined in (\ref{eqn:mercer_decomp}).
\iffalse
Denote
\begin{align*}
    g_{\star}(x) = \sum_{j \geq 2} \theta_j \phi_j(x)\quad \text{ and } \quad
     \hat{g}_{\ell}(x) = \sum_{j=2}^{N} \ell_j \bar{z}_j \phi_j(x),
\end{align*}
then we have
\begin{align*}
    f_{\star}(x) &= \theta_1 + g_{\star}(x)\\
    \hat{f}_{\ell}(x) &= 
    \bar{z}_1 +
    \hat{g}_{\ell}(x) = \frac{1}{n} \sum_{i=1}^{n} y_i + \hat{g}_{\ell}(x).
\end{align*}
\fi
Recall that when $d \geq \mathfrak{C}$, from Lemma \ref{lemma:calcul_N} we have
\begin{align*}
     \mu_0 &> \mu_1 > \cdots > \mu_{p+1} > \max_{j \geq p+2} \mu_{j},
   \end{align*}
hence we have $\lambda_1 = \mu_0$, $\lambda_2 = \mu_1$, and $\phi_1=1$.
Moreover, since for any $j \geq 2$, $\phi_j$ is orthogonal to $\phi_1 \equiv 1$, we have $\mathbb{E}_{x}(g_{\star}(x))=\mathbb{E}_{x}(\hat{g}_{\ell}(x))=0$. Therefore, 
    \begin{equation}
    \begin{aligned}
       &~ \mathbb{E}\left[
    \|\hat{f}_{\ell}-f_{\star}\|_{L^2}^2
    \right]
    =
    \mathbb{E} \left[\int \left(\hat{g}_{\ell}(x)-g_{\star}(x) + \left(\frac{1}{n} \sum_{i=1}^{n} y_i - \theta_1\right)\right)^2 \rho_{\calX}(x) ~\mathsf{d}x \right]\\
    =&~
    \mathbb{E} \left[\int \left(\hat{g}_{\ell}(x)-g_{\star}(x)\right)^2 \rho_{\calX}(x) ~\mathsf{d}x \right]\\
    &~+
    2 \mathbb{E}\left[\left(\frac{1}{n} \sum_{i=1}^{n} y_i - \theta_1\right) \int \hat{g}_{\ell}(x)-g_{\star}(x) \rho_{\calX}(x) ~\mathsf{d}x\right]
    +
    \mathbb{E}\left(\frac{1}{n} \sum_{i=1}^{n} y_i - \theta_1\right)^2\\
    {=}&~
    \mathbb{E}\left[
    \|\hat{g}_{\ell}-g_{\star}\|_{L^2}^2
    \right]
    +
    \mathbb{E}\left(\frac{1}{n} \sum_{i=1}^{n} y_i - \theta_1\right)^2.
    %\leq &~ \mathcal{D}^{\star}(1+o_d(1)),
    \end{aligned}
    \end{equation}

Denote $\mathbf{I} = \left(\frac{1}{n} \sum_{i=1}^{n} y_i - \theta_1\right)^2 - \sigma^2/n$. Since 
$\mathbb{E}(y_i \mid x_i)=\theta_1 + g_{\star}(x_i)$ and $\text{Var}(y_i \mid x_i)\leq \sigma^2$, we have
\begin{equation*}
    \begin{aligned}
        &~\mathbb{E}\mathbf{I} =
        \mathbb{E}\left(\mathbb{E}\left[~\mathbf{I} \mid \{x_1, \cdots, x_n\} \right]\right)
        \leq 
        \frac{1}{n^2} \mathbb{E}\sum_{i=1}^n g_{\star}^2(x_i)
        =
        \frac{1}{n}\sum_{j=2}^{\infty.} \theta_j^2
        \leq
    \frac{\mu_1^{s}}{n} \sum_{j=2}^{\infty} \lambda_j^{-s}\theta_{j}^2
        \leq  \frac{\mu_1^{s}}{n} R.
    \end{aligned}
\end{equation*}
Therefore, from Corollary \ref{coroll:quantities}, for any $\varepsilon > 0$, there exist a constant $D_{\varepsilon, 1}$ only depending on $\varepsilon$ and  ${\mathfrak{C}}$ defined in Lemma \ref{lemma:calcul_N}, such that for any $d > D_{\varepsilon, 1}$, and
    for any regression function $f_{\star} \in \sqrt{R}[\calB]^{s}$, we have
$$
\mathbb{E}\left(\frac{1}{n} \sum_{i=1}^{n} y_i - \theta_1\right)^2 
\leq 
\frac{\sigma^2}{n} + \frac{\mu_1^{s}}{n} R \leq \mathcal{D}^{\star}\varepsilon.
$$

On the other side, since $\mathbb{E}_{x}(g_{\star}(x))=0$, from Theorem \ref{thm:upper_1} and Theorem \ref{thm:upper_2}, there exist a constant $D_{\varepsilon}$ only depending on $\varepsilon$ and  ${\mathfrak{C}}$ defined in Lemma \ref{lemma:calcul_N}, such that for any $d > D_{\varepsilon}$, and
    for any regression function $f_{\star} \in \sqrt{R}[\calB]^{s}$, we have
$$
\mathbb{E}\left[
    \|\hat{g}_{\ell}-g_{\star}\|_{L^2}^2
    \right] \leq \mathcal{D}^{\star}(1+\varepsilon), 
$$
hence when $d \geq \mathfrak{C}$, by the definition of $o_d(1)$, we have
$$
\sup_{
        \rho_{f_{\star}} \in \mathcal{P}
        %f_{\star} \in \sqrt{R}[\calB]^{s}
        } 
    \mathbb{E}\left[
    \|\hat{f}_{\ell}-f_{\star}\|_{L^2}^2
    \right] \leq \mathcal{D}^{\star}(1+o_d(1)),
$$
finishing our proof.
\end{proof}

\section{Proof of lower bound in Theorem \ref{thm:main_pinsker_constant}}\label{sec:lower}

In this section, our goal is to show that
\begin{equation*}
        \inf_{\hat{f}} \sup_{
        \rho_{f_{\star}} \in \mathcal{P}
        %f_{\star} \in \sqrt{R}[\calB]^{s}
        } \mathbb{E}_{(X, Y) \overset{\mathcal{D}}{\sim} \rho_{f_{\star}}^{\otimes n}}\left[
    \|\hat{f}-f_{\star}\|_{L^2}^2
    \right]
     \geq
    \mathcal{D}^{\star}(1+o_d(1)).
    \end{equation*}

Denote
\begin{equation}\label{def:theta_spca_lower_bound}
    \begin{aligned}
\Theta_N &:= \left\{ \theta^N = (\theta_1, \ldots, \theta_N)^{\top} \in \mathbb{R}^{N} : \sum_{j=1}^N \lambda_j^{-s} \theta_j^2 \leq R \right\},\\
\calF_N &:= \left\{ \sum_{j=1}^N \theta_j \phi_j(\cdot) : \sum_{j=1}^N \lambda_j^{-s} \theta_j^2 \leq R \right\} \subset \sqrt{R}[\calB]^{s}.
\end{aligned}
\end{equation}
Recall that $\rho_{\calX}$ is the uniform distribution on $\mathbb{S}^{d}$. Let's denote 
\begin{align*}
    \tilde{\mathcal{P}}=\Bigg\{ \tilde{\rho}_{f_{\star}}~\bigg|~ \mbox{ joint distribution of $(x, y)$ where } x \overset{\calD}{\sim} \rho_{\calX}, y=f_{\star}(x)+\epsilon,
    \epsilon \overset{\calD}{\sim} \calN(0,\sigma^{2}),
    f_{\star}\in \calF_N \Bigg\}.
\end{align*}
It is easy to see that we have $\tilde{\mathcal{P}} \subset \calP$, the set of all the distributions $\rho_{f_{\star}}$ on $\mathcal{X} \times \mathcal{Y}$ given by (\ref{equation:true_model}) such that  Assumption 
\ref{assump_asymptotic}, 
\ref{assu:coef_of_inner_prod_kernel},  and \ref{assump_function_calss} hold for some $\alpha, \gamma>0$.
Therefore, if we denote $\tilde{\rho}_{f_{\star}}$ as the distribution in $\tilde{\mathcal{P}}$ with respect to $f_{\star}$, and denote $\mathbb{E} = \mathbb{E}_{(X, Y) \overset{\mathcal{D}}{\sim} \tilde{\rho}_{f_{\star}}^{\otimes n}}$ for notation simplicity, then we have
\begin{equation}\label{eqn:parametrized_lower_bound}
    \begin{aligned}
&~ \inf_{\hat{f}} \sup_{
        \rho_{f_{\star}} \in \mathcal{P}} \mathbb{E}_{(X, Y) \overset{\mathcal{D}}{\sim} \rho_{f_{\star}}^{\otimes n}}\left[
    \|\hat{f}-f_{\star}\|_{L^2}^2
    \right]
     \geq
     \inf_{\hat{f}} \sup_{\tilde{\rho}_{f_{\star}} \in \tilde{\mathcal{P}}} 
    \mathbb{E}
    \left[
    \|\hat{f}-f_{\star}\|_{L^2}^2
    \right]\\
    =&~
    \inf_{\hat{f}} \sup_{f_{\star}\in \calF_N} 
    \mathbb{E}
    \left[
    \|\hat{f}-f_{\star}\|_{L^2}^2
    \right]
\geq \inf_{\hat{f} \in \calF_N} \sup_{f_{\star} \in \calF_N} \mathbb{E} \|\hat{f} - f_{\star}\|_{L^2}^2 \quad \text{a.s.},\\
=&~\inf_{\hat{\theta}^N \in \Theta_N} \sup_{\theta^N \in \Theta_N} \mathbb{E} \left[ \sum_{j=1}^N (\hat{\theta}_j - \theta_j)^2 \right]
:= \mathbf{I},
\end{aligned}
\end{equation}
where the second inequality is because for all $f_{\star}\in\mathcal{F}_N$ and all estimator $\hat{f}$, there exists a random function $\hat{f}_{\mathcal{F}_N} \in\mathcal{F}_N$ such that $\|\hat{f}-f\|_2^2\geq\|\hat{f}_{\mathcal{F}_N}-f\|_2^2$ almost surely. 
For readers' convenience, we borrow the corresponding explanation from \cite{tsybakov2008introduction} as follows:
In fact, if the realization $\{(x_i,y_i)\}_{i \leq n}$ is such that $\hat{f} \in L^2$, it is sufficient to take as estimator $\sum\nolimits_{j=1}^{N}\hat{\theta}_j\phi_j$ the $L^2$ projection of $\hat{f}$ on $\mathcal{F}_N$ (indeed, $\mathcal{F}_N$ is a closed convex set in $L^2$). If $\hat{f}\notin L^2$, then $\|\hat{f}-f_{\star}\|_{L^2}^2=+\infty$ and $\|\hat{f}-f_{\star}\|_{L^2}^2\geq\|\hat{f}_{\mathcal{F}_N}-f_{\star}\|_{L^2}^2$ is trivial for all $\hat{f}_{\mathcal{F}_N} \in \mathcal{F}_N$.

\subsection{Parametric case}\label{subsec_para_case}

When $\gamma < s$, we obtain the following lower bound.

\begin{theorem}\label{thm_lower_para}
    Suppose the same conditions as Theorem \ref{thm:main_pinsker_constant}.
    When $d \geq {\mathfrak{C}}$,
   a sufficiently large constant defined in Lemma \ref{lemma:calcul_N}, 
    if $\gamma <s$, then we have
    \begin{equation*}
        \inf_{\hat{f}} \sup_{
        \rho_{f_{\star}} \in \mathcal{P}
        %f_{\star} \in \sqrt{R}[\calB]^{s}
        } \mathbb{E}_{(X, Y) \overset{\mathcal{D}}{\sim} \rho_{f_{\star}}^{\otimes n}}\left[
    \|\hat{f}-f_{\star}\|_{L^2}^2
    \right]
     \geq
    \mathcal{D}^{\star}(1+o_d(1)).
    \end{equation*}
\end{theorem}
\begin{proof}
From the definition we have $p=\lfloor\frac{\gamma}{s+1} \rfloor=0$.
From Lemma \ref{lemma:calcul_N} we know that either $q=0$ and $N=1$, or $q=1$ and $N=d+1$. 
    We first consider the case $p=0$, $q=0$ and $N=1$. From (\ref{eqn:parametrized_lower_bound}) we have
    $$
    \inf_{\hat{f}} \sup_{
        \rho_{f_{\star}} \in \mathcal{P}
        %f_{\star} \in \sqrt{R}[\calB]^{s}
        } \mathbb{E}_{(X, Y) \overset{\mathcal{D}}{\sim} \rho_{f_{\star}}^{\otimes n}}\left[
    \|\hat{f}-f_{\star}\|_{L^2}^2
    \right]
    \geq
    \mathbf{I} = \inf_{\hat{\theta}_1\in \Theta_1} \sup_{\theta_1 \in \Theta_1} \mathbb{E} \left[ (\hat{\theta}_1 - \theta_1)^2 \right],
    $$
    where $\Theta_1 = \{\theta_1: \theta_1^2 \leq R\mu_0^s \}$.

For any $\theta_1 \in \Theta_1$, note that we have $y_i  \overset{\mathcal{D}}{\sim}_{i.i.d.} \calN(\theta_1, \sigma^2)$, and it is a well-known result that we have
$$
\inf_{\hat{\theta}_1(\{y_i\}_{i =1}^{n})} \sup_{\theta_1 \in \Theta_1} \mathbb{E} \left[ (\hat{\theta}_1(\{y_i\}_{i =1}^{n}) - \theta_1)^2 \right] = 
\sup_{\theta_1 \in \Theta_1} \mathbb{E} \left[ \left(\frac{1}{n}\sum_{i=1}^{n}y_i - \theta_1\right)^2 \right]
=\frac{\sigma^2}{n},
$$
see, e.g., page 121 in \cite{shao2003mathematical}. 
\iffalse
We also provide an alternative proof as follows.

For any unbiased estimator $\tilde{\theta}_1(\{y_i\}_{i \leq n})$ of $\theta_1$, from Cramér-Rao lower bound, we have
\begin{align*}
    \text{MSE}(\tilde{\theta}_1(\{y_i\}_{i \leq n})):=\mathbb{E} \left[ (\tilde{\theta}_1(\{y_i\}_{i \leq n}) - \theta_1)^2 \right] \geq 
    \mathbb{E} \left[ \left(\frac{1}{n}\sum_{i\leq n} y_i - \theta_1 \right)^2 \right] = \frac{\sigma^2}{n}.
\end{align*}
Notice that for any estimator $\hat{\theta}_1$ of $\theta_1$, if $\hat{\theta}_1$ is biased, then we have
$$
\text{MSE}(\hat{\theta}_1) \geq \text{MSE}(\hat{\theta}_1 - \mathbb{E}\hat{\theta}_1 + \theta_1) \geq \frac{\sigma^2}{n}.
$$
Therefore, for any estimator $\hat{\theta}_1$ of $\theta_1$, we have
$$
 \sup_{\theta_1 \in \Theta_1} \mathbb{E} \left[ (\hat{\theta}_1 - \theta_1)^2 \right] \geq \frac{\sigma^2}{n},
$$
\fi
Therefore, we have
\begin{align*}
    \inf_{\hat{f}} \sup_{
        \rho_{f_{\star}} \in \mathcal{P}
        %f_{\star} \in \sqrt{R}[\calB]^{s}
        } \mathbb{E}_{(X, Y) \overset{\mathcal{D}}{\sim} \rho_{f_{\star}}^{\otimes n}}\left[
    \|\hat{f}-f_{\star}\|_{L^2}^2
    \right] \geq \inf_{\hat{\theta}_1\in \Theta_1} \sup_{\theta_1 \in \Theta_1} \mathbb{E} \left[ (\hat{\theta}_1 - \theta_1)^2 \right] = \frac{\sigma^2}{n} \overset{\text{Corollary } \ref{coroll:quantities}}{\sim} \calD^{\star}.
\end{align*}

Now we consider the case $p=0$, $q=1$ and $N>1$.
From Lemma \ref{lemma:calcul_N}, when $d \geq {\mathfrak{C}}$,
   a sufficiently large constant defined in Lemma \ref{lemma:calcul_N}, we have $\lambda_1=\mu_0$.
From 
(\ref{eqn:quantities_2_bound_ii_coro}) and (\ref{eqn:quantities_2_bound_ii}) we have $\ell_1 \sim 1$ and $\sum_{j=2}^{N}\ell_j \leq N \ell_{2} =  O_d(d \cdot d^{\gamma-s-1})=o_d(1)$. Hence we have 
$$
\calD^{\star} \sim \frac{\sigma^2}{n}.
$$
Therefore, we have
\begin{equation*}
    \begin{aligned}
\inf_{\hat{f}} \sup_{
        \rho_{f_{\star}} \in \mathcal{P}
        %f_{\star} \in \sqrt{R}[\calB]^{s}
        } \mathbb{E}_{(X, Y) \overset{\mathcal{D}}{\sim} \rho_{f_{\star}}^{\otimes n}}\left[
    \|\hat{f}-f_{\star}\|_{L^2}^2
    \right]
    \geq
    \mathbf{I}
\geq 
\inf_{\hat{\theta}_1 \in \Theta_1} \sup_{\theta_1 \in \Theta_1} \mathbb{E} \left[ (\hat{\theta}_1 - \theta_1)^2 \right] \geq \frac{\sigma^2}{n} \sim \calD^{\star},
\end{aligned}
\end{equation*}
finishing the proof.
\end{proof}

\iffalse
\begin{remark}
    When $s/2 < \gamma < s$, similar to the above proof, we have \(\sum_{j \geq 2} \ell_j = o_d(1)\) and \(\ell_1 \sim 1\), and thus we have $\calD^{\star} \sim \frac{\sigma^2}{n}$.
\end{remark}
\fi

\subsection{Non-parametric case}

When $\gamma \geq s$, we have the following lower bound.

\begin{theorem}\label{thm_lower_non_para}
    Suppose the same conditions as Theorem \ref{thm:main_pinsker_constant}.
    When $d \geq {\mathfrak{C}}$,
   a sufficiently large constant defined in Lemma \ref{lemma:calcul_N}, 
    if $\gamma \geq s$, then we have
    \begin{equation*}
        \inf_{\hat{f}} \sup_{
        \rho_{f_{\star}} \in \mathcal{P}
        %f_{\star} \in \sqrt{R}[\calB]^{s}
        } \mathbb{E}_{(X, Y) \overset{\mathcal{D}}{\sim} \rho_{f_{\star}}^{\otimes n}}\left[
    \|\hat{f}-f_{\star}\|_{L^2}^2
    \right]
     \geq
    \mathcal{D}^{\star}(1+o_d(1)).
    \end{equation*}
\end{theorem}

\begin{proof}
Fix any \(\delta \in (0, 1)\). Let
\begin{equation}\label{def:v_j_and_s_j}
    v_j^2 = \frac{\sigma^2\ell_j}{n\kappa^{\star}\lambda_j^{-s/2}} \quad \text{ and } \quad s_j^{2} = (1-\delta)v_j^2, \quad j =1,2,\ldots N.
\end{equation}
Denote \(\varphi(\cdot)\) as the p.d.f. of \(\mathcal{N}(0, 1)\), and \(\mu_s(t)=s^{-1}\varphi(t/s)\) as the p.d.f. of \(\mathcal{N}(0, s^2)\).
Suppose 
\begin{align*}
\theta^N \overset{\mathcal{D}}{\sim} \calN \left(\mathbf 0, \text{diag} \left(s_1^2, \ldots, s_N^2\right) \right),
\end{align*}
then we have $\mu(\theta^N) = \prod_{j=1}^N \mu_{s_j}(\theta_j)$.
Hence, from (\ref{eqn:parametrized_lower_bound}) we have
\begin{equation}\label{eqn:decomp_lower_bound}
    \begin{aligned}
    &~\inf_{\hat{f}} \sup_{
        \rho_{f_{\star}} \in \mathcal{P}
        %f_{\star} \in \sqrt{R}[\calB]^{s}
        } \mathbb{E}_{(X, Y) \overset{\mathcal{D}}{\sim} \rho_{f_{\star}}^{\otimes n}}\left[
    \|\hat{f}-f_{\star}\|_{L^2}^2
    \right]
    \geq 
\mathbf{I} =
\inf_{\hat{\theta}^N \in \Theta_N} \sup_{\theta^N \in \Theta_N} \mathbb{E} \left[ \sum_{j=1}^N (\hat{\theta}_j - \theta_j)^2 \right]\\
\geq &~ 
\underbrace{\inf_{\hat{\theta}^N \in \Theta_N} \sum_{k=1}^N \mathbb{E} \left[ \int_{\mathbb{R}^{N}} (\hat{\theta}_k - \theta_k)^2 \mu(\theta^N) ~\mathsf{d}\theta^N \right]}_{\mathbf{I}^{\star}}
-
\underbrace{\sup_{\hat{\theta}^N \in \Theta_N} \sum_{k=1}^N \mathbb{E} \left[ \int_{\mathbb{R}^{N} \backslash \Theta_N} (\hat{\theta}_k - \theta_k)^2 \mu(\theta^N) ~\mathsf{d}\theta^N \right]}_{\mathbf{r}^{\star}}.
\end{aligned}
\end{equation}

\subsection{Lower bound of $\mathbf{I}^{\star}$}
Notice that we have
\begin{equation}\label{eqn:bound_of_I_star}
\begin{aligned}  
 \mathbf{I}^{\star} \geq &~ \sum_{k=1}^N \inf_{\hat{\theta}_k} \int_{\mathbb{R}^{N}} \mathbb{E} \left[ (\hat{\theta}_k - \theta_k)^2 \right] \mu(\theta^N) ~\mathsf{d}\theta^N\\
 =&~\sum_{k=1}^{N} \inf_{\hat{\theta}_{k}}\bbE_{X}\int_{\bbR^{N}}\bbE_{\epsilon} \left[ (\hat{\theta}_k - \theta_k)^2 \mid (x_{1},\cdots,x_{n})\right] \mu(\theta^N) ~\mathsf{d}\theta^N\\
 \overset{\text{Fatou's lemma }}{\geq} &~\sum_{k=1}^{N} \bbE_{X}\inf_{\hat{\theta}_{k}}\int_{\bbR^{N}}\bbE_{\epsilon} \left[ (\hat{\theta}_k - \theta_k)^2 \mid (x_{1},\cdots,x_{n})\right] \mu(\theta^N) ~\mathsf{d}\theta^N\\
 =&~\sum_{k=1}^{N} \bbE_{X}\inf_{\hat{\theta}_{k}}\bbE_{\epsilon,\theta^{N}} \left[ (\hat{\theta}_k - \theta_k)^2 \mid (x_{1},\cdots,x_{n})\right]\\
\overset{\mathbf{(A)}}{\geq}  &~ \sum_{k=1}^N \mathbb{E}_X \frac{s_k^2\sigma^2}{\sigma^2+\sum_{i=1}^n \phi_k^2(x_i) s_k^2}   \geq   \sum_{k=1}^N \frac{s_k^2\sigma^2}{\mathbb{E}_X \left( \sigma^2+\sum_{i=1}^n \phi_k^2(x_i) s_k^2 \right)} \\
=&~ \sum_{k=1}^{N}\frac{(1-\delta)v_{k}^{2}\sigma^{2}}{\sigma^{2}+n(1-\delta)v_{k}^{2}} \geq 
(1-\delta)  \sum_{k=1}^N \frac{v_k^2\sigma^2/n}{\sigma^2/n + v_k^2 }\\
\overset{(\ref{def:v_j_and_s_j})}{=}&~
(1-\delta)\frac{\sigma^2}{n}\sum\limits_{j=1}^N\frac{\ell_j}{\ell_j+\kappa^{\star}\lambda_j^{-\frac{s}{2}}}\\
=&~
(1-\delta)\frac{\sigma^2}{n}\sum\limits_{j=1}^N {\ell_j}
=
(1-\delta)\mathcal{D}^{\star},
\end{aligned}
\end{equation}
where 
\iffalse
the inequality $\mathbf{(A)}$ follows from an elementary result that for any non-negative measurable function \( f : \calX_1 \times \calX_2 \to [0, \infty] \), we have
\[
\inf_{x_1 \in \calX_1} \int_{\calX_2} f(x_1, x_2) \, ~\mathsf{d}\rho_{\calX_2}(x_2) \geq \int_{\calX_2} \inf_{x_1 \in \calX_1} f(x_1, x_2) \, ~\mathsf{d}\rho_{\calX_2}(x_2);
\]
and 
\fi
the inequality $\mathbf{(A)}$ follows from the following arguments.
%in the third line follows from  $\mathbb{E}Z^{-1} \geq (\mathbb{E}Z)^{-1}$.
%Now we begin to prove (\ref{eqn:main_part_in_lower_bound}).

\iffalse
\textbf{Proof of $\mathbf{(A)}$:} 

1. Define \( g(x_2) = \inf_{x_1 \in \calX_1} f(x_1, x_2) \). Since \( f \) is a non-negative measurable function, \( g \) is also measurable.
2. For any \( x_1 \in \calX_1 \) and \( x_2 \in \calX_2 \), we have \( g(x_2) \leq f(x_1, x_2) \). Thus,
   \[
   \int_{\calX_2} g(x_2) \, \mathsf{d}\rho_{\calX_2}(x_2) \leq \int_{\calX_2} f(x_1, x_2) \, \mathsf{d}\rho_{\calX_2}(x_2).
   \]
3. Since the above inequality holds for any \( x_1 \in \calX_1 \), we get
   \[
   \int_{\calX_2} g(x_2) \, \mathsf{d}\rho_{\calX_2}(x_2) \leq \inf_{x_1 \in \calX_1} \int_{\calX_2} f(x_1, x_2) \, \mathsf{d}\rho_{\calX_2}(x_2).
   \]
4. Therefore, we have
   \[
   \int_{\calX_2} \inf_{x_1 \in \calX_1} f(x_1, x_2) \, \mathsf{d}\rho_{\calX_2}(x_2) = \int_{\calX_2} g(x_2) \, \mathsf{d}\rho_{\calX_2}(x_2) \leq \inf_{x_1 \in \calX_1} \int_{\calX_2} f(x_1, x_2) \, \mathsf{d}\rho_{\calX_2}(x_2).
   \]

\hfill \qedsymbol
\fi

For notation simplicity, let's denote $\hat{\theta}_k(\{y_i\}) = \hat{\theta}_k(\{x_i, y_i\})$ when $x_i$'s are given. For any $k, 1\leq k \leq N$, we have
\begin{align*}
    &~
    \inf_{\hat{\theta}_k} \int_{\mathbb{R}^{N}} \mathbb{E}_{\epsilon} \left[ (\hat{\theta}_k - \theta_k)^2 \mu(\theta^N) ~\mathsf{d}\theta^N \right]\\
    = &~
    \inf_{\hat{\theta}_k(\cdot)} \int_{\mathbb{R}^{N}} \int_{\mathbb{R}^{n}} \left(\hat{\theta}_k(\{y_i\}) - \theta_k\right)^2 \prod_{i=1}^{n} \mu_{\sigma^2}(\epsilon_i) \prod_{j=1}^N \mu_{s_j}(\theta_j)
    ~\mathsf{d}\epsilon_i  ~\mathsf{d}\theta_j\\
    \geq &~
    \int_{\mathbb{R}^{N-1}} \underbrace{\left[ \inf_{\hat{\theta}_k(\cdot)} \int_{\mathbb{R}} \int_{\mathbb{R}^{n}} \left(\hat{\theta}_k(\{y_i\}) - \theta_k\right)^2 \prod_{i=1}^{n} \mu_{\sigma^2}\left(y_i - \sum\nolimits_{j=1}^N \theta_j \phi_j(x_i)\right) \mu_{s_k}(\theta_k) \, ~\mathsf{d}y_i \, ~\mathsf{d}\theta_k \right]}_{\Delta} \\
    &~
    \cdot\prod_{j \neq k} \mu_{s_j}(\theta_j) ~\mathsf{d}\theta_1 \, \ldots \, ~\mathsf{d}\theta_{N}.
\end{align*}

Notice that
\begin{align*}
y_i | (\{x_i\}, \theta_1, \cdots, \theta_{k-1}, \theta_k, \theta_{k+1}, \cdots, \theta_N) &= \phi_k(x_i) \theta_k + \underbrace{\sum_{j \neq k} \phi_j(x_i) \theta_j}_{\Delta_i} + \epsilon_i,
\end{align*}
hence from Lemma \ref{lem:modified_lemma_3_4} we have $\Delta = \left(\frac{1}{s_k^2}+\frac{\sum_{i\leq n} \phi_k^2(x_i)}{\sigma^2}\right)^{-1} = \frac{s_k^2 \sigma^2}{\sigma^2 + \sum_{i\leq n} \phi_k^2(x_i) s_k^2}$. 

Therefore, we have
$$
\inf_{\hat{\theta}_{k}}\bbE_{\epsilon,\theta^{N}} \left[ (\hat{\theta}_k - \theta_k)^2 \mid (x_{1},\cdots,x_{n})\right] \geq  \frac{s_k^2 \sigma^2}{\sigma^2 + \sum_{i=1}^{n} \phi_k^2(x_i) s_k^2}.
$$

\begin{lemma}\label{lem:modified_lemma_3_4}
Let $c, c_1, \cdots, c_n, \Delta_1, \cdots, \Delta_n$ be $2n+1$ constants. Consider a statistical model with $n$ Gaussian observations:
$$
t_i = c_i a + \Delta_i + \epsilon_i, \quad i = 1, \cdots, n,
$$
where $a \in \mathbb{R}$, $\epsilon_i \overset{\mathcal{D}}{\sim}_{i.i.d.} N(0, \sigma^2)$. 
For an estimator $\hat{a}=\hat{a}(t_1, \cdots, t_n)$ of the parameter $a$, define its squared risk $\mathbb{E}\left[(\hat{a}-a)^2\right]$, as well as its Bayes risk with respect to the prior distribution $\calN\left(0, c^2\right)$:
\begin{equation*}
\begin{aligned}
\mathcal{R}^B(\hat{a})=\mathbb{E}\left[(\hat{a}(t_1, \cdots, t_n)-a)^2\right],
\end{aligned}
\end{equation*}
If we define the Bayes estimator as the minimizer of the Bayes risk among all estimators:
$$
\hat{a}^B=\arg \min _{\hat{a}} \mathcal{R}^B(\hat{a}),
$$
then we have
\begin{align*}
    \mathcal{R}^B\left(\hat{a}^B\right) = \left(\frac{1}{c^2}+\frac{\sum_{i=1}^{n} c_i^2}{\sigma^2}\right)^{-1}.
\end{align*}
\end{lemma}

\begin{proof}
    Denote $\mathbf{t}=(t_1, \cdots, t_n)$. The posterior distribution of $a$ is
    \begin{align*}
        a \mid \mathbf{t} \overset{\mathcal{D}}{\sim} \mathcal{N}\left(\mu^{\prime}, \sigma^{\prime 2}\right),
    \end{align*}
where
$
\sigma^{\prime 2}=\left(\frac{1}{c^2}+\frac{\sum_{i=1}^{n} c_i^2}{\sigma^2}\right)^{-1}
$
and
$
\mu^{\prime}= \frac{\sigma^{\prime 2}}{\sigma^2} \sum_{i=1}^n c_i\left(t_i-\Delta_i\right) 
$. Therefore, the Bayes estimator $\hat{a}^B$ is the mean of the posterior distribution:
$$
\hat{a}^B=\frac{\sigma^{\prime 2}}{\sigma^2} \sum_{i=1}^n c_i\left(t_i-\Delta_i\right),
$$
and the Bayes Risk is 
$$
\mathcal{R}^B\left(\hat{a}^B\right)=\sigma^{\prime 2} = \left(\frac{1}{c^2}+\frac{\sum_{i=1}^{n} c_i^2}{\sigma^2}\right)^{-1}.
$$
\end{proof}

\subsubsection{\bf Upper bound of $\mathbf{r}^{\star}$}\label{subsec_bound_of_resi}

The technique we use to bound the residual part $\mathbf{r}^{\star}$ in (\ref{eqn:decomp_lower_bound}) is quite standard, and we first recall some results in \cite{tsybakov2008introduction}.

\begin{lemma}[Restate (3.45) and (3.48) in \cite{tsybakov2008introduction}]\label{lem:restate_residual_control}
    We have
    \begin{equation}
        \begin{aligned}
            \mathbf{r}^{\star} &\leq 6 \lambda_1^s R\sqrt{\mathbb{P}_\mu(\mathbb{R}^{N} \backslash \Theta_N)}\\
            \mathbb{P}_\mu(\mathbb{R}^{N} \backslash \Theta_N) &\leq \exp{\left(-\frac{\delta^2}{8(1-\delta)^2}\frac{\sum\nolimits_{j=1}^Ns_j^2\lambda_j^{-s}}{\max_{1\leq j\leq N}s_j^2\lambda_j^{-s}}\right)},
        \end{aligned}
    \end{equation}
    where $\Theta_N$ is defined in  (\ref{def:theta_spca_lower_bound}), $\mu$ is the p.d.f. of $\theta^N \overset{\mathcal{D}}{\sim} \calN \left(\mathbf 0, \text{diag} \left(s_1^2, \ldots, s_N^2\right) \right)$, and $s_j$'s are defined in (\ref{def:v_j_and_s_j}).
\end{lemma}

Recall that we have $N = \sum_{k=0}^{q} N(d, k)$ for $q \in \{p, p+1\}$.
Since $\gamma \geq s > s/2$, from Lemma \ref{lemma:calcul_N} we have (i) $p \geq 1$ or (ii) $p=0$ and $q=p+1$. Hence, from the definition of $s_j$ and $\kappa^{\star}$, we have
\begin{align*}
    \sum_{j=1}^N s_j^2\lambda_j^{-s} &= 
    (1-\delta)\frac{\sigma^2}{n\kappa^{\star}}\sum_{j=1}^N \frac{\ell_j}{\lambda_j^{s/2}} = (1-\delta)R,
\end{align*}
and
\begin{align*}
    \max_{1\leq j\leq N}s_j^2\lambda_j^{-s}
    = (1-\delta){\sigma^2}  \max_{1\leq j\leq N} \frac{\ell_j}{n\lambda_j^{s/2}\kappa^{\star}}
    \overset{\text{Proposition } \ref{prop_calcul_max_ell}}{=}
    (1-\delta) O_d(d^{-\beta}),
\end{align*}
where $\beta = \min\{1, \gamma-s/2\}>0$ is a constant only depending on $\gamma$ and $s$.

Combining with Lemma \ref{lem:restate_residual_control}, we have
\begin{align*}
    \mathbf{r}^{\star} \leq 6 K_{\max}^s R 
    \exp{\left(-\frac{\delta^2 R}{16(1-\delta)^2} \Omega_d(d^{\beta})\right)}
    = o_d(\calD^{\star}).
\end{align*}

Finally, from  (\ref{eqn:decomp_lower_bound}) and (\ref{eqn:bound_of_I_star}), we have
\begin{align*}
    \inf_{\hat{f}} \sup_{f_{\star}\in \sqrt{R}[\calB]^{s}} \mathbb{E} \|\hat{f} - f_{\star}\|_{L^2}^2 \geq (1+o_d(1))\mathcal{D}^{\star} - \delta\calD^{\star},
\end{align*}
and the proof is completed by making $\delta$ tend to $0$.
\end{proof}

\end{document}